\documentclass[pdflatex,sn-mathphys-num]{sn-jnl}
\usepackage{graphicx}%
\usepackage{multirow}%
\usepackage{amsmath,amssymb,amsfonts}%
\usepackage{amsthm}%
\usepackage{mathrsfs}%
\usepackage[title]{appendix}%
\usepackage{xcolor}%
\usepackage{textcomp}%
\usepackage{manyfoot}%
\usepackage{booktabs}%
\usepackage{algorithm}%
\usepackage{algorithmic}
\usepackage{subfig}
\usepackage{listings}%
\usepackage{array}
\usepackage{hhline}

\usepackage{bm}

\newtheorem{remark}{Remark}[section]
\newtheorem{theorem}{Theorem}[section]  

\newtheorem{lemma}{Lemma}[section]

\newtheorem{corollary}[theorem]{Corollary}
\newtheorem{assumption}[theorem]{Assumption}
\theoremstyle{thmstyleone}%
\theoremstyle{thmstyletwo}%

\theoremstyle{thmstylethree}%

\begin{document}

\title[Article Title]{Group-Sparse Regularized CP Decomposition for Incomplete Tensors with Recovery Error Bounds}

\author[1]{Mingxin Chen}\email{cmx02@hnu.edu.cn}

\author*[1]{Minru Bai}\email{minru-bai@hnu.edu.cn}

\author[1]{Qiyuan Wei}\email{qywei@hnu.edu.cn}
\author[1]{Zihao Wang}\email{wzh97@hnu.edu.cn}

\affil*[1]{School of Mathematics, Hunan University, Lushan Road (S), Yuelu District, Changsha, 410082, China}


\abstract{Recovering high-order tensors from incomplete and noisy measurements is an ill-posed inverse problem arising in various applications, including chemometrics, signal processing, and machine learning. However, many existing approaches rely on low-rank surrogate functions, which may not provide an explicit and interpretable factor representation. Although tensor CANDECOMP /  \allowbreak PARAFAC(CP)  decomposition offers a structured and interpretable alternative, its application to incomplete tensors is challenging because the underlying CP rank is generally unknown and needs to be estimated from partially observed data. To address these issues, we propose a group-sparse regularization model for CP decomposition of incomplete tensors, which simultaneously achieves rank estimation, tensor factorization, and tensor completion. We further establish an error bound between the optimal solution of the model and its ground truth, which explicitly characterizes the recovery error in terms of the sampling size, the initial CP rank, the true CP rank, and other relevant problem parameters. To solve the proposed model, we develop a block coordinate proximal gradient algorithm with extrapolation and provide convergence analysis. Moreover, we design an adaptive tuning strategy for the $\ell_{2,0}$ regularization parameter. Numerical experiments on both synthetic and real-world data demonstrate the robustness and effectiveness of the proposed model and algorithm.}




\maketitle

\section{Introduction}\label{sec1}
The recovery of high-order tensors from partial and noisy observations is a classical problem with broad applications in chemometrics, signal processing, computer vision, and data mining. Low-rank tensor models provide an effective way to capture the underlying multilinear structure for data recovery. Among various tensor decomposition schemes, CANDECOMP/PARAFAC (CP) decomposition is a widely used low-rank representation that decomposes a high-order tensor as a sum of a finite number of rank-one tensors. The quality of the resulting representation is closely related to the CP rank.  Focusing on the third-order tensor which is the most widely applied high-order tensor form in engineering practice, this paper mainly investigates the rank CP decomposition problem for incomplete tensors. Given a noisy data tensor $\mathcal{Y} \in \mathbb{R}^{I \times J \times K} $, only its entries indexed by the set \(\Omega\) are observed, the task is to find minimal positive integer $R$, and matrices $\mathbf{A} \in \mathbb{R}^{I \times R} $, $\mathbf{B} \in \mathbb{R}^{ J \times R}$ and  $\mathbf{C} \in \mathbb{R}^{ K \times R}$ such that
\begin{equation}
	\mathcal{Y}_{ijk} = \sum_{r=1}^{R} a_{ir}b_{jr}c_{kr} + \mathcal{E}_{ijk}, \; (i,j,k) \in \Omega,
	\label{origin_problem}
\end{equation}
where  $\mathcal{E} $ represents the noisy tensor, $a_{ir}$, $b_{jr}$ and $c_{kr}$ represent the $(i, r)$th, $(j, r)$th, $(k, r)$th element of $\mathbf{A} $, $\mathbf{B} $ and  $\mathbf{C}$, respectively. This problem $(\ref{origin_problem})$ is closely related to the problem of component separation in chemometrics. When 
the matrices $\mathbf{A}  $ and $\mathbf{B} $ have column unit constraints, i.e.,
\begin{equation}
	\sum_{i=1}^{I}a_{ir}^2 = 1, \; \;\;\sum_{j=1}^{J}b_{jr}^2 = 1,  \; \;\; r=1,\dots,R,
	\nonumber
\end{equation}
then problem $(\ref{origin_problem})$ provides an
interpretable latent factor representation.
The matrices $\mathbf{A}  $ and $\mathbf{B} $ are called the qualitative profiles matrices, and $\mathbf{C}$ represent concentration profiles matrix \cite{wu2024}. The elements $a_{ir}$ and $b_{jr}$ usually represent the normalized intensities of the $r$-th component in the instrument
channels $i$ and $j$, respectively, and
$c_{kr}$ represents the relative concentration of the $r$-th component in the $k$-th sample. $R$ denotes the number of components. 

Classical CP decomposition methods were generally developed under the assumption that the tensor data are fully observed and that the CP rank is known. The most widely studied CP decomposition model is the unconstrained optimization problem of minimizing the least squares loss function, and several algorithms have been proposed to solve this optimization problem. Harshman \cite{harshman1970} proposed the parallel factor analysis  (PARAFAC) algorithm. Kolda et al. \cite{Kodal2018} developed randomized CP decomposition algorithms for handling large-scale dense and sparse tensors. Recently, constrained CP decomposition models have received increasing attention because they can incorporate prior information from practical applications. Common constraints include orthogonality, nonnegativity, and column unit constraints on the factor columns. For example, Wang et al. \cite{wang2024inexactly} developed a CP decomposition model with nonnegative and column unit constraints for the component separation problem in complex chemical systems, and proposed an accelerated inexact block coordinate descent algorithm to solve the proposed model. However, the true CP rank is often unknown and difficult to estimate in practice. To address this issue, Wang et al. \cite{wang2018} employed \(\ell_1\)-regularization on the scaling coefficients to obtain an approximate convex surrogate for the CP rank, established an \(\ell_1\) regularized CP decomposition model, and developed a proximal alternating minimization algorithm to solve the model.

On the other hand, tensor completion is to recover the missing entries from partially observed data.
Over the past decade, a wide range of tensor completion methods have been developed to tackle this issue by tensor rank surrogate methods. Convex relaxation is one of the most widely studied approaches to tensor completion, which aims to approximate the NP-hard tensor rank minimization problem by replacing the rank function with a nuclear norm based surrogate \cite{liu2012,mu2014}. Although nuclear?norm?based models are computationally tractable, the nuclear norm tends to over-shrink large singular values and may produce biased solutions. To mitigate the estimation bias, Bai et al. \cite{bai2016adaptive} developed an adaptive correction approach for higher-order tensor completion and demonstrated that a suitable estimator for the correction term could reduce the error bound of the corrected model, while the corrected model mainly addresses missing-data problems in noise-free settings. To approximate rank function more accurately, nonconvex rank surrogates have also been studied. Li et al. \cite{li2018low} established a nonconvex $\ell_p$-norm relaxation model for low Tucker rank tensor recovery, demonstrated its ability to recover tensors at lower sampling ratios than the convex nuclear-norm relaxation model, and employed the alternating direction method of multipliers (ADMM) to solve the resulting optimization problem. Recently, Zhao et al.\cite{zhao2022} proposed a nonconvex equivalent surrogate model with nonconvex equivalent surrogate functions of tensor rank and sparsity measure for robust tensor completion(RTC) problems, and proved that the optimal solution of the proposed nonconvex formulation shares identical global optimal solutions with RTC problem. The advantage of the above rank surrogate-based methods is that missing tensors can be recovered without prior knowledge of their rank. However, rank-surrogate-based tensor completion methods do not necessarily provide the explicit component-wise latent factor representation offered by CP decomposition, which is important in applications requiring interpretable factor structures. 

To address this limitation, some researchers have explored CP-decomposition-based methods for recovering missing entries. CP decomposition provides an interpretable latent factor representation, making it particularly suitable for applications such as chemical component separation. In the field of chemometrics, G. Tomasi et al. \cite{tomasi2005} proposed  an algorithm (INDAFAC) derived from PARAFAC based on the idea of the Levenberg-Marquardt(LM) method to recover incomplete data, which vectorizes the missing entries and employs the LM method to solve the resulting least-squares problem. More recently, Dong et al.\cite{dong2024} developed a novel second-order calibration algorithm for missing values (PARAFACM), which combines the idea of evaluating multivariate linear regression as a set of individual univariate linear regressions with the iterative least squares principle. Besides, Cai et al.\cite{cai2022nonconvex} studied low rank tensor completion problem for noisy third-order symmetric tensors based on the CP decomposition and derived the corresponding error bounds. Although these methods achieve satisfactory performance in computing the CP decomposition from incomplete tensors, their effectiveness relies on the assumption that the CP rank is known, and some of the corresponding algorithms lack theoretical convergence analysis. In particular, the error bounds established in \cite{Jain2017} and \cite{cai2022nonconvex} are derived under the assumption that the CP rank is known. Therefore, these results do not cover the case where the true CP rank is unknown and the prescribed initial CP rank may exceed it.
To address the problem of rank estimation, Zhao et al.\cite{Zhao2015} and Jiang et al.\cite{jiang2022l} incorporated rank estimation into CP decomposition in the presence of missing entries. Zhao et al.\cite{Zhao2015} formulated CP factorization within a hierarchical probabilistic framework and adopted a fully Bayesian approach by imposing sparsity-inducing priors on the latent factors together with appropriate hyperpriors over all hyperparameters, thereby enabling automatic rank determination. However, the method in \cite{Zhao2015} may become ineffective when the proportion of missing entries exceeds 0.5. Jiang et al.\cite{jiang2022l} applied $\ell_1$ regularization to the scaling coefficients, which serves as an approximate surrogate for the CP rank. Very recently, Wang et al. \cite{wang2025} 
proved that the CP rank can be exactly 
estimated by minimizing the group sparsity of any one of the factor matrices under the unit length constraints on the columns of the other factor matrices, and proposed a CP
decomposition model with group sparsity regularization, which integrates rank estimation and tensor decomposition as an optimization problem. Wang et al.\cite{wang2018} also discuss the probabilistic consistency of the sparsity result and suggest a way to choose the $\ell_1 $-regularization parameter for practical computation.
Note that these methods in \cite{wang2025} and \cite{wang2018}  are designed for tensors without missing data. In the presence of missing entries, the choice of the regularization parameter becomes more sensitive because the available data provide limited information about the underlying CP structure. Therefore, developing a theoretically derived adaptive parameter selection strategy for CP decomposition with incomplete tensors is particularly important.

These observations raise one important open questions. Can we extend the group-sparse CP rank estimation framework \cite{wang2025} to cope with incomplete and noisy tensor data, such that rank estimation, factorization and tensor completion can be performed jointly?

To tackle these challenges, we develop a group-sparse regularization model for CP decomposition of incomplete tensors, which can achieve rank estimation, tensor decomposition, and missing-entry completion simultaneously. Then, we establish the error bound between the optimal solution of the proposed model and its true solution under general noise assumptions and further derive an explicit error bound under the additive Gaussian noise setting  without the prior of the true CP rank. Our error bound additionally accounts for the initial CP rank and depends jointly on the initial CP rank, the true CP rank, the number of observed entries, and other relevant problem parameters. To solve this model, we proposed a block-coordinate proximal gradient descent algorithm, enhanced with an extrapolation step after each block update. We prove that every limit point of the resulting sequence is a stationary point of the proposed model. In addition, we design an adaptive tuning strategy for the nonconvex $\ell_{2,0}$ regularization parameter. Extensive numerical experiments on synthetic low-rank tensors and real chemometric datasets verify that our method is capable of handling over-specified initial component numbers and low sampling rates, and demonstrate the practical robustness of the proposed model and algorithm.

\subsection{Contributions}
The main contributions of this paper are summarized below, with clear distinction from the prior group-sparse CP work for fully observed tensors \cite{wang2025}:

$\bullet$ We propose a CP decomposition model with group sparsity regularization for incomplete tensors, which simultaneously achieves rank estimation, tensor factorization, and tensor completion from partially sampled entries. This extends the existing  group-sparse CP theory from complete tensor to incomplete observation settings.

$\bullet$ We establish an error bound between tan optimal solution of the proposed model and the ground truth under general noise assumptions under general noise assumptions, and give explicit expressions for Gaussian additive noise. The resulting error bound characterizes the dependence of the recovery error on the initial CP rank, the true CP rank, the number of observed entries, and other relevant problem parameters. Furthermore, we design an adaptive tuning strategy for the nonconvex $\ell_{2,0}$ regularization parameter, enabling its adaptive selection under incomplete observations. These theoretical results are new and specific to the incomplete tensor setup.

$\bullet$ We propose a block-coordinate proximal gradient descent algorithm with verification of the descent property and extrapolation to solve the proposed model and provide convergence analysis. This algorithmic convergence result serves as a supporting technical tool rather than the primary theoretical innovation. Numerical experiments on real-world chemical data and randomly generated low rank tensors demonstrate that our proposed method can overcome the problem of overestimated initial number of component across different sampling rates. The experimental results further validate the effectiveness of the proposed model and algorithm.

\subsection{Organization}
The remainder of the paper is structured in the following way. Some notation and definitions are introduced in Section 2. In Section 3, we propose the regularization CP decomposition model for incomplete and noisy tensors. We establish an error bound between the optimal solution of the model and its ground truth in Section 4.  We develop a block-coordinate proximal gradient descent algorithm with extrapolation to solve this optimization problem and analyze the convergence of the proposed algorithm in Section 5. In Section 6, we conduct extensive numerical experiments on randomly generated low rank tensors and real-world data to demonstrate the robustness and effectiveness of the proposed model and algorithms. We summarize this paper in Section 7.

\section{Preliminaries}
Before proceeding, we first present some definitions and notations. Tensors are denoted by calligraphic letters, e.g., $\mathcal{X}$.
Matrices are denoted by boldface capital letters,e.g., $\mathbf{X} $. Vectors are denoted by bold lowercase letters, e.g.,$\textbf{x}$, and scalars are denoted by ordinary letters, e.g., $x$. The field of real numbers is denoted as $\mathbb{R}$. Sets are denoted by hollow letters,e.g., $\mathbb{X}$. For a third-order tensor $\mathcal{X} \in \mathbb{R}^{I \times J \times K}$, we denote its $(i,j,k)$ th entry and mode-$i$ unfolding by $\textbf{X}_{ijk}$ and $\mathbf{X}_{(i)} $, respectively. The $i$ th column of the matrix $\textbf{X}$ is denoted as $\textbf{x}_i$. The transpose of matrix  $\mathbf{X}$ is denoted by $\mathbf{X}^T$. We use the symbols $*$, $ \circ$,  $\bigotimes$ and $ \bigodot $ to represent the Hadamard product, outer product, Kronecker product and  Khatri-Rao product, respectively. The Frobenius norm of tensor $\mathcal{X}$ and matrix $\mathbf{X}$ is denoted by $\Vert \mathcal{X} \Vert$ and  $\Vert \mathbf{X} \Vert$,  respectively. The $Euclidean$ norm of vector $\mathbf{x} $ is denoted by $\Vert \mathbf{x} \Vert$. The infinity norms of tensor $\mathcal{X}$ and matrix $\mathbf{X}$, denoted by $\Vert \mathcal{X} \Vert_{\infty}$ and $\Vert \mathbf{X} \Vert_{\infty}$ are defined as $\Vert \mathcal{X} \Vert_{\infty} = \max_{i,j,k}|\textbf{X}_{ijk}|$ and $\Vert \mathbf{X} \Vert_{\infty} = \max_{i,j}|\textbf{X}_{ij}|$, respectively. Considering a group structure for $\mathbf{X}$ based on prior knowledge $\mathbf{X}=(\mathbf{x}_1,\dots,\mathbf{x}_m) \in \mathbb{R}^{n \times m}$, the $\ell_{2,0}$ norm of $\mathbf{X} $ is defined as $\Vert \mathbf{X} \Vert_{2,0}:=\sum_{i=1}^{m} \Vert \mathbf{x}_i \Vert ^0 $, where $0^0 = 0$. In this paper, the column of the matrix $\mathbf{X}$ are defined groups. \\

	$\textsc{Definition 2.1} $ (\textit{Moreau envelopes and proximal mappings \cite{Bolte2013}  }) \textit{For a proper, lower semicontinuous function }$f$ : $\mathbb{R}^n \to \mathbb{R}$ \textit{and parameter value} $\lambda > 0$, \textit{the Moreau envelope function} $e^{f}_{\lambda}$ \textit{and the proximal mapping} $prox^{f}_{\lambda}$ \textit{are defined as} 
\begin{equation*}
	\begin{split}
		& e^{f}_{\lambda}(x):=\mathop{\inf} \limits_{w}\{f(x)+ \frac{1}{2 \lambda} \Vert w -x \Vert^2 \},\\
		& 
		prox^{f}_{\lambda}(x):=\mathop{\text{argmin}} \limits_{w}\{f(x)+ \frac{1}{2 \lambda} \Vert w -x \Vert^2 \}. \\
	\end{split}  
\end{equation*}

$\textsc{Definition 2.2}  $ (\textit{Kullback-Leibler divergence} \cite{cover2012elements}  ) \textit{Let the random variable $y $  have probability density functions $ p_{x_1}(y) $ and $ p_{x_2}(y)  $
	under parameters $x_1 $ and $x_2  $, respectively. The Kullback?Leibler divergence from $ p_{x_1}(y) $ and $ p_{x_2}(y)  $ is defined as} 
\begin{equation}
	D \left(p_{x_1}(y) \Vert  p_{x_2}(y)  \right) = \mathbb{E}_{p_{x_1}} \left[  \ln \frac{p_{x_1}(y)}{p_{x_2}(y)}  \right],
\end{equation}
where $ \mathbb{E}_{p_{x_1}} \left[ \cdot \right]  $ denote expectation with respect to the probability density function $p_{x_1} . $

$\textsc{Definition 2.3} $ (\textit{Hellinger affinity} \cite{cover2012elements} ) \textit{Let the random variable $y $  have probability density functions $ p_{x_1}(y) $ and $ p_{x_2}(y)  $
	under parameters $x_1 $ and $x_2  $, respectively.  The Hellinger affinity from $ p_{x_1}(y) $ and $ p_{x_2}(y) $  is defined as }
\begin{equation}
	H \left(p_{x_1}(y) \Vert  p_{x_2}(y)  \right) = \mathbb{E}_{p_{x_1}} \left[  
	\sqrt{\frac{p_{x_2}(y)}{p_{x_1}(y)}}
	\right].
\end{equation}
The  Hellinger affinity satisfies symmetry, i.e.
\begin{equation}
	H \left(p_{x_1}(y) \Vert  p_{x_2}(y)  \right) = \mathbb{E}_{p_{x_2}} \left[  
	\sqrt{\frac{p_{x_1}(y)}{p_{x_2}(y)}}
	\right],
\end{equation}
where $ \mathbb{E}_{p_{x_1}} \left[ \cdot \right]  $ and $ \mathbb{E}_{p_{x_2}} \left[ \cdot \right]  $ denote expectation with respect to the probability density function $p_{x_1}(y)  $ and  $ p_{x_2}(y)$, respectively.

\section{A group sparse-based CP decomposition model for incomplete tensors}
In this section, we will establish a model for problem $(\ref{origin_problem})$ by investigating CP decomposition for incomplete tensors when CP rank is unknown. By leveraging the equivalence between the CP rank and the minimization of the group sparsity on any of the factor matrices under the unit length constraints on the columns of the others \cite[Theorem 3.1]{wang2025}, we proposed a group sparse-based CP decomposition model for incomplete tensors for problem $(\ref{origin_problem})$, formulated as
\begin{equation}
	\begin{array}{cc}		 \mathop{\text{min}}\limits_{\mathbf{A},\mathbf{B},\mathbf{C}} \frac{1}{2}|| {\mathcal{Y}}_{\Omega} -([\![\mathbf{A},\mathbf{B},\mathbf{C}] \! ])_{\Omega} ||^2_F + \lambda \Vert  \mathbf{C} \Vert_{2,0}   \\
		s.t. \; \Vert \mathbf{a}_{j}\Vert = 1,\, \Vert \mathbf{b}_{j} \Vert = 1, \, j=1,2,\dots,R,    
	\end{array}
	\label{4}
\end{equation}
where $ [\! [\mathbf{A},\mathbf{B},\mathbf{C}] \! ]:=\sum_{r=1}^{R}\mathbf{a}_r \circ \mathbf{b}_r \circ \mathbf{c}_r $, $\Vert \mathbf{C} \Vert_{2,0}:=\sum_{i=1}^{R} \Vert \mathbf{c}_i \Vert ^0$, $\mathbf{a}_r \in \mathbb{R}^{I} $, $\mathbf{b}_r \in \mathbb{R}^{J} $ and $\mathbf{c}_r \in \mathbb{R}^{K} $ represent the $r$-th column of the factor matrix $\mathbf{A} $, $\mathbf{B} $ and $\mathbf{C} $, respectively, 
$\Omega  $ represents the index set of observed entries, and the tensor ${\mathcal{X}_{\Omega}}$ is partially observed over the index set $\Omega$. The positive integer \(R\) is an initial rank that is assumed to be no smaller than the rank of the underlying noiseless tensor corresponding to \(\mathcal{Y}\). 
The scalar $\lambda>0 $ is the regularization parameter.
An equivalent form of the model $(\ref{4})$ is given by 
\begin{equation}
	\begin{array}{cc}
		\mathop{\text{min}}\limits_{\mathbf{A},\mathbf{B},\mathbf{C}} \frac{1}{2}||\mathcal{M}*(\mathcal{Y}-[ \! [\mathbf{A},\mathbf{B},\mathbf{C}] \! ]) ||^2_F + \lambda \Vert \mathbf{C}\Vert_{2,0}   \\
		s.t. \; \Vert \mathbf{a}_{j} \Vert = 1,\; \Vert \mathbf{b}_{j} \Vert = 1,j=1,2,\dots,R,     			
	\end{array}
	\label{CP MV}
\end{equation}
where the entries of of $\mathcal{M}$ are defined as follows:
\begin{equation}
	\begin{array}{ccc}
		\mathcal{M}_{ijk}  = \begin{cases}  1  \quad & (i, j, k) 	 \in \Omega, \\ 
			0  \quad & (i, j, k) 	 \notin \Omega.  \\ 
		\end{cases} 
	\end{array}
	\label{M}
	\nonumber
\end{equation}
Our model $(\ref{4})$ can perform rank estimation, tensor decomposition, and tensor completion simultaneously. The proposed model yields three interpretable factor matrices $\mathbf{A},\mathbf{B},\mathbf{C} $ from the incomplete tensor $ {\mathcal{X}} $ and the complete tensor is then recovered by fusing these factors, thereby achieving tensor completion simultaneously. By making use of the indicator function, we can rewrite the model $(\ref{CP MV})$ into unconstrained optimization problem as follows:
\begin{equation}
	\begin{array}{cc}
		\mathop{\text{min}}\limits_{\mathbf{A},\mathbf{B},\mathbf{C}} G(\mathbf{A},\mathbf{B},\mathbf{C}) \! = \!
		\frac{1}{2} \Vert \mathcal{M}*(\mathcal{Y}-[ \![\mathbf{A},\mathbf{B},\mathbf{C}] \! ]) \Vert ^2_F \! + \! \lambda \Vert \mathbf{C} \Vert_{2,0} \!+\! \delta_{\mathbb{A}}(\mathbf{A}) \!+ \! \delta_{\mathbb{B}}(\mathbf{B}),
	\end{array}
	\label{spe_unconstrain}
\end{equation}
where $\delta_{\mathbb{S}}$ is the indicator function of the set $\mathbb{S}$, and $\mathbb{A} := \{\mathbf{A} \in \mathbb{R}^{I \times R} \, | \,\Vert \mathbf{a}_r \Vert = 1 \}$, $\mathbb{B} := \{\mathbf{B} \in \mathbb{R}^{J \times R} \, \allowbreak
|  \, \Vert \mathbf{b}_r \Vert = 1 \}. $   
\begin{remark}
	If $ \Omega $ is a fully sampled set, model $(\ref{spe_unconstrain}) $ reduces to the model  from \cite{wang2025} when the regularization parameter $\kappa = 0$.
\end{remark}

	\section{Error bound} 
In this section, we establish an error bound between an optimal solution of the generalized model presented below and the corresponding ground-truth tensor when the true CP rank is unknown.

To this end, let $\mathcal{Y}\in\mathbb{R}^{n_1\times\cdots\times n_N}$ be an $N$th-order data tensor, of which only the entries indexed by $\Omega$ are observed. For the factor matrices
\[
\mathbf{A}^{(i)}
:=
\bigl(\mathbf{a}_1^{(i)},\ldots,\mathbf{a}_R^{(i)}\bigr)
\in\mathbb{R}^{n_i\times R},
\qquad i=1,\ldots,N,
\]
define
\[
[\![\mathbf{A}^{(1)},\ldots,\mathbf{A}^{(N)}]\!]
:=
\sum_{r=1}^{R}
\mathbf{a}_r^{(1)}\circ\mathbf{a}_r^{(2)}
\circ\cdots\circ\mathbf{a}_r^{(N)}.
\]
With this notation, the proposed model in $(\ref{4})$ can be naturally extended to the $N$th-order case as follows:
\begin{equation}
	\begin{array}{cc}
		\mathop{\text{min}}\limits_{\mathbf{A}^{(1)},\dots,\mathbf{A}^{(N)}} \frac{1}{2}|| {\mathcal{Y}}_{\Omega} -([\![\mathbf{A}^{(1)},\dots,\mathbf{A}^{(N)}] \! ])_{\Omega} ||^2_F + \lambda \Vert  \mathbf{A}^{(N)} \Vert_{2,0}   \\
		s.t. \; \Vert \mathbf{a}_{j}^{(i)}\Vert = 1, \; i=1,\dots,N-1, \; j=1,2,\dots,R.
	\end{array}
	\label{extend_model}
\end{equation}

\subsection{Assumptions and the MLE model}
Assume tensor $ \mathcal{X}^{*} \in $ 
$\mathbb{R}^{n_1 \times \cdots n_N } $ is a true tensor whose CP rank is $\bar{R} $. We present three basic assumptions.

\textbf{(i)} $ \mathcal{X}^{*} $ admits the following CP decomposition:
\begin{equation}\label{X*}
	\begin{array}{clc}
		\mathcal{X}^*=[\![\mathbf{A}^{(1)^*}, \cdots, \mathbf{A}^{(N)^*}]\!]=\sum_{r=1}^{R_1} \mathbf{a}^{(1)^*}_r \circ \cdots \circ \mathbf{a}^{(N)^*}_r,
	\end{array}
\end{equation}
where $\mathbf{A}^{(i)^*}\in \mathbb{R}^{n_i \times {R_1}}(i=1,\cdots,N-1)$ are the column-normalized matrix. $\mathbf{A}^{(N)^*}\in \mathbb{R}^{n_N \times {R_1}}$ satisfies element-wise boundedness, i.e., there exists a positive constant $L$, so that $ \Vert \mathbf{A}^{(N)^*} \Vert_{\infty} \leq L $. According to $(\ref{X*})$, we have $  \Vert  \mathcal{X}^*  \Vert_{\infty}  \leq L $. 

\textbf{(ii)} Suppose that we observe a noisy version of $ \mathcal{X}^{*}  $ in some random subset of the entries $\Omega \subset [n_1] \times \cdots [n_N] $, which is generated by using an independent Bernoulli model with parameter $ \rho \in (0,1] $. The specific procedure is as follows: $\mathbf{(1)}$ generate $n_1 \cdots n_N $ i.i.d. Bernoulli random variables
$ m_{i_1,\cdots,i_N} $ with $Pr(m_{i_1,\cdots,i_N}=1) = \rho $, $ \forall i_1,\cdots,i_N $; $\mathbf{(2)}$ obtain the set $ \Omega = \{ (i_1,\cdots,i_N) |  \allowbreak m_{i_1,\cdots,i_N}=1  \}$. Under the additive Gaussian noise model and conditioned on $\Omega $, we observe the following noisy measurements at the positions of $\Omega $
\begin{equation}
	\begin{array}{clc}
		\mathcal{Y}_{i_1, \cdots, i_N}=\mathcal{X}_{i_1, \cdots, i_N}^*+\mathcal{E}_{i_1, \cdots, i_N},\;(i_1, \cdots, i_N) \in \Omega,
	\end{array}
\end{equation}
where $ \mathcal{E}_{i_1, \cdots, i_N} $ are i.i.d noise terms and $\mathcal{Y}_{i_1, \cdots, i_N}  $ are observed entries.

\textbf{(iii)} Conditional on the sampling set \(\Omega\), the observations \(\mathcal{Y}_{i_1,\ldots,i_N}\) are mutually independent, with their distributions parameterized by the corresponding entries of the true tensor \(\mathcal{X}^*\). Consequently, the joint probability density function of \(\mathcal{Y}_{\Omega}\) is given by
\begin{equation}
	p_{\mathcal{X}^{*}_{\Omega}(\mathcal{Y}_{\Omega})} := \prod_{(i_1, \cdots, i_N)\in \Omega}
	p_{\mathcal{X}^{*}_{i_1, \cdots, i_N}}\!\left(\mathcal{Y}_{i_1, \cdots, i_N}\right),
	\label{y_defin}
\end{equation}
where $p_{\mathcal{X}^{*}_{i_1, \cdots, i_N}}\!\left(\mathcal{Y}_{i_1, \cdots, i_N}\right) $ is the probability density function of observation $\mathcal{Y}_{i_1, \cdots, i_N} $ which depends on the probability density function of the noise and is parametrized by $\mathcal{X}^{*}_{i_1, \cdots, i_N} $.

Let  $\beta \geq 1$ and $ R \geq \bar{R}$,
define $n_{\max} = \max \{ n_1,\dots, $
$n_N\}
, \gamma =
2^{\left\lceil\log_2\left(n_{\max }\right)^\beta\right\rceil} $, where
$ {\left\lceil \cdot \right\rceil} $
represents the round-up function.

Now, we construct the candidate estimator set of $ \mathcal{X}^{*}_{\Omega} $, i.e.,
\begin{equation}\label{X_candi}
	\begin{array}{rcl}
		\mathbb{X}=\big\{\mathcal{X}= [\![\mathbf{A}^{(1)},\dots,\mathbf{A}^{(N)}] \! ] \mid \mathbf{A}^{(i)} \in \mathbb{A}_i, \; i=1,\cdots,N, \; \Vert \mathcal{X}\Vert_{\infty} \leq 2L   \big\},
	\end{array}
\end{equation}   
where $\mathbb{A}_{i} (i=1,\dots,N-1  )$ is the set of matrices $\mathbf{A}^{(i)}\in\mathbb{R}^{n_i\times R}$, whose columns have unit Euclidean norm and whose entries are discretized into $\gamma$ equally spaced levels over the interval $[-1,1]$, and $\mathbb{A}_{N}$ is the set of matrices $\mathbf{A}^{(N)}  \in \mathbb{R}^{n_N \times R}$, whose column vectors have zero Euclidean norm or are discretized into $ \gamma  $ equally spaced levels within the interval $[-L,L] $.

By Assumption \textbf{(i)}$-$\textbf{(iii)}, the extended model $ ((\ref{extend_model})$ can be rewritten as follows
\begin{equation}
	\label{our_model}
	\min_{\mathcal{X} = [\![\mathbf{A}^{(1)},\dots,\mathbf{A}^{(N)}] \! ] \in \mathbb{X}} \frac{1}{2 \sigma^2} \Vert   \mathcal{Y}_{\Omega}  -  \mathcal{X}_{\Omega}    \Vert^2 - \ln \frac{1}{(2 \pi \sigma^2)^{\Omega/2} } + \kappa \Vert \mathbf{A}^{(N)} \Vert_{2,0},
\end{equation}
where $\kappa = \lambda/ \sigma^2 $. If $ \ln p_{\mathcal{X}_{\Omega}}(\mathcal{Y}_{\Omega}) = - (\frac{1}{2 \sigma^2} \Vert   \mathcal{Y}_{\Omega}  -  \mathcal{X}_{\Omega}    \Vert^2 - \ln \frac{1}{(2 \pi \sigma^2)^{\Omega/2} }),$ the model $(\ref{our_model})$ can be rewritten as 
\begin{equation}
	\label{MLE_model}
	\min_{\mathcal{X} = [\![\mathbf{A}^{(1)},\dots,\mathbf{A}^{(N)}] \! ] \in \mathbb{X}}  - \ln p_{\mathcal{X}_{\Omega}}(\mathcal{Y}_{\Omega})  + \kappa \Vert \mathbf{A}^{(N)} \Vert_{2,0},
\end{equation}
where $p_{\mathcal{X}_{\Omega}}(\mathcal{Y}_{\Omega})$ denotes probability density function. Model $(\ref{MLE_model})$ is a maximum likelihood estimation (MLE) model that generalizes model $(\ref{our_model})$. In particular, when the probability density function $p_{\mathcal{X}_{\Omega}}(\mathcal{Y}_{\Omega})$ corresponds to additive Gaussian noise, model $(\ref{MLE_model})$ reduces to model (\ref{our_model}). Next, we establish an error bound between the optimal solution of model $(\ref{MLE_model})$ and its ground truth under general noise. Then, we obtain an error bound of model $(\ref{extend_model})$ by assuming the observation $\mathcal{Y}_{\Omega} $  are corrupted by additive Gaussian noise.

\begin{remark} 
	The value of $\beta $ should not be chosen excessively large. If $\beta $ is chosen sufficiently large, the discreteness assumption on the factor matrices in $\mathbb{X}  $ can be approximately eliminated. This is because continuous optimization algorithms executed on a computer are inherently limited by finite precision arithmetic.
\end{remark}

\subsection{Analysis of error bound}
In this subsection, we derive the theoretical error bound. Our proof strategy is primarily inspired by those in \cite{soni2016}, but differs from them in several respects.

Define the "good" set $\mathcal{G}_{\eta,\rho} := \mathcal{G}_{\eta,\rho}(\mathcal{X}^{*},\mathcal{X}) $  as the subset of all possible sampling sets $\Omega \in [n_1]\times[n_2]\times\cdots\times[n_N] $ such that, for every $\mathcal{X}\in\mathbb{X} $,
$$
D \left( p_{\mathcal{X}^*_{\Omega}} || p_{\mathcal{X}_{\Omega}}  \right) \leq \frac{3\rho}{2} D \left( p_{\mathcal{X}^*} || p_{\mathcal{X}}  \right) +
\frac{4 E_{KL}}{3} \left[ \ln \left( \frac{1}{\sigma}   \right) + \operatorname{pen}(\mathcal{X}) \ln2 \right]
$$
and
$$
\!-2 \! \ln  H\left(p_{\mathcal{X}^*_{\Omega}},p_{\mathcal{X}_{\Omega}}\right)
\!\geq \! \frac{\rho}{2} \left( \!-2\! \ln  H\left( p_{\mathcal{X}^*},p_{\mathcal{X}}\right)        -\!  \frac{4 E_{H}}{3} \!\! \left[ \ln \left( \! \frac{1}{\sigma} \!  \right) \! +\! \operatorname{pen}(\mathcal{X}) \ln2    \right] \right)\!, 
$$
where $\rho = l/(n_1 \dots n_{N}) $ denotes the ratio of measured entries to total entries, $ l \leq n_1 \dots n_{N} $ is the expected total number of measurements, and $ \eta \in (0,1) $ is an arbitrary constant. $D ( p_{\mathcal{X}^*_{\Omega}} || p_{\mathcal{X}_{\Omega}}  )$  represents the Kullback$-$Leibler divergence from $p_{\mathcal{X}^*_{\Omega}} $ and $ p_{\mathcal{X}_{\Omega}} $, $H ( p_{\mathcal{X}^*},p_{\mathcal{X}}) $ represents the Hellinger affinity from $p_{\mathcal{X}^*_{}} $ and $ p_{\mathcal{X}_{}} $.

The following Lemma $\ref{good_set}$ and Lemma $\ref{lemma_43}$ are the Nth-order generalizations of Lemma 11   and  Lemma 8  in \cite{soni2016}, respectively. Lemma $\ref{lemma_43}$ is proved by applying  Lemma $\ref{good_set}$. Their proofs follow the same arguments as those in \cite{soni2016}, with S replaced by $\Omega$ and $n_1n_2$ replaced by $n_1\cdots n_N$. For brevity, we omit the proofs.

\begin{lemma}
	\label{good_set}
	Let $\mathbb{X} $ be a countable collection of candidate estimates $\mathcal{X} $ of  $\mathcal{X}^{*} \in \mathbb{R}^{n_1 \times \cdots \times n_N} $, with corresponding penalties  $\operatorname{pen} \left( \mathcal{X} \right)  $ satisfying the Kraft inequality $(\ref{KM inequality}) $. Let $ \Omega $ be a random subset of 
	$[n_1] \times \cdots \times [n_{N}]$ generated according to the Erd\H{o}s-R\'enyi model. For any fixed
	$ \eta,\rho \in (0,1) $, we have $\mathbb{P}\left[ \Omega \notin \mathcal{G}_{\eta,\rho}  \right]  \leq 2 \eta .$
\end{lemma}

\begin{lemma}
	\label{KM_ineq}
	Let $\mathbb{X} $ be a countable collection of candidate estimates $\mathcal{X}$  of $\mathcal{X}^{*} \in \mathbb{R}^{n_1 \times \cdots \times n_N} $, each with a corresponding penalty term $pen(\mathcal{X}) \geq 1 $, so that the collection of penalties satisfies
	the summability condition:
	\begin{equation}\label{KM inequality}
		\sum_{\mathcal{X} \in \mathbb{X}} 2^{-\operatorname{pen}(\mathcal{X})} \leq 1.
	\end{equation}
	Given a positive integer $l$ satisfying $ k \leq  l \leq n_1  \cdots n_N$ and set $\rho = \frac{l}{n_1  \cdots n_N} $. Subsequently, generate $ n_1  \cdots n_N$ independent Bernoulli$(\rho)$ random variables $ \{ \Omega_{i_1,\cdots,i_N}\}$ so that each  
	$(i_1,i_2,\dots,i_{N}) 
	\in [n_1]\times[n_2] \times \cdots \times [n_{N}]$ is included in $\Omega$ independently with probability $\eta $.
	An index $(i_1,\cdots,i_N)$ is included in the sampling set $\Omega$ if $ \Omega_{i_1,\cdots,i_N} =1$, and is excluded otherwise. For the given sampling set $ \Omega$,
	the observations $ \mathcal{Y}_{\Omega} $ are mutually independent, with joint distribution $\mathcal{Y}_{\Omega} \sim p_{\mathcal{X}^{*}_{\Omega}} =\prod\limits_{(i_1,\cdots,i_N) \in \Omega} p_{\mathcal{X}^{*}_{(i_1,\cdots,i_N)}}  $. Let $E_{KL}$ be an upper bound for the maximum KL-divergence 
	\begin{equation}
		E_{KL} \geq  \max\limits_{\substack{\mathcal{X}\in\mathbb{X}}} \max\limits_{\substack (i_1,\cdots,i_N)} D\!\left( p_{X^{*}_{i_1,\cdots,i_N}} \,\|\, p_{X_{i_1,\cdots,i_N}} \right),
		\label{EKL_BOUND}
	\end{equation}
	and choose $ \zeta $ such that $ \zeta \geq (1+ \frac{2 E_{KL}}{3}) 2 \ln2 $. Then,  the following complexity-penalized maximum likelihood estimator
	\begin{equation}
		\widehat{\mathcal{X} }^ {\zeta}(\Omega,\mathcal{Y}_{\Omega}) = \arg\min_{\mathcal{X}\in \mathbb{X}} \{ -\ln p_{\mathcal{X}_{\Omega}}(\mathcal{Y}_{\Omega}) + \zeta \operatorname{pen}(\mathcal{X})   \}
		\label{X_armin}
	\end{equation}
	satisfies the following error bound:
	\begin{equation*}
		\begin{array}{rl}
			\frac{\mathbb{E}_{\Omega,\mathcal{Y}_{\Omega}}[-2 \ln(H(p_{\widehat{\mathcal{X} }^ {\zeta}},p_{\mathcal{X}^{*}} ))]}{n_1  \cdots n_N} 
			\leq  3 \min \limits_{\substack{\mathcal{X}\in\mathbb{X} }}  \left\{   \frac{D( p_{\mathcal{X}^{*}} \,\|\, p_{\mathcal{X}} ) }{n_1  \cdots n_N} + ( \zeta +\frac{4 E_{KL} \ln2 }{3} )\frac{\operatorname{pen} (\mathcal{X})}{l}   \right\}
			+ \frac{8 E_{KL} \ln l}{l},  
		\end{array}
	\end{equation*}
	\label{lemma_43}
\end{lemma}  
where the expectation $ \mathbb{E}_{\Omega,\mathcal{Y}_{\Omega}}[\cdot] $ is with respect to the joint
distribution of $\Omega $ and $ \mathcal{Y}_{\Omega}$.

Applying Lemma $\ref{lemma_43}$, we establish the error bound between the optimal solution of the extended model $(\ref{extend_model})$ and its ground truth under general noise. 
\begin{theorem}
	\label{Bernoulli sampling}
	Let $ \Omega $ be sampled according to the independent Bernoulli sampling model with parameter $ \rho = \frac{l}{n_1  \cdots n_N}$, and let 
	$\mathcal{Y}_{\Omega}$ denote the corresponding observations defined in $(\ref{y_defin}) $. Denote by $ E_{KL} $ any upper bound satisfying $(\ref{EKL_BOUND})$ on the maximum Kullback?Leibler divergence between \(p_{\mathcal{X}^*_{i_1,\ldots,i_N}}\) and \(p_{\mathcal{X}_{i_1,\ldots,i_N}}\) for $\mathcal{X}\in\mathbb{X} $, where $\mathcal{X}$ is defined in $(\ref{X_armin})  $. If $\kappa \geq
	2\widehat{N}\left(1+\frac{2E_{\mathrm{KL}}}{3}\right)
	(\beta+1) \allowbreak \ln(n_{\max}) $, then the regularized constrained maximum likelihood estimator 
	$$
	\widehat{\mathcal{X}}_{}^{\kappa}=\underset{\mathcal{X}= [\![\mathbf{A}^{(1)}, \cdots, \mathbf{A}^{(N)}]\!] \in \mathbb{X}}{\arg\min} \{  -\ln p_{\mathcal{X}_{\Omega}}\left(\mathcal{Y}_{\Omega} \right)+\kappa \Vert \mathbf{A}^{(N)}\Vert _{2,0}    \} 
	$$
	satisfies the following error bound:
	\begin{equation*}
		\begin{array}{rl}
			\frac{\mathbb{E}_{\Omega,\mathcal{Y}_{\Omega}}[-2 In(H(p_{\widehat{\mathcal{X} }^ {\kappa}},p_{\mathcal{X}^{*}} ))]}{n_1  \cdots n_N} 
			&\leq  3 \min \limits_ { \substack{\mathcal{X}\in\mathbb{X} }}  \Big\{   \frac{D( p_{\mathcal{X}^{*}} \,\|\, p_{\mathcal{X}} ) }{n_1  \cdots n_N}
			+  \left( \kappa +\frac{4 E_{KL} \ln2 }{3} (\beta+1) \ln \left(n_{\max }\right) \right)  \\
			&   \frac{  \hat{N}\Vert \mathbf{A}^{(N)}\Vert _{2,0} +R }{l} \Big\} 
			+   \frac{8 E_{KL} In l}{l},    
		\end{array}
	\end{equation*}
\end{theorem}
where $\hat{N} =\big(\sum_{i=1}^{N}n_i-N+1\big)  $  and  $ \kappa = \zeta 
\hat{N} \log_2 (\gamma) $.

\textit{Proof}. We utilize Lemma $\ref{lemma_43} $ to complete this proof. To make the penalty 
$ \text{pen}(\mathcal{X}) \geq 1 $ satisfy Kraft-Mcmillan inequality $(\ref{KM inequality})$, 
we consider encoding element-wise of $\mathbf{A}^{(i)}, i=1,\cdots,N$ as follows:

$\bullet$ For factor matrix $\mathbf{A}^{(N)} $, we employ a two-step procedure for encoding. Firstly, we encode the locations of the groups using a binary  number of length R. Then we encode the entries of the non-zero group using  $ \log_2(\gamma) $ bits. In this case, a total of $R+ n_{N} \Vert \mathbf{A}^{(N)} \Vert_{2,0} \log_2(\gamma) $ bits is used to encode $\mathbf{A}^{(N)} $.

$\bullet$ For factor matrix  $\mathbf{A}^{(i)}, i=1,\cdots,N-1$,  each element is encoded using $\log_2(\gamma)$ bits. The columns 
corresponding to the zero columns of $\mathbf{A}^{(N)}$ do not need 
to be encoded by the definition of the CP decomposition. So the total number of bits required to encode the elements in 
$\mathbf{A}^{(i)}$, $i=1, \allowbreak \ldots, 
N-1$, is
$
\sum_{i=1}^{N-1}
(n_i-1)\left\|\mathbf{A}^{(N)}\right\|_{2,0}\log_2(\gamma).
$

We define the penalty term $\text{pen}(\mathcal{X}) $ as 
\begin{equation}
	\operatorname{pen}(\mathcal{X})= \hat{N}\Vert \mathbf{A}^{(N)}\Vert _{2,0} \log_2 (\gamma)+{ R},
	\label{pen(x)}
\end{equation}
where $\hat{N} =\sum_{i=1}^{N}n_i-N+1$. Let $ {\mathbb{X}}' $  be a set of all tensors $\mathcal{X}  $ whose factor matrices $ A^{(i)} \in \mathbb{A}^{(i)} $, .i.e. 
\begin{equation}
	\begin{array}{rcl}
		\mathbb{{X}}'=\big\{\mathcal{X}= [\![\mathbf{A}^{(1)},\dots,\mathbf{A}^{(N)}] \! ] \mid \mathbf{A}^{(i)} \in \mathbb{A}_i, \; i=1,\cdots,N  \big\}.
	\end{array}
\end{equation}   

According to the uniqueness of the above encoding strategy and 
Kraft-Mcmillan inequality \cite{cover2012elements} , we have $\sum_{\mathcal{X} \in  {\mathbb{X}}' } 2^{-\operatorname{pen}(\mathcal{X})} \leq 1$, which implies
$ \sum_{\mathcal{X} \in \mathbb{X}} 2^{-\operatorname{pen}(\mathcal{X})}       \leq  \sum_{\mathcal{X} \in {\mathbb{X}}'} 2^{-\operatorname{pen}(\mathcal{X})}   
\leq 1 $ due to $\mathbb{X} \subseteq \allowbreak \mathbb{X}'$. By substituting $(\ref{pen(x)})$ into $(\ref{X_armin})$, we obtain the estimator of randomly subsampled and noisy observations $ \mathcal{Y}_{\Omega} $ as
\begin{equation}
	\begin{array}{rl}
		\mathcal{\widehat{X}}^{\zeta}_{}\left(\mathcal{Y}_{\Omega}\right)
		=&\underset{\mathcal{X} \in \mathbb{X}}{\arg \min }\Big\{-\ln p_{\mathcal{X}_{\Omega}}\left(\mathcal{Y}_{\Omega}\right)+ \zeta 
		\hat{N} log_2 (\gamma) \Vert \mathbf{A}^{(N)}\Vert _{2,0}  + \zeta R \Big\}  \\
		=& \underset{\mathcal{X} \in \mathbb{X}}{\arg \min }\Big\{-\ln p_{\mathcal{X}_{\Omega}}\left(\mathcal{Y}_{\Omega}\right)+ \zeta 
		\hat{N} \log_2 (\gamma) \Vert \mathbf{A}^{(N)}\Vert _{2,0}  \Big\}.
	\end{array}
\end{equation}
Together with Lemma $ \ref{lemma_43}$ and $ \zeta \geq (1+ \frac{2 E_{KL}}{3}) 2 \ln2 $, one has
\vspace{-6pt}
\begin{equation*}
	\begin{array}{rr}
		\frac{\mathbb{E}_{\Omega,\mathcal{Y}_{\Omega}}[\!-2 In( H(p_{\widehat{\mathcal{X} }^ {\zeta}},p_{\mathcal{X}^{*}} \! ))]}{n_1  \cdots n_N} 
		\!\!\! \!  \!  \! \! \! \! &\leq \!   3 \! \min \limits_{\substack{\mathcal{X}\in\mathbb{X} }}  \left\{ \! \!    \frac{D( p_{\mathcal{X}^{*}} \,\|\, p_{\mathcal{X}} \! ) }{n_1  \cdots n_N}
		\! + \!  ( \zeta \! + \! \frac{4 E_{KL} \! \! \ln2 }{3} \! ) \! \frac{ \log_2 ( \gamma ) \hat{N} \! \Vert \mathbf{A}^{(\! N \!)}\Vert _{2,0} \!+ \! {R} \!  }{l} \!  \right\} 
		\! + \!  \frac{8 \! E_{KL} \! \! \ln \! l}{l} \\
		&\leq \!    3 \! \min  \limits_{\substack{\mathcal{X}\in\mathbb{X} }}  \left\{ \! \!   \frac{D( p_{\mathcal{X}^{*}} \,\|\, p_{\mathcal{X}} \! ) }{n_1  \cdots n_N}
		\!+\!  (\! \zeta \! + \! \frac{4 E_{KL} \! \! \ln2 }{3} \! ) \! \log_2 (\! \gamma \!) \frac{  \hat{N} \! \Vert \mathbf{A}^{( \! N \! )}\Vert _{2,0} \! + \! {R}  \! }{l} \!  \right\} 
		\! +\!   \frac{8 \! E_{KL}  \!  \! \ln \!  l}{l}.   \\
		
	\end{array}
\end{equation*}
Set $\kappa=\zeta\widehat{N}\log_2(\gamma)$. By the definition of $\gamma$, we have $\log_2(\gamma)\leq(\beta+1)\ln(n_{\max})/ \ln 2.$
Consequently, if $ \kappa \geq
2\widehat{N}\left(1+\frac{2E_{\mathrm{KL}}}{3}\right)
(\beta+1)\ln(n_{\max}) $  , then the estimator 
\begin{equation*}
	\widehat{\mathcal{X}}_{}^{\kappa}=\underset{\mathcal{X}= [\![\mathbf{A}^{(1)}, \cdots, \mathbf{A}^{(N)}]\!] \in \mathbb{X}}{\arg\min} \{  -\ln p_{\mathcal{X}_{\Omega}}\left(\mathcal{Y}_{\Omega} \right)+\kappa \Vert \mathbf{A}^{(N)}\Vert _{2,0}    \} 
\end{equation*}
satisfies the error bound
\begin{equation*}
	\begin{array}{rl}
		\frac{\mathbb{E}_{\Omega,\mathcal{Y}_{\Omega}}[\!-2 In(H(p_{\widehat{\mathcal{X} }^ {\kappa}},p_{\mathcal{X}^{*}}\! ))]}{n_1  \cdots n_N} 
		\! \!  \! \!  \! \! &\leq \! 3 \! \min \limits_{\substack{\mathcal{X}\in\mathbb{X} }}  \left\{ \! \! \!   \frac{D( p_{\mathcal{X}^{*}} \,\|\, p_{\mathcal{X}} ) }{n_1  \cdots n_N}
		\! +  \!  \left(  \!  \kappa  \!+  \!\frac{4 E_{KL} \hat{N} (\beta +1) \ln \left(n_{\max }\right)  }{3}  \! \!  \right)  \! \!  \frac{   \Vert \mathbf{A}^{(N)} \! \Vert _{2,0} \! + \! \frac{{R}}{\hat{N}}   }{l}  \! \!  \right\} \\
		&+  \frac{8 E_{KL} \ln l}{l}.
	\end{array}
\end{equation*}
This completes the proof.  \hfill  $\square  $

Finally, we establish the error bound between the optimal solution of the model $(\ref{extend_model})$ and its ground truth under the assumption of additive Gaussian noise when the CP rank is unknown. The probability density function of $ \mathcal{Y}_{\Omega} $ is defined as
\begin{equation}\label{mo1gus}
	p_{\mathcal{X}_{\Omega}^{*}}(\mathcal{Y}_{\Omega})=\frac{1}{(2 \pi \sigma^2)^{|\Omega|/2}}  \exp \big(-\frac{1}{2 \sigma^2}\left\Vert \mathcal{Y}_{\Omega}-\mathcal{X}_{\Omega}^{*}   \right\Vert ^2\big),
\end{equation}
where observation $ \mathcal{Y}_{\Omega} $ follows a multivariate Gaussian distribution of dimension $ \vert \Omega  \vert $, 
whose mean corresponds to the tensor entries at the sampled locations and covariance matrix is $ \sigma^2 I_{|\Omega|}$, where $ I_{|\Omega|}  $ denotes the identity matrix of size $ \vert \Omega  \vert \times \vert \Omega  \vert  $.
Define  
\begin{equation}
	\beta := 1+ \frac{ \ln \big( 2R  \big(  \sum_{i=0}^{N-1} \binom{N}{i}   \big) +1   \big)}{\ln (n_{\max})}        .
	\label{beta}
\end{equation}
Applying Theorem $\ref{Bernoulli sampling}$, we can establish the error bound of model $(\ref{extend_model})$ in the discrete setting under additive Gaussian noise in the following theorem.
\begin{theorem}
	Let $ \beta $ be given by $(\ref{beta}) $. Suppose that the joint probability density function of $\mathcal{Y}_{\Omega}$ is given by $(\ref{mo1gus})$, and that 
	\begin{equation}
		\kappa =   2\hat{N}  (1+ \frac{2 E_{KL}}{3}) (\beta+1) \ln \left(n_{\max }\right).
		\label{kapp_va}
	\end{equation}
	Then, the estimator
	\begin{equation*}
		\widehat{\mathcal{X}}_{}^{\kappa}=\underset{\mathcal{X}= [\![\mathbf{A}^{(1)}, \cdots, \mathbf{A}^{(N)}]\!] \in \mathbb{X}}{\arg\min} \{  \frac{1}{2\sigma^2 } \Vert   \mathcal{Y}_{\Omega}  -  \mathcal{X}_{\Omega}    \Vert^2 + \kappa  \Vert \mathbf{A}^{(N)}\Vert _{2,0}    \} 
	\end{equation*}
	satisfies the error bound
	\begin{equation*}
		\begin{array}{rl}
			\frac{\mathbb{E}_{\Omega,\mathcal{Y}_{\Omega}^{}}[ \Vert  \widehat{\mathcal{X}}^{\kappa} - \mathcal{X}^{*}   \Vert_{F}^2  ]}{n_1  \cdots n_N} \! 
			\leq \!     \left(    (24 \sigma^2 \! + \! 64 L^2)  (\beta+1) \ln \left(n_{\max }\right)   \right)    \frac{ \hat{N}   \Vert \mathbf{A}^{(N)}\Vert _{2,0} +{R}}{l}   
			+  \frac{64L^2 \ln l}{l} \! + \! \frac{3L^2}{2 l^{2/N}}.    
		\end{array}
	\end{equation*}
	\label{error_bound}
\end{theorem}

\textit{Proof}. Fix a tensor $ \mathcal{X}^{*} $, and let $ \mathcal{X} \in \mathbb{X}$ be arbitrary. By the definition of the Hellinger affinity and the KL divergence and under the assumption of additive Gaussian noise, we have
\begin{equation*}\label{KLH}
	\begin{array}{ll}
		&D(p_{X_{i_1,\cdots,i_N}^*}\Vert p_{X_{i_1,\cdots,i_N}})=\frac{(X_{i_1,\cdots,i_N}^*-X_{i_1,\cdots,i_N})^2}{2\sigma^2},\\
		&-2\ln H(p_{X_{i_1,\cdots,i_N}},p_{X_{i_1,\cdots,i_N}^*})=\frac{(X_{i_1,\cdots,i_N}^*-X_{i_1,\cdots,i_N})^2}{4\sigma^2},   \forall (i_1,\cdots,i_N) \in \Omega.
	\end{array}
\end{equation*}
This implies that 
\begin{equation*}\label{KLH1}
	\begin{array}{ll}
		&D\left(p_{\mathcal{X}^*}\Vert p_{\mathcal{X}}\right)=\left\Vert \mathcal{X} -\mathcal{X}^*\right\Vert _F^2/2\sigma^2,\\
		&-2\ln H\left(p_{\mathcal{X}},p_{\mathcal{X}^*}\right)=\left\Vert \mathcal{X} -\mathcal{X}^* \right\Vert _F^2/4\sigma^2.
	\end{array}
\end{equation*}

Since the entries of both \(\mathcal{X}^{*}\) and \(\mathcal{X}\) are bounded in magnitude by \(L\), we may choose $ E_{KL} = \allowbreak  2L^2/ \sigma^2$ in Theorem $ \ref{Bernoulli sampling}$.
Using the value of  $\kappa $ given by $(\ref{kapp_va}) $ and applying Theorem $\ref{Bernoulli sampling}$, the group sparsity penalized maximum likelihood estimate satisfies the following mean square error bound
\begin{equation*}
	\begin{array}{rr}
		\frac{\mathbb{E}_{\Omega,\mathcal{Y}_{\Omega}^{}}[ \Vert  \widehat{\mathcal{X}}^{\kappa} \! -  \! \mathcal{X}^{*}   \Vert_{F}^2  ]}{n_1  \cdots n_N} 
		\! \!\! \leq \!  6  \!  \min \limits_{\substack{\mathcal{X}\in\mathbb{X} }}  \left\{  \! \!  \frac{\Vert \mathcal{X} \!- \mathcal{X}^{*}    \Vert_{F}^2 }{n_1  \cdots n_N}
		\! \! + \! \!  \left( \!  2\kappa \sigma^2  \! \!  + \! \! \frac{16L^2 \!  \hat{N} ( \! \beta \!+ \! 1 \!) \ln \left( \! n_{\max }\right)  }{3}  \! \!  \right) \!\!   \frac{   \Vert \! \mathbf{A}^{( \!N \!)} \! \Vert _{2,0}\!  +\! \frac{R}{\hat{N}}   }{l}  \! \!  \right\} 
		\! \! + \! \!  \frac{64L^2 \! \ln \! l }{l}.
	\end{array}   
\end{equation*}
For any $ \mathcal{X} \in \mathbb{X} ,$ we have
\begin{equation}
	\begin{array}{rl}
		\frac{\mathbb{E}_{\Omega,\mathcal{Y}_{\Omega}^{}}[ \Vert \widehat{\mathcal{X}}^{\kappa} -\mathcal{X}^{*} \!  \Vert_{F}^2  ]}{n_1  \cdots n_N} 
		\! \!\! \leq \!\! 6 \!   \left\{ \! \!   \frac{\Vert \mathcal{X} \!- \mathcal{X}^{*}    \Vert_{F}^2 }{n_1  \cdots n_N}
		\! + \!  \left( \!  2\kappa \sigma^2 \!+ \! \frac{16L^2 \! \hat{N} (\beta \!+ \! 1) \! \ln \left(n_{\max }\right)  }{3} \! \!   \right) \! \!   \frac{   \Vert \mathbf{A}^{(\! N \!) }\Vert _{2,0} \!+\! \frac{R}{\hat{N}}   }{l}  \! \! \right\}
		\!+\!  \frac{64   L^2 \! \ln \! l}{l}.
	\end{array}   
	\label{ten_ineq}
\end{equation}

Next, we construct the nearest discretized tensor
$\mathcal{X}_1^*=[\![\mathbf{A}_1^{(1)*},\ldots,\mathbf{A}_1^{(N)*}]\!]$ of $\mathcal{X}^{*} $ such that $\mathcal{X}_1^* \in \mathbb{X}$. For each $i=1,\ldots,N-1$, every entry of
$\mathbf{A}_1^{(i)*}$ is chosen as the nearest discretized approximations
to the corresponding entry of $\mathbf{A}^{(i)*}$. In the last factor
matrix, the entries in each nonzero column of $\mathbf{A}_1^{(N)*}$ are the nearest discretized approximation to the corresponding entry of $\mathbf{A}^{(N)*}$. Obviously,
$\|\mathbf{A}_1^{(N)*}\|_{2,0}
=\|\mathbf{A}^{(N)*}\|_{2,0}$. Next, we establish an error bound for $   \Vert \mathcal{X}_1^*-\mathcal{X}^* \Vert_{\infty}.  $ To this end, we consider the matricization of the tensor across the N-th dimension as follows:
\begin{equation*}
	\begin{array}{ll}
		\left\Vert \mathcal{X}_1^*-\mathcal{X}^*\right\Vert _{\infty} \! \! 
		=  \! \big\Vert (\mathbf{A}_1^{(N \!- \! 1)^*}\odot\cdots\odot\mathbf{A}_1^{(1)^*})\mathbf{A}_1^{{(N)^*}^T} \! \! \!\!\! - \!(\mathbf{A}^{(N \!- \!1)^*}\odot\cdots\odot\mathbf{A}^{(1)^*})\mathbf{A}^{{(N)^*}^T}\big\Vert _{\infty}.
	\end{array}
	\label{true_X}
\end{equation*}
Using the identities $\mathbf{A}^{(i)^*}_1=\mathbf{A}^{(i)^*}+\Delta_{\mathbf{A}}^{(i)},i=1,\cdots,N$, together $\Vert \mathbf{X}\odot\mathbf{Y}\Vert _{\infty}\leq  \Vert \mathbf{X}\Vert _{\infty}\Vert \mathbf{Y}\Vert _{\infty}$,  $\Vert \mathbf{X}\mathbf{Y}\Vert _{\infty}\leq R\Vert \mathbf{X}\Vert _{\infty}\Vert \mathbf{Y}\Vert _{\infty}$ and the triangle inequality for the infinity norm, we deduce that
\begin{equation*}
	\begin{array}{clc}
		\left\Vert \mathcal{X}_1^*-\mathcal{X}^*\right\Vert _{\infty} \leq 
		R\Big(\prod_{i=1}^{N}(\left\Vert \Delta_{\mathbf{A}}^{(i)}\right\Vert _{\infty}+\Vert \mathbf{A}^{(i)*}\Vert _{\infty})-\prod_{i=1}^{N}\Vert \mathbf{A}^{(i)*}\Vert _{\infty}\Big).
	\end{array}
\end{equation*}
By the definition of $\mathbb{X}  $, we have $\Vert \mathbf{A}^{(i)}\Vert _{\infty}\leq 1,i=1,\cdots,N-1$, and $\Vert \mathbf{A}^{(N)}\Vert _{\infty}\leq L$. 
Furthermore, utilizing  the fact that $\left\Vert \Delta_\mathbf{A}^{(i)}\right\Vert _{\infty} \leq \frac{\Vert \mathbf{A}^{(i)}\Vert _{\infty}}{\gamma-1}\leq \frac{1}{\gamma-1},i=1,\cdots,N-1$, and $\left\Vert \Delta_\mathbf{A}^{(N)}\right\Vert _{\infty} \leq \frac{L}{\gamma-1}$, we obtain 
\begin{equation*}
	\begin{array}{rl}
		\left\Vert \mathcal{X}_1^*-\mathcal{X}^*\right\Vert _{\infty} & \leq R\big((\frac{1}{\gamma-1}+1)^{N-1}(\frac{L}{\gamma-1}+L)-L\big)\\
		& \leq RL\big((1+\frac{1}{\gamma-1})^N-1\big)       \\
		& =  \frac{RL}{\gamma-1}   \left[ \binom{N}{0} \big( \frac{1}{\gamma-1} \big)^{N-1} + \cdots +  \binom{N}{N-1}  \right] \\
		& \leq  \frac{RL \big(  \sum_{i=0}^{N-1} \binom{N}{i}   \big)  }{\gamma-1}, 
	\end{array}
\end{equation*}
where the last inequality holds due to $\gamma \geq 2 $. By the choice of $\beta $ in $(\ref{beta})  $, we have 
$\gamma \geq  \big(  \sum_{i=0}^{N-1} \binom{N}{i}   \big) \allowbreak  + 1 .$
Consequently, 
$\Vert \mathcal{X}_1^*-\mathcal{X}^*    \Vert_{\infty} \leq L/2$. Combining this bound with \(\|\mathcal X^*\|_\infty\leq L\) and applying the triangle inequality gives
$ \Vert   \mathcal{X}_1^*  \Vert_{\infty} \leq \frac{3L}{2} < 2L$, which implies  $\mathcal{X}_1^* \in \mathbb{X} .$ Furthermore, we have
\begin{equation}
	\gamma \geq (n_{max})^{\beta} =n_{max}\Big(2R \big(  \sum_{i=0}^{N-1} \binom{N}{i}   \big) + 1 \Big) \geq n_{max} \Big(2R \big(  \sum_{i=0}^{N-1} \binom{N}{i}   \big)\Big) + 1,
\end{equation}
where the first inequality holds due to the definition of $\gamma $ and the equality follows from the properties of logarithms and the definition of $\beta$. Therefore, we have
\begin{equation}
	\gamma - 1 \geq n_{max} \Big(2R \big(  \sum_{i=0}^{N-1} \binom{N}{i}   \big) \Big). 
	\label{gamma_ineq}
\end{equation}
Using the inequality $(\ref{gamma_ineq})$, we obtain
\begin{equation}
	\begin{array}{cc}
		\frac{\Vert \mathcal{X}^{*}_{1} - \mathcal{X}^{*}   \Vert_{F}^2}{n_1  \cdots n_N} \leq \Vert \mathcal{X}^{*}_{1} - \mathcal{X}^{*}   \Vert_{\infty}^2  \leq   \left(     \frac{RL \big(  \sum_{i=0}^{N-1} \binom{N}{i}   \big)  }{\gamma-1}      \right)^2 \leq \frac{L^2}{4 n_{max}^2} \leq \frac{L^2}{4 l^{2/N}},
		\label{4_21}
	\end{array}
\end{equation}
where the first inequality holds since the root mean square of the tensor entries does not exceed their maximum absolute value. The last inequality follows from the sampling number $l$ satisfying $l \leq n_{\max}^{N}$. Due to $\mathcal{X}_1^{*} \in \mathbb{X}$ and $(\ref{4_21})$, substituting $ \mathcal{X}_1^*  $ for $  \mathcal{X}$ in $ (\ref{ten_ineq})$, we have
\begin{equation}
	\begin{array}{rl}
		\frac{\mathbb{E}_{\Omega,\mathcal{Y}_{\Omega}^{}}[ \Vert  \widehat{\mathcal{X}}^{\kappa} - \mathcal{X}^{*}   \Vert_{F}^2  ]}{n_1  \cdots  n_N} 
		\!\!\! \! \! &\leq    6  \left\{  \frac{L^2}{4 l^{2/N}}
		\!+ \! \left( 2\kappa \sigma^2 \! + \! \frac{16L^2  \hat{N} (\beta\! + \! 1) \ln \left(n_{\max }\right)   }{3}  \right)   \frac{   \Vert \mathbf{A}^{(N)^{*}}\Vert _{2,0}\! + \!   \frac{R}{\hat{N}}}{l}   \right\} 
		\!+\!  \frac{64L^2 \ln l}{l}  \\
		& = \! \!     \left(    (24 \sigma^2 \!+ \! 64 L^2)  (\beta+1) \ln \left(n_{\max }\right)  \right)    \frac{ \hat{N}    \Vert \mathbf{A}^{(N)^{*}}\Vert _{2,0} \! + \! {R}  }{l} \! + \! \frac{64L^2 \ln l}{l}  +\frac{3L^2}{2 l^{2/N}}, \\

	\end{array}
	\label{error bound}
\end{equation}
where the equality holds due to the value of $\kappa $ given by $(\ref{kapp_va})$.
This completes the proof.   \hfill  $\square$ 
\begin{remark}
	It follows from $(\ref{error bound})$ that the error bound depends primarily on the sampling size, the true CP rank of the tensor, and the initial CP rank. Increasing the sampling size yields a smaller error bound and higher reconstruction accuracy. Moreover, the error bound can be further reduced when the underlying tensor is intrinsically low rank. 
\end{remark}

	\section{Block-coordinate descent algorithm and convergence analysis}
The section introduces block-coordinate proximal gradient descent algorithm to solve the problem $(\ref{4})$. We first propose block-coordinate proximal gradient descent algorithm in a more general problem setting in subsection 5.1, analyze its convergence in subsection 5.2, and finally apply it to problem $(\ref{4})$. 
\subsection{Block-coordinate descent algorithms with extrapolation}
In this subsection, we introduce block-coordinate descent algorithms with extrapolation to solve the unconstrained optimization problem in the form of 
\begin{equation}
	\begin{array}{cc}
		\mathop{\text{min}}\limits_{\mathbf{X},\mathbf{Y},\mathbf{Z}} F(\mathbf{X},\mathbf{Y},\mathbf{Z}) =
		f(\mathbf{X},\mathbf{Y},\mathbf{Z}) + r_1(\mathbf{X}) + r_2(\mathbf{Y}) + r_3(\mathbf{Z}),
	\end{array}
	\label{gene_Cp_MODEL}
\end{equation}
where $ (\mathbf{X},\mathbf{Y},\mathbf{Z}) \in \mathbb{W}:=(\mathbb{X},\mathbb{Y},\mathbb{Z}) $
with $\mathbf{X} = (\mathbf{x}_{1},\dots,\mathbf{x}_{R})$, $\mathbf{Y} = (\mathbf{y}_{1},\dots,\mathbf{y}_{R}) $ and $\mathbf{Z} = (\mathbf{z}_{1},\dots,\mathbf{z}_{R}) $ being a group of vectors in Euclidean spaces (real infinite dimension) $ \mathbb{X} = \mathbb{X}_{1} \times \dots \times \mathbb{X}_{R} $,     $ \mathbb{Y} = \allowbreak \mathbb{Y}_{1} \times \dots \times \mathbb{Y}_{R} $ and  $ \mathbb{Z} = \mathbb{Z}_{1} \times \dots \times \mathbb{Z}_{R} $, respectively. $ f:\mathbb{W} \rightarrow \mathbb{R} $ is smooth (not necessarily convex) function, and these regularization term $r_1(\mathbf{X})  $, $r_2(\mathbf{Y})  $ and $r_3(\mathbf{Z})  $ are proper lower semicontinuous (not necessarily convex) function with the block-separable structure that  $ r_1(\mathbf{X}) = \sum_{i=1}^{R}r_{1,i}(\mathbf{x}_{i}) $,  $ r_2(\mathbf{Y}) = \sum_{j=1}^{R}r_{2,j}(\mathbf{y}_{j}) $ and $ r_3(\mathbf{Z}) = \sum_{k=1}^{R}r_{3,k}(\mathbf{z}_{k}) $.

Algorithm $\ref{algorithm1}$, designed to solve  problem $(\ref{gene_Cp_MODEL})$, is a block-coordinate proximal gradient descent method with extrapolation, where the inner subproblems are approximately solved by executing Algorithm $\ref{Algorithm2}$ for a finite number of iterations.
Here, $\mathbf{x}_{j} $, $\mathbf{y}_{j} $ and $\mathbf{z}_{j} $  represent the $j$-th block of $\mathbf{X}$, $\mathbf{Y}$ and $\mathbf{Z}$,  respectively. m is defined as the maximum number of inner iterations that is used to solve the subproblem involving all block variables. In Algorithm $\ref{algorithm1}$, we denote the value of $\mathbf{X}$ after the $k$-th outer iteration by $\mathbf{X}^{k+1}_{} $, and represent the value of $\mathbf{x}_{j}$ in the $k$-th outer iteration and after the $u$-th inner iteration by using $\mathbf{x}_{j}^{k+1,u}$. This can be expressed similarly for $\mathbf{Y}  $ and $\mathbf{Z} $.  For the sake of convenience, we denote 
\begin{equation}
	\begin{array}{cc}
		\mathbf{x}_{<j}:=(\mathbf{x}_{1},\dots,\mathbf{x}_{j-1}),\mathbf{y}_{<j}:=(\mathbf{y}_{1},\dots,\mathbf{y}_{j-1}),\mathbf{z}_{<j}:=(\mathbf{z}_{1},\dots,\mathbf{z}_{j-1}), \\
		\mathbf{x}_{ \geq j}:=(\mathbf{x}_{j},\dots,\mathbf{x}_{R}),\mathbf{y}_{ \geq j}:=(\mathbf{y}_{j},\dots,\mathbf{y}_{R}),\mathbf{z}_{ \geq j}:=(\mathbf{z}_{j},\dots,\mathbf{z}_{R}).
	\end{array}
\end{equation}
Algorithm $\ref{Algorithm2}$ is used to implement the function sub\_BCPGD in Algorithm $\ref{algorithm1}$. 
In particular, we observe that Algorithm $\ref{Algorithm2}$ is a block-coordinate proximal gradient descent algorithm, carried out on the subproblems: 
\begin{equation}
	\begin{array}{cc}
		\mathop{\text{min}}\limits_{\mathbf{X}} 
		\big\{f(\mathbf{X},\mathbf{Y},\mathbf{Z}) + r_1(\mathbf{X}) \big\}
		,\mathop{\text{min}}\limits_{\mathbf{Y}} 
		\big\{f(\mathbf{X},\mathbf{Y},\mathbf{Z}) +  r_2(\mathbf{Y}) \big\} ,\mathop{\text{min}}\limits_{\mathbf{Z}} 
		\big\{f(\mathbf{X},\mathbf{Y},\mathbf{Z}) +r_3(\mathbf{Z})\big\}.
	\end{array}
	\label{subproblem}
\end{equation}
At the $u$-th circle of iterations $(0 \leq u \leq m) $, we have each $ \textbf{x}_{j}^{u+1} $, $ \textbf{y}_{j}^{u+1} $ and $ \textbf{z}_{j}^{u+1} $  by
\begin{equation}
	\begin{array}{cc}
		\textbf{x}_{j}^{u+1} \in  \mathop{\text{argmin}}\limits_{\mathbf{\mathbf{x}_{j}}} \{ \langle  \nabla
		_{\mathbf{x}_{j}}f(\mathbf{x}_{<i}^{u+1},\mathbf{x}_{ \geq j}^{u},\mathbf{Y},\mathbf{Z}), \mathbf{x}_{j}-\mathbf{x}_{j}^{u}    \rangle    
		+  \frac{\epsilon+l_{1,j}^{u}}{2} \Vert \mathbf{x}_{j}-\mathbf{x}_{j}^{u}   \Vert^2    
		+ r_{1,j}(\mathbf{x}_{j})  \},   \\

		\textbf{y}_{j}^{u+1} \in  \mathop{\text{argmin}}\limits_{\mathbf{\mathbf{y}_{j}}} \{ \langle  \nabla
		_{\mathbf{y}_{j}}f(\mathbf{X},\mathbf{y}_{<j}^{u+1},\mathbf{y}_{ \geq j}^{u},\mathbf{Z}), \mathbf{y}_{j}-\mathbf{y}_{j}^{u}    \rangle    
		+  \frac{\epsilon+l_{2,j}^{u}}{2} \Vert \mathbf{y}_{j}-\mathbf{y}_{j}^{u}   \Vert^2    
		+ r_{2,j}(\mathbf{y}_{j})  \},  \\

		\mathbf{z}_{j}^{u+1} \in  \mathop{\text{argmin}}\limits_{\mathbf{\mathbf{z}_{j}}} \{ \langle  \nabla
		_{\mathbf{z}_{j}}f(\mathbf{X},\mathbf{Y},\mathbf{z}_{<j}^{u+1},\mathbf{z}_{ \geq j}^{u}), \mathbf{z}_{j}-\mathbf{z}_{j}^{u}    \rangle    
		+  \frac{\epsilon+l_{3,j}^{u}}{2} \Vert \mathbf{z}_{j}-\mathbf{z}_{j}^{u}   \Vert^2    
		+ r_{3,j}(\mathbf{z}_{j})  \},  \\
		
	\end{array}
	\label{proximal}
\end{equation}
where $ l_{1,j}^{u}$,$ l_{2,j}^{u}$ and $ l_{3,j}^{u}$ are the Lipschitz constant of the gradient in $(\ref{proximal})$, which can be written a proximal mapping. 

Now, we describe in detail how Algorithms $\ref{algorithm1}$ and $\ref{Algorithm2}$ are utilized to solve the proposed model. Let 
$ \mathbf{X} = \mathbf{A}$, $\mathbf{Y} = \mathbf{B}$, $ \mathbf{Z} = \mathbf{C}$, $\mathbf{x}_{j}=\mathbf{a}_{j}$, $ \mathbf{y}_{j}=\mathbf{b}_{j}$, $\mathbf{z}_{j}=\mathbf{c}_{j}$,
\vspace{-6pt}
\begin{equation}
	\begin{array}{cc}
		f(\mathbf{A},\mathbf{B},\mathbf{C}) = \frac{1}{2} \Vert \mathcal{M}*(\mathcal{Y}-[ \![\mathbf{A},\mathbf{B},\mathbf{C}] \! ]) \Vert ^2_F,
		\label{ob_fun}
	\end{array}
\end{equation}
and
\begin{equation}
	\begin{array}{cc}   
		r_{1,j}(\mathbf{a}_j) = \sigma_{\mathbb{A}_j}(\mathbf{A}) , \,   
		r_{2,j}(\mathbf{b}_j) = \sigma_{\mathbb{B}_j}(\mathbf{B}), \, 
		r_{3,j}(\mathbf{c}_j) = \lambda \Vert \mathbf{c}_j \Vert^{0}, \, j=1,\dots,R. 
	\end{array}
	\label{in_fun}
\end{equation}
For notational convenience, we set
\begin{equation}
	\begin{array}{cc}
		\mathbf{D}^{(1)} = (\mathbf{C} \bigodot \mathbf{B}), \, \mathbf{D}^{(2)} = \mathbf{C}  \bigodot \mathbf{A}, \, \mathbf{D}^{(3)} = \mathbf{B} \bigodot \mathbf{A}.
	\end{array}
\end{equation}
By the definition of matrix Khatri-Rao product, we know that the column number of each $  \mathbf{D^{(i)}}  $ is R. Hence, $  \mathbf{D^{(i)}}  $ can be written $ \mathbf{D^{(i)}} = (\mathbf{d}^{(i)}_{1},\dots,\mathbf{d}^{(i)}_{R}).$ We can calculate that
\begin{equation}
	\begin{array}{cc}
		\begin{cases}
			\nabla_{\mathbf{a}_j}f(\mathbf{A},\mathbf{B},\mathbf{C}) = [ \mathbf{M}_{(1)} *(\mathbf{a}_j \mathbf{d}^{(1)T}_{j}+ \sum_{\substack{r=1 \\ r \ne j}}^{R} \mathbf{a}_r \mathbf{d}^{(1)T}_{r}- \mathbf{Y}_{(1)})    ]\mathbf{d}^{(1)}_{j},     \\

			\nabla_{\mathbf{b}_j}f(\mathbf{A},\mathbf{B},\mathbf{C}) = [ \mathbf{M}_{(2)} *(\mathbf{b}_j \mathbf{d}^{(2)T}_{j}+ \sum_{\substack{r=1 \\ r \ne j}}^{R} \mathbf{b}_r \mathbf{d}^{(2)T}_{r}- \mathbf{Y}_{(2)})    ]\mathbf{d}^{(2)}_{j},     \\
			
			\nabla_{\mathbf{c}_j}f(\mathbf{A},\mathbf{B},\mathbf{C}) = [ \mathbf{M}_{(3)} *(\mathbf{c}_j \mathbf{d}^{(3)T}_{j}+ \sum_{\substack{r=1 \\ r \ne j}}^{R} \mathbf{c}_r \mathbf{d}^{(3)T}_{r}- \mathbf{Y}_{(3)})    ]\mathbf{d}^{(3)}_{j}  .  \\
		\end{cases}
	\end{array}
	\label{gradient}
\end{equation}
From $(\ref{gradient})$, the Lipschitz constant of $\nabla_{\mathbf{a}_j}f(\mathbf{A}, \! \mathbf{B},\!\mathbf{C}),     \nabla_{\mathbf{b}_j}f(\mathbf{A}, \! \mathbf{B},\! \mathbf{C})$ and $\nabla_{\mathbf{c}_j}f(\mathbf{A},\!\mathbf{B},\! \mathbf{C}) $ are in the form of, respectively
\begin{equation}
	\begin{array}{cc}
		l_{i,j} := \mathbf{d}^{(i)T}_{j} \mathbf{d}^{(i)}_{j},i=1,2,3.    
	\end{array}
\end{equation}

Furthermore, we need to compute explicitly the proximal mappings of $\sigma_{\mathbb{A}_{j}}(\cdot), 
\sigma_{\mathbb{B}_{j}}(\cdot) ,j=1,\dots,R$ and
$ \lambda \Vert\cdot \Vert^{0} $ in Step 4, 6, 8 of Algorithm $\ref{algorithm1}$, which are represented by
\begin{equation}
	\begin{array}{ccc}
		(\operatorname{prox}_{\gamma}^{ \lambda \Vert \cdot \Vert^{0}}(\mathbf{x})) = \begin{cases}  \mathbf{x}, & \text{if} \quad  \Vert \mathbf{x} \Vert > \sqrt{2 \lambda \gamma},    \\
			\{\textbf{0},\mathbf{x} \}   &  \text{if} \quad \Vert \mathbf{x} \Vert = \sqrt{2 \lambda \gamma},  \\
			\textbf{0}, & \text{if} \quad \Vert \mathbf{x} \Vert  < \sqrt{2 \lambda \gamma},  \\
		\end{cases} 
	\end{array}
	\label{prox_20}
\end{equation}
and 
\begin{equation}
	\begin{array}{ccc}
		(\operatorname{prox}_{\gamma}^{ \sigma_{\mathbb{X}_{j}}(\cdot) }(\mathbf{x})) = \begin{cases}  \{ \mathbf{e}| \Vert \mathbf{e} \Vert =1  \}, & \text{if} \quad  \Vert \mathbf{x} \Vert  =0,  \\
			\frac{\mathbf{x}}{\Vert \mathbf{x} \Vert }  &  \text{otherwise},  \\
		\end{cases} 
	\end{array}
	\label{prox_in}
\end{equation}
where $\mathbb{X}$ is $  \mathbb{A}$ or $ \mathbb{B}$ and $\mathbf{e} $ represents standard basis vector of the same dimension as $\mathbf{x}$. 
\begin{algorithm}[!h]
	\caption{ (BCPGD) A  block-coordinate proximal gradient descent algorithm with
		extrapolation for solving $(\ref{gene_Cp_MODEL})$}
	\label{algorithm1}
	
	\begin{algorithmic}[1]
		
		\REQUIRE{Set the positive integers $m_{i}>0,i=1,2,3$, and the real parameter $\epsilon >0$. Let $\bar{\mathbf{W}}_0:=\mathbf{W}_0 $, $m := max \{ m_{1},m_{2},m_{3}\}$, $\rho = \frac{\epsilon}{m}$,  $\varsigma := 0.99 $ and $\lambda_{0} := \hat{\lambda} $.}
		\ENSURE{$\mathbf{W}^* := (\mathbf{X}^{k+1},\mathbf{Y}^{k+1},\mathbf{Z}^{k+1}).$}
		\FOR{$ k=0,1,\dots $}
		\STATE \quad 
		$ \lambda_{k} := max(\varsigma*\lambda_{k-1},10^{-6}).$
		
		\STATE  \quad Update $ \mathbf{X}^{k+1} := \operatorname{sub\_BCPGD}((\mathbf{X}^{k},\bar{\mathbf{Y}}^{k},\bar{\mathbf{Z}}^{k}),m_1)$ by Algorithm $\ref{Algorithm2}$.
		
		\STATE \quad Update  $  \Tilde  {\mathbf{X} }^{k+1} := \mathbf{X}^{k+1} + \alpha_1 (\mathbf{X}^{k+1}-\bar{\mathbf{X}}^{k}).$

		\STATE  \quad  Update $\mathbf{Y}^{k+1} := \operatorname{sub\_BCPGD}( (\Tilde{\mathbf{X}}_{}^{k+1},{\mathbf{Y}}_{}^{k},\bar{\mathbf{Z}}^{k}),m_2)$ by Algorithm $\ref{Algorithm2}$.
		\STATE  \quad Update $\Tilde{ \mathbf{Y} }_{}^{k+1} := \mathbf{Y}^{k+1} + \alpha_2 (\mathbf{Y}^{k+1}-\bar{\mathbf{Y}}_{}^{k}).$
		\STATE  \quad Update $\mathbf{Z}^{k+1}:=\operatorname{sub\_BCPGD} ((\Tilde{\mathbf{X}}_{}^{k+1},\Tilde{\mathbf{Y}}^{k+1},{\mathbf{Z}}^{k}),m_3,\lambda_{k})$ by Algorithm $\ref{Algorithm2}$.
		\STATE \quad  Update $\Tilde { \mathbf{Z}}^{k+1} := \mathbf{Z}^{k+1} + \alpha_3 (\mathbf{Z}^{k+1}-\bar{\mathbf{Z}}^{k}).$
		
		\STATE \quad \textbf{if} $F(\mathbf{W}^{k+1}) > F(\mathbf{W}^{k})- \frac{\rho}{2}\Vert \mathbf{W}^{k+1}-\mathbf{W}^{k}  \Vert^2 $

		\STATE \;  \quad Update $\mathbf{X}^{k+1} :=\operatorname{sub\_BCPGD}((\mathbf{X}^{k},\mathbf{Y}^{k},\mathbf{Z}^{k}),m_1)$ by Algorithm $\ref{Algorithm2}$.
		\STATE \;  \quad Update $ \bar{\mathbf{X}}_{}^{k+1} := \mathbf{X}^{k+1} + \alpha_1 (\mathbf{X}^{k+1}-\bar{\mathbf{X}}^{k}).$

		\STATE  \;  \quad Update $\mathbf{Y}^{k+1} :=\operatorname{sub\_BCPGD}((\mathbf{X}^{k+1},\mathbf{Y}^{k},\mathbf{Z}^{k}),m_2)$ by Algorithm $\ref{Algorithm2}$.
		\STATE \;  \quad Update $\bar{\mathbf{Y}}^{k+1} := \mathbf{Y}^{k+1} + \alpha_2 (\mathbf{Y}^{k+1}-\bar{\mathbf{Y}}^{k}).$
		\STATE \;  \quad Update $\mathbf{Z}^{k+1}:=\operatorname{sub\_BCPGD}((\mathbf{X}^{k+1}, \! \!  \mathbf{Y}^{k+1}, \! \mathbf{Z}^{k}), \! m_3, \! \lambda_{k})$ by Algorithm $\ref{Algorithm2}$.
		\STATE \;  \quad Update $\bar{\mathbf{Z}}^{k+1} := \mathbf{Z}^{k+1} + \alpha_3 (\mathbf{Z}^{k+1}-\bar{\mathbf{Z}}^{k}).$
		\STATE \quad  \textbf{else}
		\STATE \quad  \quad Update $\bar{\mathbf{W}}^{k+1} := \Tilde{ \mathbf{W} }^{k+1}.$
		\STATE \quad  \textbf{end if}  
		\ENDFOR
	\end{algorithmic}
\end{algorithm}

\begin{algorithm}[!h]
	\caption{(sub\_BCPGD) A block-coordinate proximal gradient descent algorithm for
		solving subproblems $(\ref{subproblem})$ of Algorithm $\ref{algorithm1} $}
	\label{Algorithm2}
	
	\begin{algorithmic}[1]
		
		\REQUIRE{The initial point $ (\mathbf{V}_{1},\mathbf{V}_{2},\mathbf{V}_{3}) := \mathbf{W} = (\mathbf{X},\mathbf{Y},\mathbf{Z}) \in \mathbb{W},$ the max iteration number m, the index $i$ , and the real parameter $ \epsilon > 0$.}
		\ENSURE{$\mathbf{V}^{m+1}_{i} := (\mathbf{v}^{m_{}+1}_{i,1},\ldots,\mathbf{v}^{m_{}+1}_{i,R})$.}
		
		\STATE \textbf{Set} $\mathbf{V}^{(0)}_{i} = \mathbf{V}_{i}$.
		\FOR{$u=0,\dots,m$}
		\FOR{$j=1,\dots,R$}
		\STATE \qquad $ \mathbf{v}_{i,j}^{u+1} \in  \operatorname{prox}^{r_{i,j}}_{\frac{1}{\epsilon+l_{i,j}^{u}}}(\mathbf{v}_{i,j}^{u}-\frac{1}{\epsilon+l_{i,j}^{u}} \bigtriangledown_{\mathbf{v}_{i,j}}f(\mathbf{V}_{<i}^{},\mathbf{v}_{i,<j}^{u+1},\mathbf{v}_{i,\geq j}^{u},\mathbf{V}^{}_{>i})) .$
		\ENDFOR
		\ENDFOR
	\end{algorithmic}
\end{algorithm}

\begin{remark}
	Compared with Algorithm 4.1 in \cite{wang2025}, Algorithm  $\ref{algorithm1}$ adopts a different strategy for selecting the extrapolation parameter and requires less computational time.
\end{remark}

	\subsection{Convergence analysis of Algorithm}
This subsection presents the convergence analysis of Algorithm $\ref{algorithm1}$. We make the following assumptions.
\begin{assumption}
	\label{a_4}
	\textbf{(i)} $inf_{\mathcal{X}_1 \times \mathcal{X}_2 \times \mathcal{X}_3 } f > -\infty,inf_{\mathcal{X}_i} r_{i} > -\infty,i=1,2,3,  $ and the proximal mapping of each $r_i$ can be calculated directly.
	
	\textbf{(ii)} The partial gradient $\nabla_{\mathbf{X}} f(\mathbf{X},\mathbf{Y},\mathbf{Z})$, $\nabla_{\mathbf{Y}} f(\mathbf{X},\mathbf{Y},\mathbf{Z})$ and $\nabla_{\mathbf{Z}} f(\mathbf{X},\mathbf{Y},\mathbf{Z})$
	are globally Lipschitz when the other two variables are fixed. Their Lipschitz constants are denoted by $l_x,l_y,l_z$, respectively.
	
	\textbf{(iii)} There exist positive numbers $l_x^{+},l_x^{-},l_y^{+},l_y^{-}$ and $l_z^{+},l_z^{-}$ such that
	\begin{equation}
		\begin{array}{cc}
			sup\{l_x^{k},k \in N \} \leq l_x^{+}, sup\{l_y^{k},k \in N \} \leq l_y^{+},sup\{l_z^{k},k \in N \} \leq l_z^{+} ,\\ 
			inf\{l_x^{k},k \in N \} \geq l_x^{-}, inf\{l_y^{k},k \in N \} \geq l_y^{-},inf\{l_z^{k},k \in N \} \geq l_z^{-}. \nonumber
		\end{array}
	\end{equation}
	
	\textbf{(iv)} $\nabla f$ is Lipschitz continuous on bounded subsets of $ \mathcal{X}_1 \times \mathcal{X}_2 \times \mathcal{X}_3  $, i.e., for each bounded subset $\mathbb{S}$ of $ \mathcal{X}_1 \times \mathcal{X}_2 \times \mathcal{X}_3  $, there exists $M>0$, such that for all $\mathbf{W}^i=(\mathbf{X}^i,\mathbf{Y}^i,\mathbf{Z}^i)\in \mathbb{S},i=1,2,$ $\Vert \nabla f(\mathbf{W}^1)-\nabla f(\mathbf{W}^2)    \Vert \leq M \Vert \mathbf{W}^1- \mathbf{W}^2  \Vert $.
\end{assumption}
\begin{lemma} {\cite{bolte2014}}
	\label{dec}
	$  $ Let $f:\mathbb{R}^n \rightarrow \mathbb{R} $ be a continuously differentiable function with gradient $\nabla f$ assumed $L_{f}$-Lipschitz continuous. Then,
	\begin{equation}
		f(\mathbf{b}) \leq f(\mathbf{a}) + \langle \nabla f(\mathbf{a}),\mathbf{b}-\mathbf{a}\rangle+\frac{L_{f}}{2} \Vert \mathbf{b}-\mathbf{a} \Vert^2,\quad   \forall \mathbf{a},\mathbf{b} \in \mathbb{R}^n.   \nonumber
	\end{equation}
\end{lemma}

The following result is a descent lemma for Algorithm $\ref{algorithm1}$.  
\begin{lemma}
	\label{l_2}
	Suppose that Assumption $\ref{a_4} $ holds. Let the sequence $\{\mathbf{W}^k\}$ be generated by Algorithm $\ref{algorithm1} $ from the starting point $\mathbf{W}^0$. The following assertions hold:  \\
	\textbf{(i)} The sequence $\{F(\mathbf{W}^k)\}$ is non-increasing, and one has
	\begin{equation}
		F(\mathbf{W}^{k+1})+\frac{\rho}{2} \Vert \mathbf{W}^{k+1} -\mathbf{W}^{k})   \Vert^2 \leq   F(\mathbf{W}^{k}),\quad \forall k \geq 0.
		\label{decrease_inequality}
	\end{equation}
	\textbf{(ii)} It holds that $ \lim\limits_{k \to +\infty} \Vert \mathbf{W}^{k+1} -\mathbf{W}^{k}  \Vert = 0 $ and $ \lim\limits_{k \to +\infty} \Vert  \bar{\mathbf{W}}^{k} -\mathbf{W}^{k}  \Vert = 0 $.
\end{lemma}
The proof is provided in Appendix \ref{leamm_6_3}.

Next, we establish the convergence theorem for Algorithm $\ref{algorithm1}$. 
\begin{theorem}
	\label{limit_point}
	Suppose that Assumption $\ref{a_4}$ holds and the sequence $\{\mathbf{Z}^{k}\}_{k \in N}  $ is bounded. Let $\{\mathbf{W}^{k}\}_{k \in N}$ be the sequence generated by Algorithm $\ref{algorithm1}$ from the starting point $\mathbf{W}^{0}.$ Then, any accumulation point of  $\{\mathbf{W}^{k}\}_{k \in N}$ is a stationary point of problem $(\ref{gene_Cp_MODEL}).$ 
\end{theorem}
The proof is provided in Appendix \ref{theorem6_4}. According to Theorem $\ref{limit_point} $, it holds that each accumulation point of the sequence
generated by Algorithm $\ref{algorithm1}$ is a stationary point of problem (\ref{gene_Cp_MODEL}).
\begin{corollary}
	Supposed that the sequence  $\{\mathbf{C}^{k}\}_{k \in N}  $ is
	bounded. The sequence $\{ (\mathbf{A}^{k},\mathbf{B}^{k},\mathbf{C}^{k} )\}_{k \in N}$ generated by Algorithm $\ref{algorithm1}$  for solving $(\ref{4})  $ has at least one accumulation point, and any accumulation point of the sequence is a stationary point of problem $(\ref{4})  $. 
\end{corollary}

\textit{Proof.} Based on Theorem $\ref{limit_point}$, we can show that Assumption $\ref{a_4}$ holds. In addtion, (i) of Assumption $\ref{a_4}$ holds automatically from $(\ref{ob_fun}),(\ref{in_fun}),(\ref{prox_20}) $ and $(\ref{prox_in})$. Using $(\ref{ob_fun})$, one has
\begin{equation}
	\begin{array}{cc}
		\begin{cases}    \nabla_{\mathbf{A}}f(\mathbf{A},\mathbf{B},\mathbf{C} ) = [\mathbf{M}_{(1)}*(\mathbf{A} \mathbf{D}^{(1)T}-\mathbf{Y}_{(1)})] \mathbf{D}^{(1)},   \\
			\nabla_{\mathbf{B}}f(\mathbf{A},\mathbf{B},\mathbf{C} ) = [\mathbf{M}_{(2)}*(\mathbf{B}  \mathbf{D }^{(2)T}-\mathbf{Y}_{(2)})] \mathbf{D}^{(2)},   \\
			\nabla_{\mathbf{C}}f(\mathbf{A},\mathbf{B},\mathbf{C} ) = [\mathbf{M}_{(3)}*( \mathbf{C} \mathbf{D}^{(3)T}-\mathbf{Y}_{(3)})]\mathbf{D}^{(3)} . \\ 
		\end{cases}
	\end{array}
\end{equation}
It follows from the above equations that \textbf{(ii)} of Assumption $\ref{a_4} $ also holds.

According to the assumption and constraints, we can know that the sequence $\{(\mathbf{A}^{k},\mathbf{B}^{k},\mathbf{C}^{k} ) \allowbreak \}_{k \in N}$ is bounded. It can be inferred that $\{ (\mathbf{A}^{k},\mathbf{B}^{k},\mathbf{C}^{k} )\}_{k \in N}$ has least one accumulation point. Moreover, based on the smoothness of $f(\mathbf{A},\mathbf{B},\mathbf{C} ) $ and the boundedness of the generated sequence, we  deduce that (iii) and (iv) of Assumption $\ref{a_4}$ hold.
Hence, any accumulation point of $\{ (\mathbf{A}^{k},\mathbf{B}^{k}, \allowbreak \mathbf{C}^{k} )   \}_{k \in N}$ is a stationary point of problem $(\ref{4})  $  by Theorem $\ref{limit_point}$. This completes the proof.     \hfill  $\square   $

	\section{Numerical experiments}
In this section, the proposed method is applied to solve the component problem in chemometrics. We compare our method with four algorithms in the experiments: LRTCPR  \cite{jiang2022l}, INDAFAC \cite{tomasi2005}, BCPF \cite{Zhao2015}, and PARAFACM\cite{dong2024}. In our experiments, all the parameters for each method are optimized to ensure fair comparisons and peak performance. The source code is publicly available at \url{https://github.com/cmx269269/Tensor_CP_Decomposition}.  

\subsection{Experiments settings}
All experiments are performed on an desktop computer eqipped with AMD 7945 CPU, 64 GB of RAM and MATLAB 2024b. Each method is operated 30 times to obtain the average performance. \\
\indent For a given tensor $\mathcal{Y} \in \mathbb{R} ^{I \times J \times K} $  and
the sequence $\{\mathbf{A}^{k},\mathbf{B}^{k},\mathbf{C}^{k} \}$ generated by these methods, the relative error at $k$-th outer-iteration is defined as follows:
\begin{equation}
	ReEr_{k} = \frac{\Vert ({\mathcal{Y}}-[ \![\mathbf{A}^{k},\mathbf{B}^{k},\mathbf{C}^{k}] \! ]) \Vert}{\Vert \bar{\mathcal{Y} }\Vert} .
	\label{ReEr}
\end{equation}

The stopping criterion is based on the relative KKT residual defined by $\eta:=max   \{ \eta_1,\eta_2,\eta_3\} $, where
\begin{equation}
	\begin{array}{cc}
		
		\begin{cases}
			\eta_1:=\frac{\Vert \textbf{A} - prox_{1}^{\sigma_{\textbf{A}}(\cdot)}((\textbf{A}- \nabla_{\textbf{A}}f(\textbf{A},\textbf{B},\textbf{C})))     \Vert}{1+\Vert \textbf{A} \Vert + \Vert \nabla_{\textbf{A}}f(\textbf{A},\textbf{B},\textbf{C})) \Vert } ,\\
			\eta_2:= \frac{\Vert \textbf{B} - prox_{1}^{\sigma_{\textbf{B}}(\cdot)}((\textbf{B}- \nabla_{\textbf{B}}f(\textbf{A},\textbf{B},\textbf{C})))     \Vert}{1+\Vert \textbf{B} \Vert + \Vert \nabla_{\textbf{B}}f(\textbf{A},\textbf{B},\textbf{C})) \Vert } ,           \\
			\eta_3:= \frac{\Vert \textbf{C} - prox_{1}^{\lambda \Vert \cdot \Vert_{2,0}}((\textbf{C}- \nabla_{C}f(\textbf{A},\textbf{B},\textbf{C}))) \Vert}{1+\Vert \textbf{C} \Vert + \Vert \nabla_{\textbf{C}}f(\textbf{A},\textbf{B},\textbf{C})) \Vert }.                    \\
		\end{cases}
	\end{array}
	\nonumber
\end{equation}

We terminate the algorithms when $\eta  <10^{-9}.$

To compare the performance of different algorithms, we use root mean squared error of prediction (RMSEP) [13] to measure the recovery rates of the components in Section 7.4. The true and computed (after regression) concentration profiles matrix are denoted by $\mathbf{C}_{real} \in \mathbb{R}^{K \times R}  $ and $\mathbf{C} \in \mathbb{R}^{K \times R} $, respectively.  RMSEP is defined as follows:
\begin{equation}
	\begin{array}{cc}
		RMSEP = \sqrt{\frac{\sum_{i=1}^{K} \sum_{j=1}^{R}(\mathbf{C}_{real}(i,j)-\mathbf{C}(i,j))^2}{K \times R}  } .
	\end{array}
\end{equation}

\begin{remark}
	If the components cannot be separated correctly, we will regard the value of RMSEP as "-".
\end{remark}

\subsection{Parameters settings}
\begin{figure}[!htbp]
	\centering
	\subfloat[$\mathcal{X}_1 \quad \alpha_1$ ]{\includegraphics[width=4.2cm,height=3.6cm]{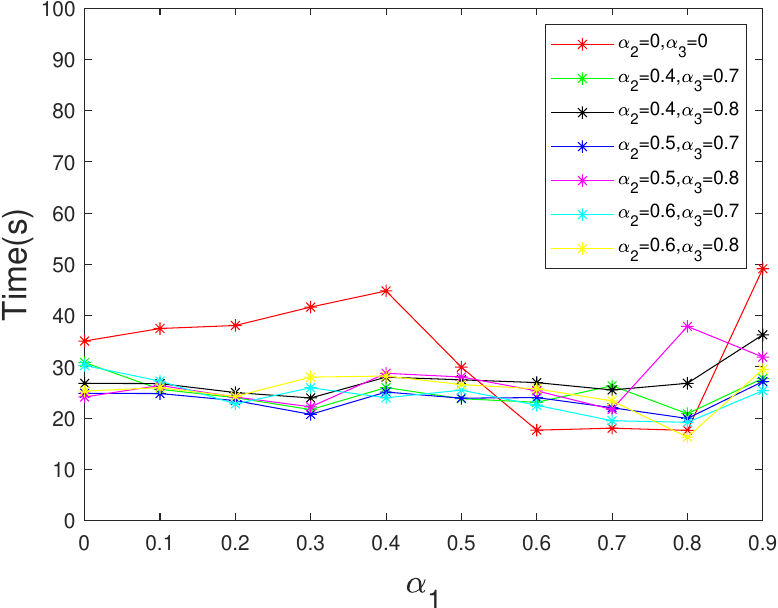}}
	\subfloat[$\mathcal{X}_1 \quad \alpha_2$]{ \includegraphics[width=4.2cm,height=3.6cm]{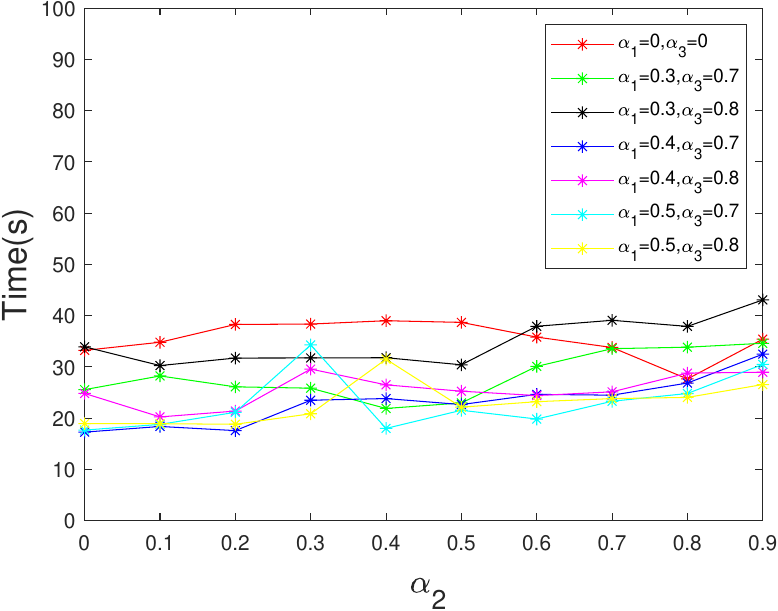}}
	\subfloat[$\mathcal{X}_1 \quad \alpha_3$ ]{ \includegraphics[width=4.2cm,height=3.6cm]{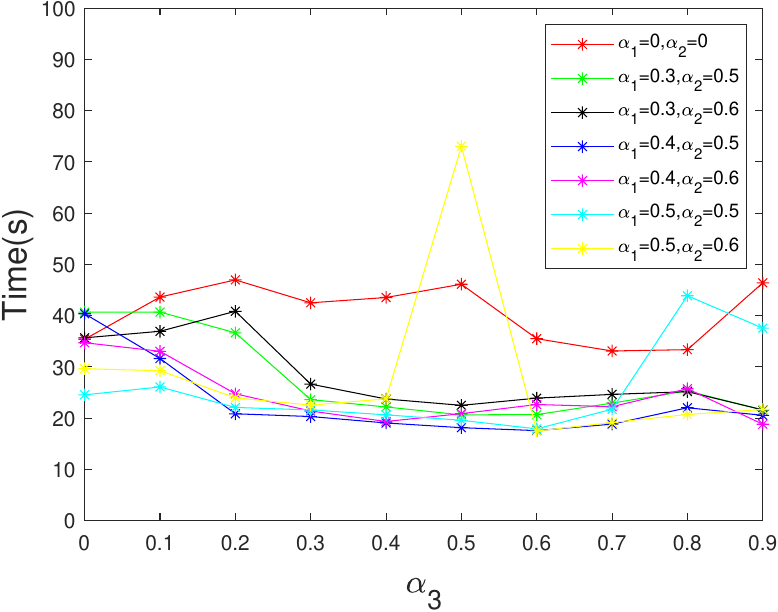}}
	
	\subfloat[$\mathcal{X}_2 \quad \alpha_1$]{\includegraphics[width=4.2cm,height=3.6cm]{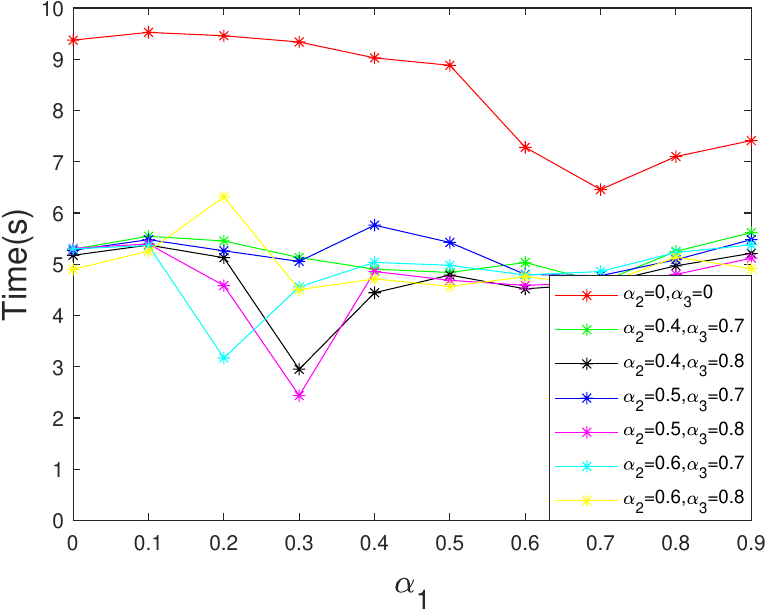}}
	\subfloat[$\mathcal{X}_2 \quad \alpha_2$]{ \includegraphics[width=4.2cm,height=3.6cm]{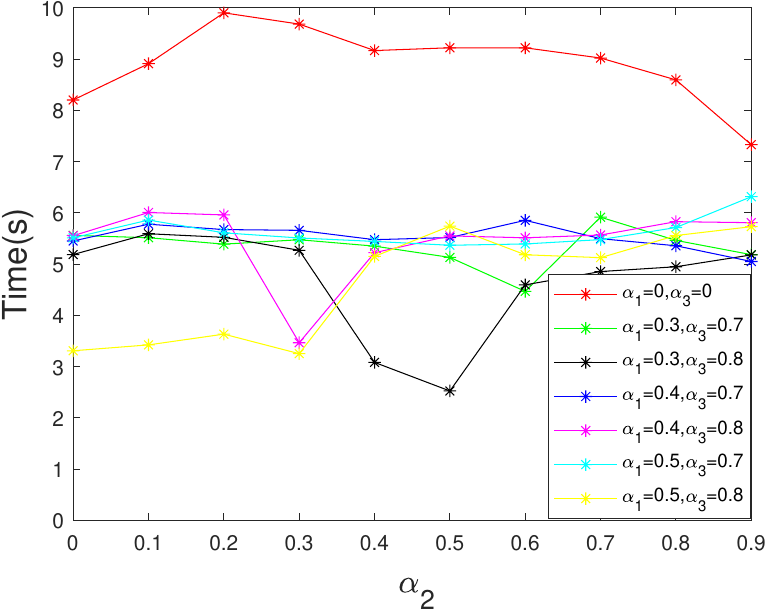}}
	\subfloat[$\mathcal{X}_2 \quad \alpha_3$]{\includegraphics[width=4.2cm,height=3.6cm]{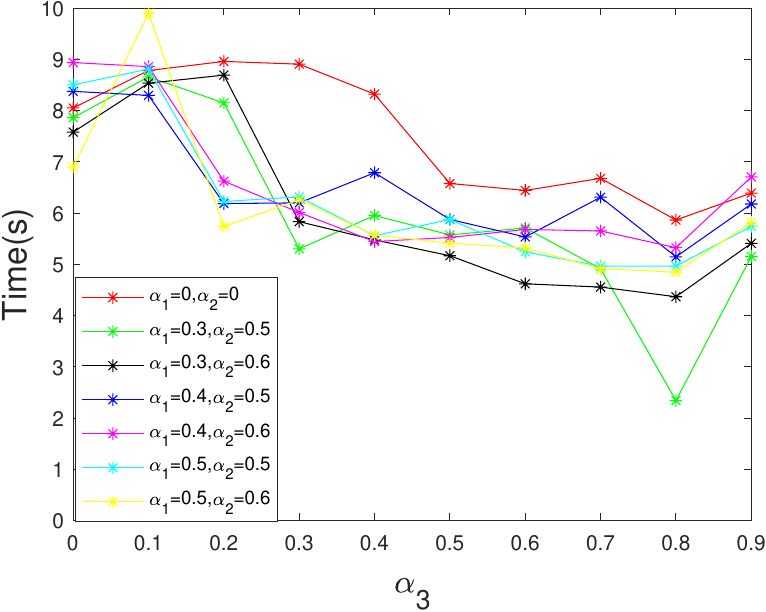}}
	\caption{ The computation time of BCPGDEM with different accelerated parameters setting for $\mathcal{X}_1 $ and $\mathcal{X}_2 $,  respectively}.  
	\label{fig1}
\end{figure}

\begin{figure}[!htbp]
	\centering
	\subfloat[$ \mathcal{X}_1 \in  \mathbb{R}^{50 \times 50 \times 20} $ ]{\includegraphics[width=5cm,height=4.2cm]{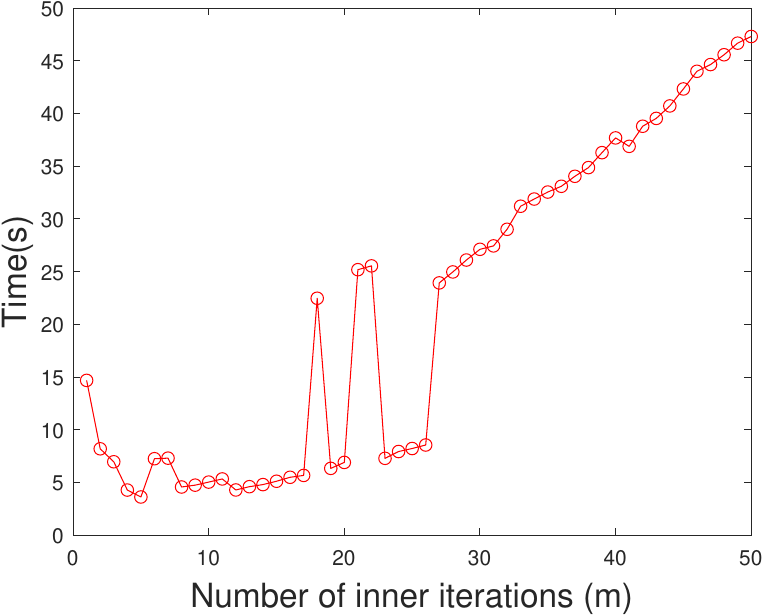}}\hspace{0.2cm}
	\subfloat[$  \mathcal{X}_2 \in \mathbb{R}^{52 \times 214 \times 16}$]{ \includegraphics[width=5cm,height=4.2cm]{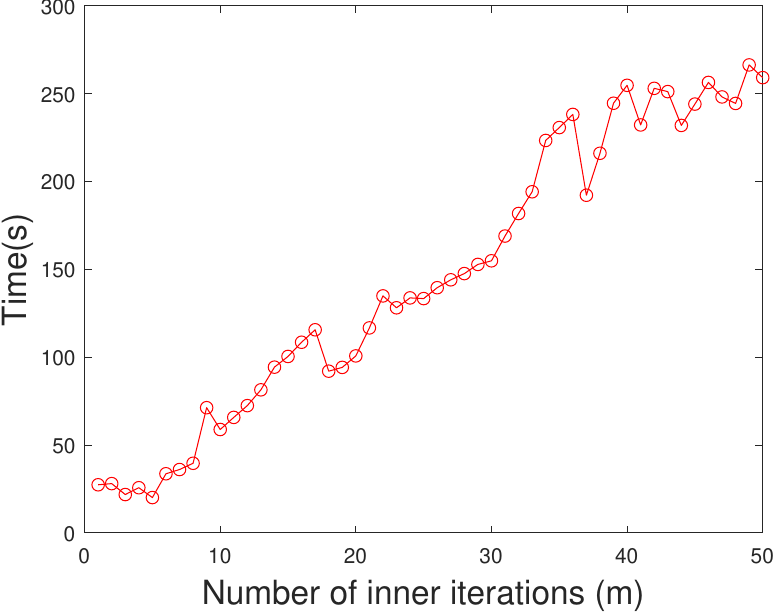}}\hspace{0.2cm}
	\caption{ \textit{The computation time of BCPGDEM under different numbers of innner iterations m setting for $\mathcal{X}_1 $ and $\mathcal{X}_2 $, respectively} }
	\label{fig2}
\end{figure}

In this section, we introduce how to choose the parameters of the proposed model, which include the regularization parameter  $\lambda $, initial points and three accelerated parameters  $\alpha_1,\alpha_2 $ and $\alpha_3 $ of the proposed algorithm, along with a  parameters $m$ which controls the numbers of inner-iterations. \\
$\textbf{Adaptive tuning strategy for the regularization parameter} $ $\bm{\lambda}$: The initial regularization parameter is set as
\[
\lambda_0 := n\sqrt{2\sigma^2\ln(200IJK)},
\]
where $\sigma^2$ can be estimated based on the residual $\|\mathcal{A}-\mathcal{Y}\|_F$, with $\mathcal{A}$ obtained from the observed tensor $\mathcal{Y}$ via generalized eigenvalue decomposition. Here, $n$ is a quantity depending on the tensor dimensions and is estimated using the following piecewise function:
\begin{equation}
	n =
	\begin{cases}
		1, & \text{if } IJK < Q,\\
		11+\lceil IJK/Q\rceil, & \text{otherwise},
	\end{cases}
	\label{n_value}
\end{equation}
where $Q=5\times 10^5$. The regularization parameter is then updated according to
$
\lambda_{k+1}=\max\{10^{-6},\,0.99\lambda_k\}
$.
\\
\textbf{Choice of acceleration parameters} $\bm{\alpha_1,\alpha_2  , \alpha_3} \textbf{:}   $ In order to evaluate the influence of the acceleration parameters on the computational cost, we use two datasets with different dimensions: a random tensor dataset $\mathcal{X}_1 \in  \mathbb{R}^{50 \times 50 \times 20} $ and a real data $\mathcal{X}_2 \in \mathbb{R}^{52 \times 214 \times 16}  $. We set the number of inner $\allowbreak $ -iterations to m=5, while fixing the initial estimated rank at $R=5$ and the sampling ratio at $SR =0.5$. The computational times of our proposed algorithm under different accelerated parameters $\alpha_1,\alpha_2 $ and $\alpha_3 $, ranging from 0 to 0.9, are illustrated in Figure $\ref{fig1}$. From these observations, we can see that the optimal three parameters are concentrated around 0,3, 0.5 and 0.8, respectively.
Hence, we select $\alpha_1 = 0.3,\alpha_2 = 0.5,\alpha_3 = 0.8$ as the accelerated parameters for the numerical experiments in the paper. 
\\
\textbf{Choice of the numbers of innner-iterations m:}
To investigate the relationship between computation time and  the numbers of innner-iterations m, we also use random tensor data $ \mathcal{X}_1$ and real data $ \mathcal{X}_2 $. The range of m is 1 to 50, and the initial estimated component number and sampling radio are the same as the previous part. The change of the computation time  with different numbers of inner iterations m are shown in Figure  $\ref{fig2}$. From these results, the computation time shows a downward trend as m increases from 1 to 5, whereas it rises with a further increase in the number of inner iterations $m $. Therefore, we adopt m=5 for all subsequent experiments. 
\\
\textbf{Initialization:} The initial point is obtained by performing generalized eigenvalue decomposition on the incomplete tensor, with its missing entries set to zero.\\

\subsection{Experiments on randomly generated low rank tensors}
In this subsection, the performance of the proposed method is assessed by using three randomly generated tensors, referred to as $\mathcal{X}_3 \in \mathbb{R}^{100 \times 50 \times 10} $, $\mathcal{X}_4 \in \mathbb{R}^{300 \times 100 \times 20} $, and $\mathcal{X}_5 \in \mathbb{R}^{400 \times 300 \times 30}  $. The true CP ranks of $\mathcal{X}_3 $,  $ \mathcal{X}_4 $ and $ \mathcal{X}_5 $ are 4, 5 and 7, respectively.
To demonstrate that the proposed model is not sensitive to R, we calculate the result of CP decomposition of them by utilizing overestimated initial component number R.  
The initial number of components \(R\) is set to 30 for each of $\mathcal{X}_3,\mathcal{X}_4,\mathcal{X}_5,$ respectively.  Table $\ref{random_tensor} $ presents tensor dimensions, the initial rank, the true rank, the computed rank, and the relative error in different sample radio  $SR =\{0.5,0.3,0.1\} .$ It can be seen that our proposed model and algorithm can successfully calculate the true rank of these tensors, which implies their effectiveness and robustness.

\begin{table}[!htbp]
	\footnotesize
	\caption{ \textit{The performance of Algorithm \ref{algorithm1} on randomly generated low rank tensors.} }
	\label{random_tensor}
	\centering
	\begin{tabular}{|c|c|c|c|c|c|c|c|}
		\hline
		tensor & Dimension & SR  &  Initial-R  & True-R & Compute-R & RelErr\\
		\hline
		\multirow{1}*{$\mathcal{X}_3$}
		& \multirow{1}*{$100 \times 50 \times 10  $} 
		
		& 0.5 & 30 & 4 & 4 & 5.0824e-08\\
		\hline
		\multirow{1}*{$\mathcal{X}_3$}
		& \multirow{1}*{$100 \times 50 \times 10  $}  & 0.3 & 30 & 4 & 4 & 7.8530e-08\\
		\hline
		\multirow{1}*{$\mathcal{X}_3$}
		& \multirow{1}*{$100 \times 50 \times 10  $}  & 0.1 & 30 & 4 & 4 & 1.9466e-07\\
		\hline
		\multirow{1}*{$\mathcal{X}_4$}
		&   \multirow{1}*{$300 \times 100 \times\ 20  $}  & 0.5 & 30 & 5 & 5 &3.2282e-08   \\
		\hline
		\multirow{1}*{$\mathcal{X}_4$}
		& \multirow{1}*{$300 \times 100 \times\ 20  $}  & 0.3 & 30 & 5 & 5 &5.1420e-08\\
		\hline
		\multirow{1}*{$\mathcal{X}_4$}
		&  \multirow{1}*{$300 \times 100 \times\ 20  $}  & 0.1 & 30 & 5 & 5 &1.6652e-07\\
		\hline
		\multirow{1}*{$\mathcal{X}_5$}
		&  \multirow{1}*{$400 \times 300 \times 30 $}  & 0.5 & 30 & 7 & 7 &2.6736e-08 \\
		\hline
		\multirow{1}*{$\mathcal{X}_5$}
		& \multirow{1}*{$400 \times 300 \times 30 $}  & 0.3 & 30 & 7 & 7 &4.4824e-08 \\
		\hline
		\multirow{1}*{$\mathcal{X}_5$}
		&  \multirow{1}*{$400 \times 300 \times 30 $}  & 0.1 & 30 & 7 & 7 &1.3574e-07 \\
		\hline
	\end{tabular} 
\end{table}
\subsection{Experiments on real data}
In this subsection, we conduct experiments on  two groups of chemical data, which  are composed of  two components data $(32 \times 39 \times 15)$ in  and macrocephalae rhizoma data $ (52 \times 214 \times 16)$ to evaluate the performance of various algorithms.

\subsubsection{Two-component data}
\begin{figure}[!htbp]
	\centering
	\subfloat[Real]{ \includegraphics[width=4.2cm,height=3.7cm]{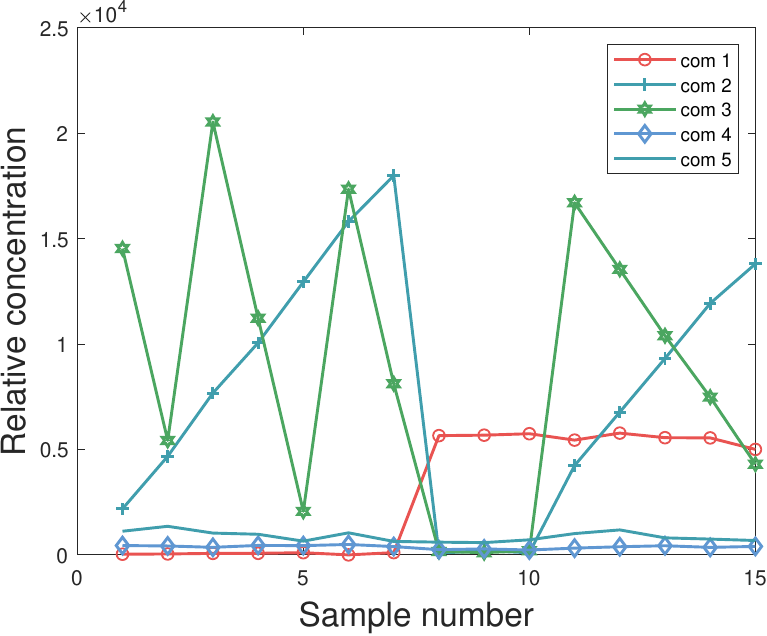}}
	\subfloat[INDAFAC]{\includegraphics[width=4.2cm,height=3.6cm]{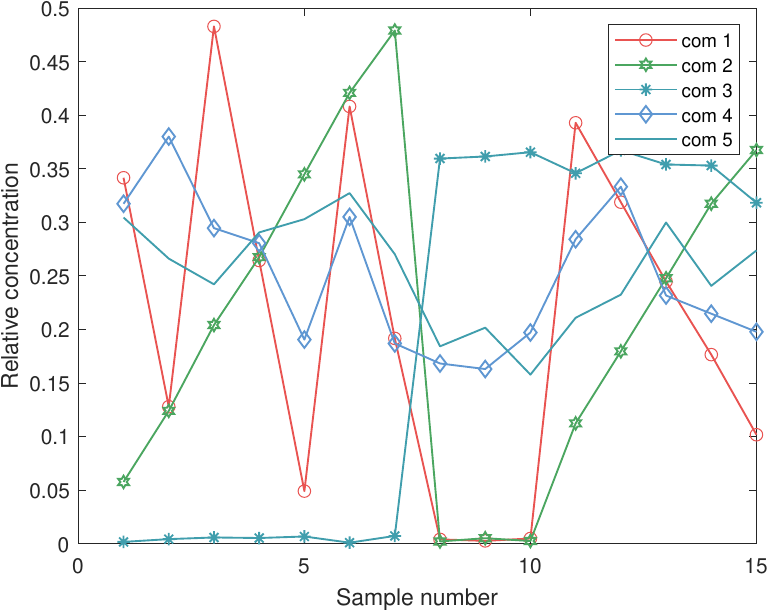}}
	\subfloat[PARAFACM]{ \includegraphics[width=4.2cm,height=3.7cm]{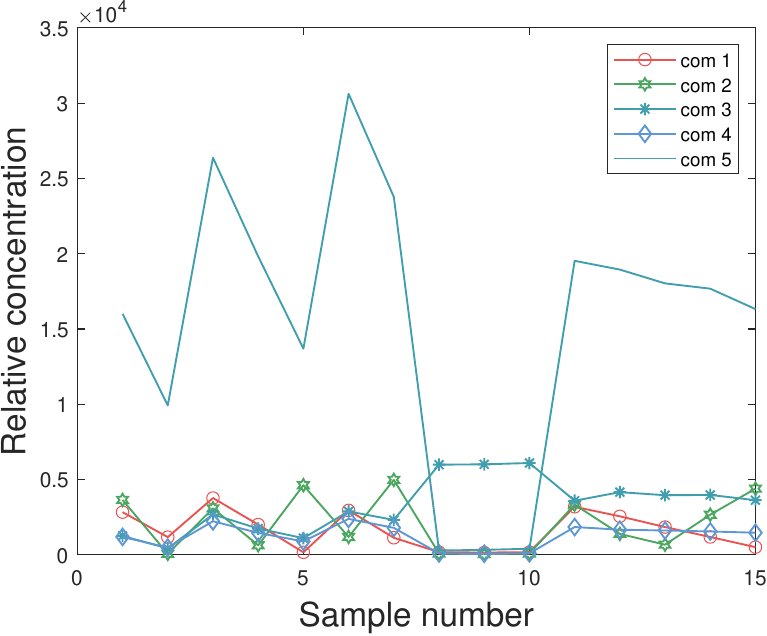}}
	
	\subfloat[LRTCPR]{\includegraphics[width=4.2cm,height=3.6cm]{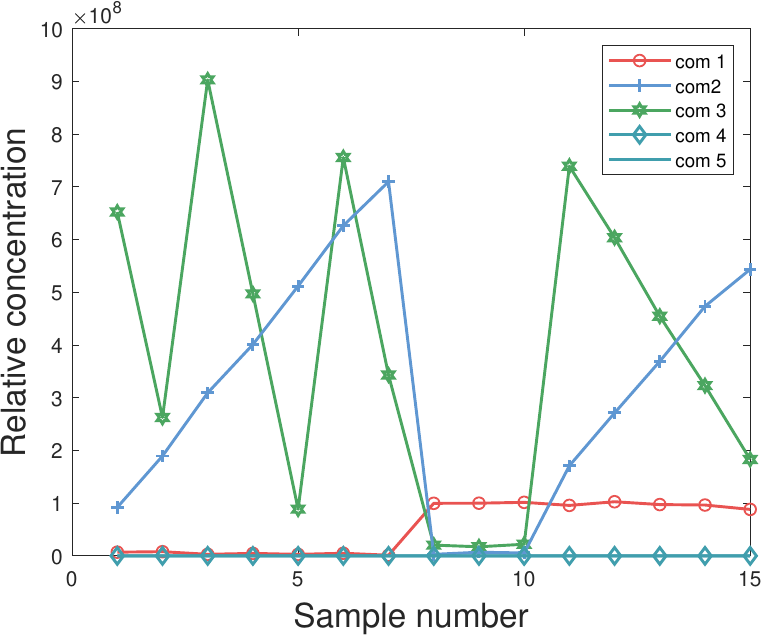}}
	\subfloat[BCPF]{\includegraphics[width=4.2cm,height=3.5cm]{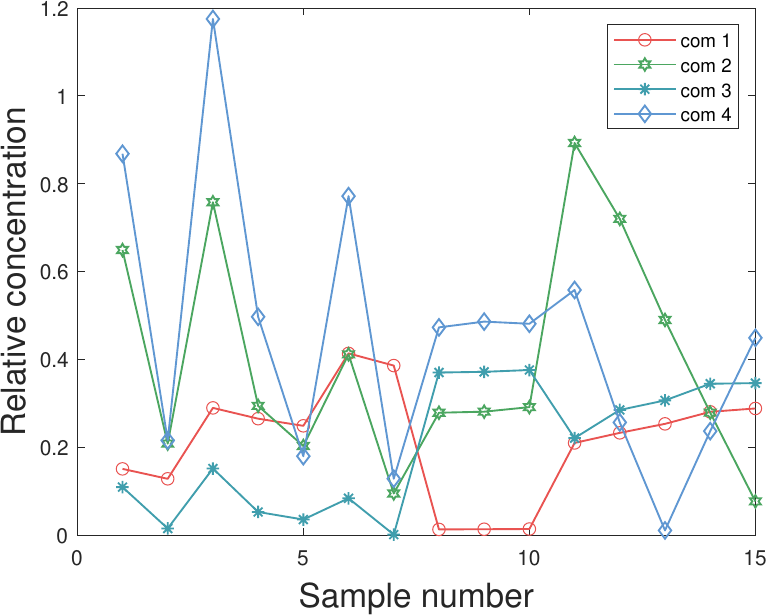}}
	\subfloat[BCPGD]{\includegraphics[width=4.2cm,height=3.6cm]{Figure/two_component/tc_03_5_C11111_deepened_thick_com2_plus_markers_continuous_thick.pdf}}
	
	\caption{ Analytical results with $R=5$ and
		$SR=0.3$ for the the two-component data. The relative concentration profiles resolved by (a) Real, (b) INDAFAC, (c) PARAFACM, (d) LRTCPR and (e)  BCPF and (f) BCPGD, respectively.} 
	\label{fig6}
\end{figure}

This data contains two components and exhibits background drift and therefore can be seen as consisting of three  components. We choose the number of components R = \{3,5,9,20\} for different sampling ratios $SR = \{0.5,0.3,0.1\}$. Table $\ref{two_com}$ reports the values of RMSEP and computation times of different algorithms. Figure $ \ref{fig6} $ demostrates the real concentration profiles and the relative concentration profiles resolved by different methods with $R = 5$ and $SR=0.3$, corresponding to C. 

Table \ref{two_com} shows that, PARAFACM, INDAFAC, LRTCPR and BCPF may not yield accurate component separation when the rank is overestimated, whereas BCPGD consistently demonstrates satisfactory results. In particular, LRTCPR can obtain acceptable results merely when the rank slightly exceeds the true rank.
Figure $ \ref{fig6} $ also shows that when R = 5, only BCPGD can separate three components exactly, underscoring the advantage of incorporating the group-sparse regularization.

\begin{table}[!htbp]
	\footnotesize
	\caption{ Comparison of root mean square error of prediction (RMSEP) of different algorithms and computation time (s) for the Two-component data in different sample radios $SR = \{0.5,0.3,0.1\}$   with R = 3,5,9,20. 
		"-" shows that the components cannot be separated correctly. } 
	\label{two_com}

	\scriptsize 
	\centering
	\begin{tabular}{|l|l|l|l|l|l|l|}
		\hline
		SR = 0.5 & Algorithm  & PARAFACM & INDAFAC & LRTCPR &BCPF & BCPGD\\
		\hline
		\multirow{1}*{R=3}  
		~ & RMSEP & 4.1619  & 4.0111 & \textbf{3.3707}  &27.9282 & 4.0065\\
		\hline
		\multirow{1}*{R=5}  
		~ & RMSEP &-  & - &7.5003 &-   &\textbf{4.5931} \\
		\hline
		
		\multirow{1}*{R=9}  
		~ & RMSEP &-  &  &- &-   &\textbf{4.3243} \\
		\hline

		\multirow{1}*{R=20} 
		~ & RMSEP & -   &- & -& - & \textbf{6.0561} \\
		\hline

		SR = 0.3  & Algorithm  & PARAFACM  & INDAFAC & LRTCPR &BCPF  &BCPGDEM\\
		\hline
		\multirow{1}*{R=3}  
		~ & RMSEP &  3.9704  & 3.8452 & 5.5178  & -& \textbf{3.8374}\\
		\hline
		\multirow{1}*{R=5}  
		~ & RMSEP &-  &- & 5.6479 &-  & \textbf{4.1049}\\
		\hline
		
		\multirow{1}*{R=9}  
		~ & RMSEP &-  & -&- &-   &\textbf{4.1247} \\
		\hline

		\multirow{1}*{R=20}  
		~ & RMSEP & -   &  &- & - & \textbf{4.1249} \\
		\hline

		SR = 0.1 & Algorithm  & PARAFACM  & INDAFAC & LRTCPR & BCPF & BCPGDEM\\
		\hline
		\multirow{1}*{R=3}  
		~ & RMSEP &  - & 4.3489 & 13.0055  & -& \textbf{3.7340} \\
		\hline
		\multirow{1}*{R=5}  
		~ & RMSEP & -  & - &4.2417& - & \textbf{3.4202}  \\
		\hline
		
		\multirow{1}*{R=9}  
		~ & RMSEP &-  &  &- &-   &\textbf{3.5648} \\
		\hline
		
		\multirow{1}*{R=20}  
		~ & RMSEP & -   &- &-&- &\textbf{3.1937} \\
		\hline
		
	\end{tabular} 
	
\end{table}
	\subsubsection{Macrocephalae rhizoma data} 
\begin{figure}[!htbp]
	\centering
	\subfloat[Real]{\includegraphics[width=4.2cm,height=3.6cm]{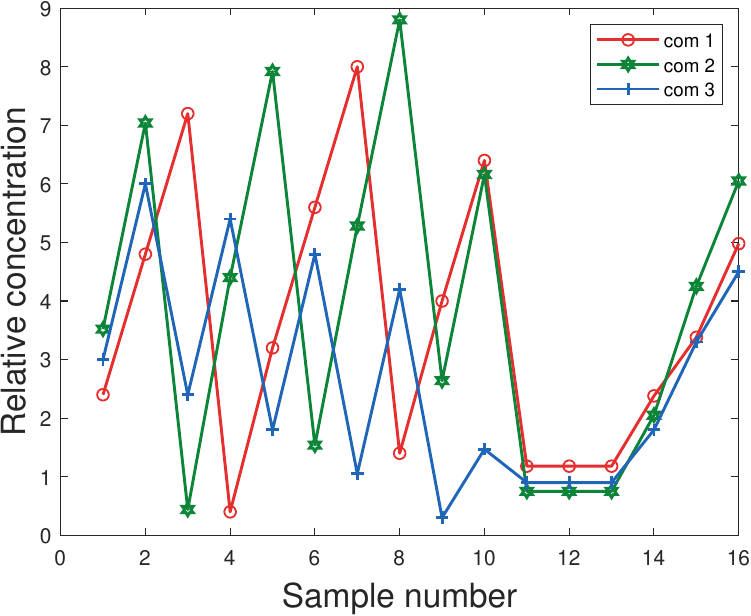}}\hspace{0cm}
	\subfloat[INDAFAC]{\includegraphics[width=4.2cm,height=3.6cm]{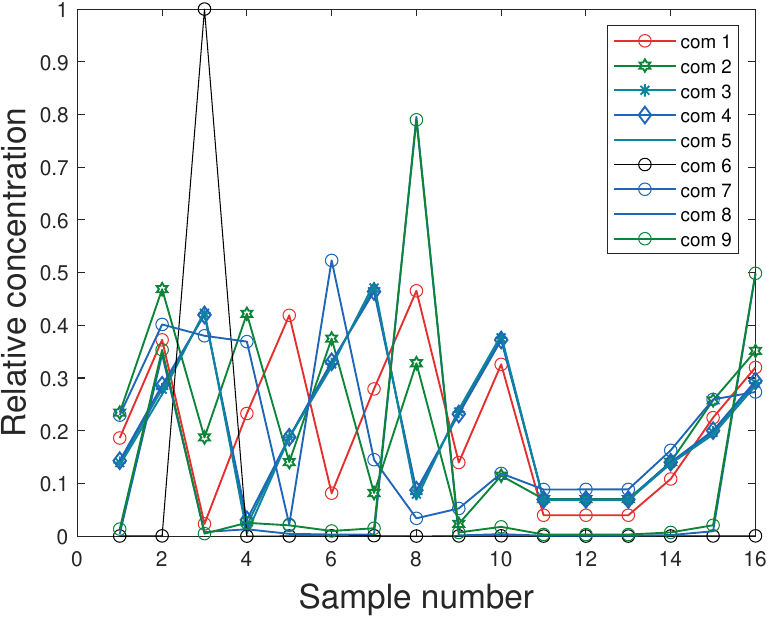}}\hspace{0cm}
	\subfloat[PARAFACM]{\includegraphics[width=4.2cm,height=3.6cm]{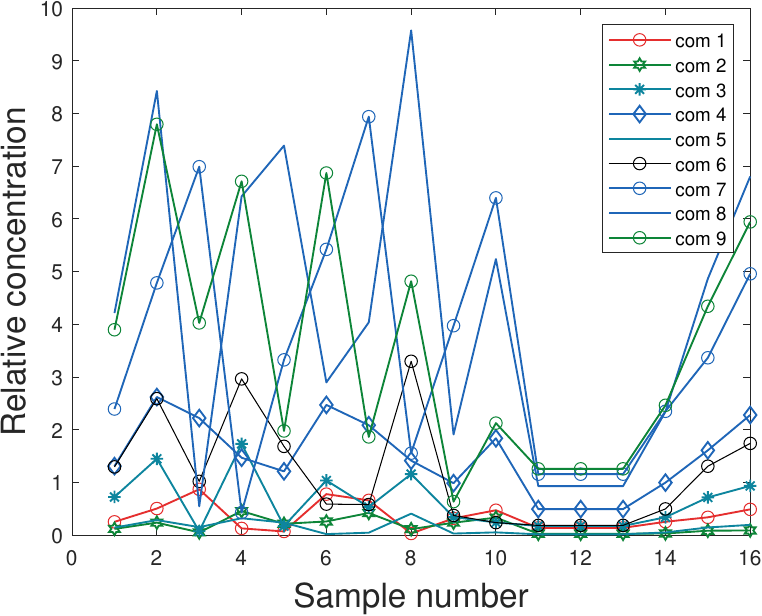}}\hspace{0cm}
	
	\subfloat[LRTCPR]{\includegraphics[width=4.2cm,height=3.6cm]{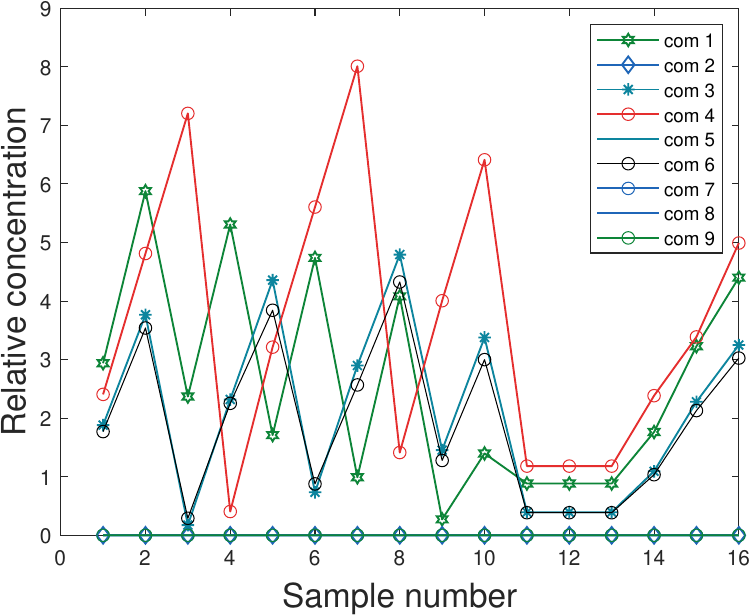}}
	\subfloat[BCPF]{ \includegraphics[width=4.2cm,height=3.6cm]{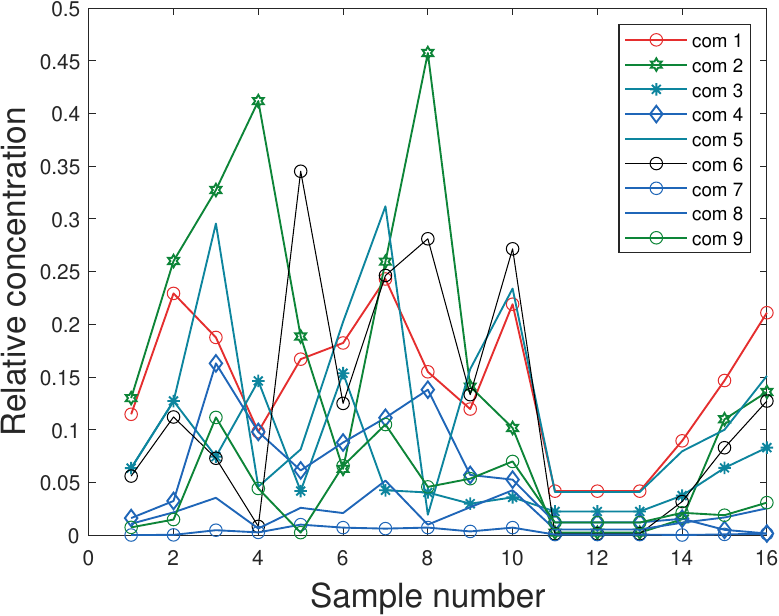}}\hspace{0cm}
	\subfloat[BCPGD]{\includegraphics[width=4.2cm,height=3.6cm]{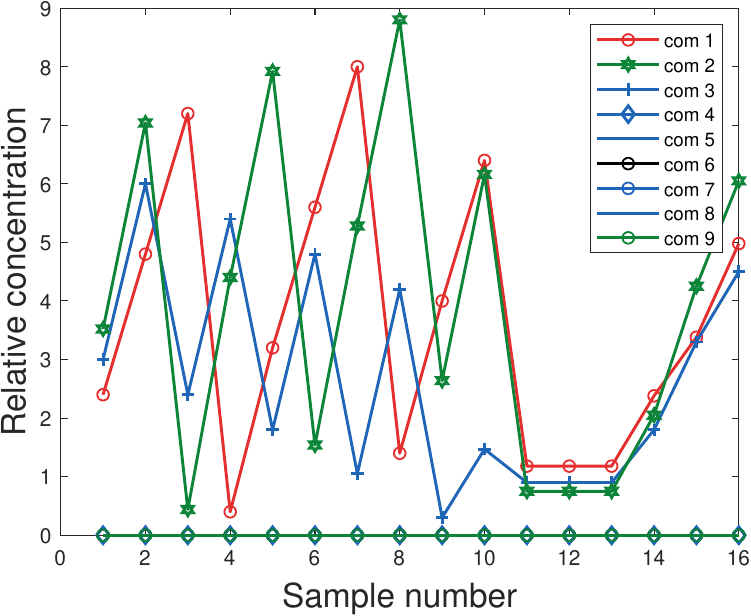}}
	\caption{ Analytical results with $R=9$ and $SR = 0.3 $ for the macrocephalae rhizoma data. The relative concentration profiles resolved by (a) Real, (b) INDAFAC, (c) PARAFACM, (d) LRTCPR (e)  BCPF and (f) BCPGD, respectively. }
	\label{fig3}
\end{figure}

\begin{figure}[!htbp]
	\centering
	\subfloat[Real]{\includegraphics[width=4.2cm,height=3.6cm]{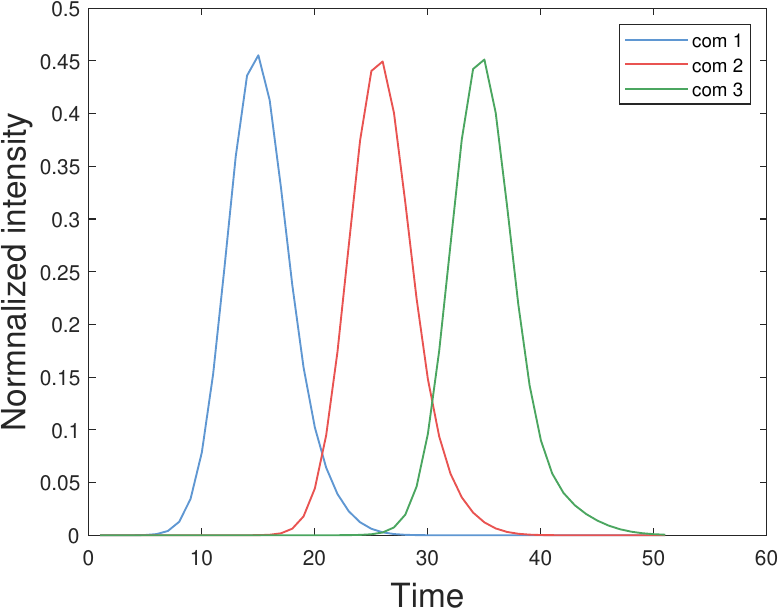}}
	\subfloat[INDAFAC]{\includegraphics[width=4.2cm,height=3.6cm]{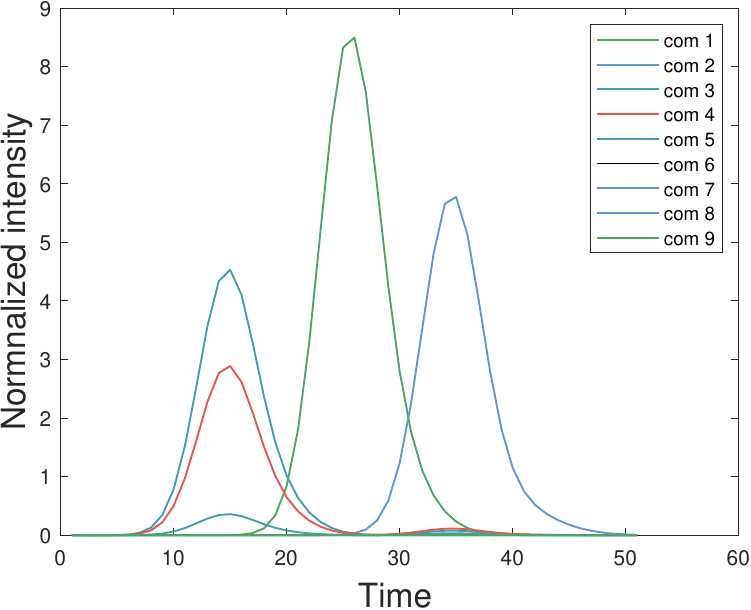}}
	\subfloat[PARAFACM]{\includegraphics[width=4.2cm,height=3.6cm]{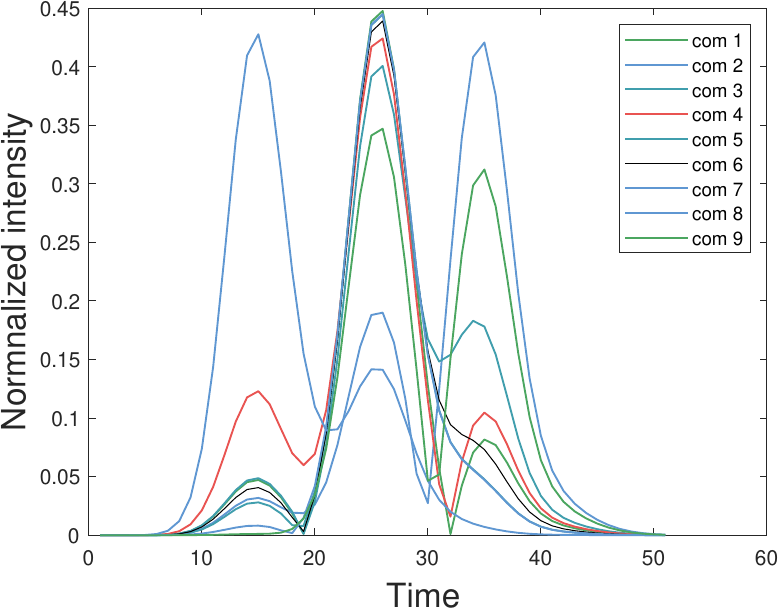}}
	
	\subfloat[LRTCPR]{\includegraphics[width=4.2cm,height=3.6cm]{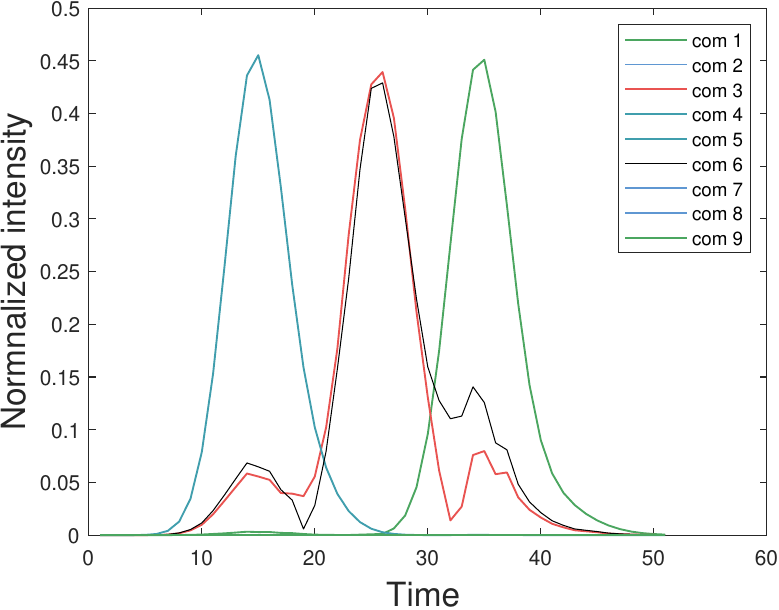}}
	\subfloat[BCPF]{\includegraphics[width=4.2cm,height=3.6cm]{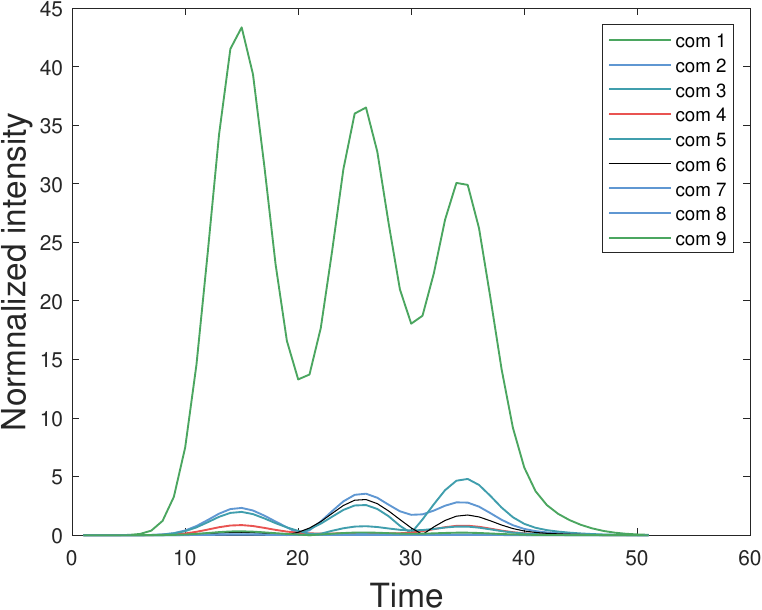}}
	\subfloat[BCPGD]{\includegraphics[width=4.2cm,height=3.6cm]{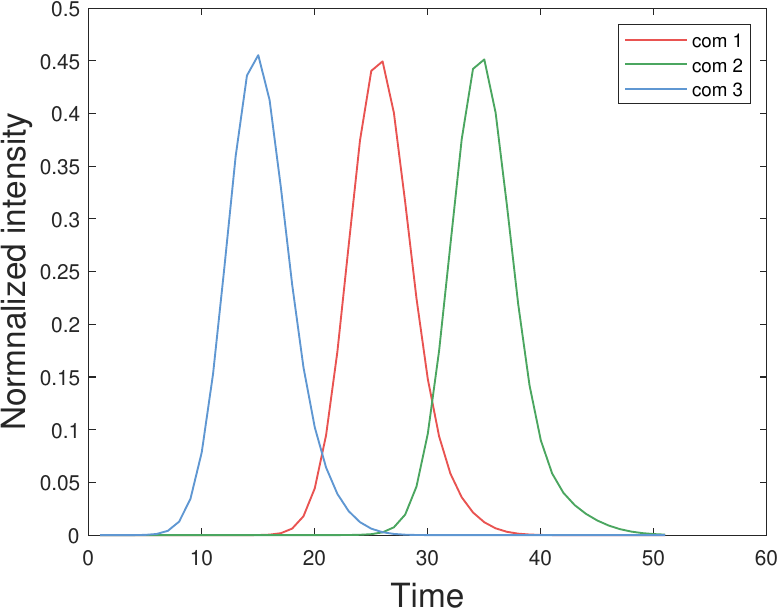}}
	
	\caption{Analytical results with $R=9$ and $SR=0.3$ for the macrocephalae rhizoma data. The normalized chromatographic profiles resolved by (a) Real, (b) INDAFAC, (c) PARAFACM, (d) LRTCPR ,(e)  BCPF and (f) BCPGD, respectively.}
	\label{fig4}
\end{figure}

\begin{figure}[!htbp]
	\centering
	\subfloat[Real]{\includegraphics[width=4.2cm,height=3.6cm]{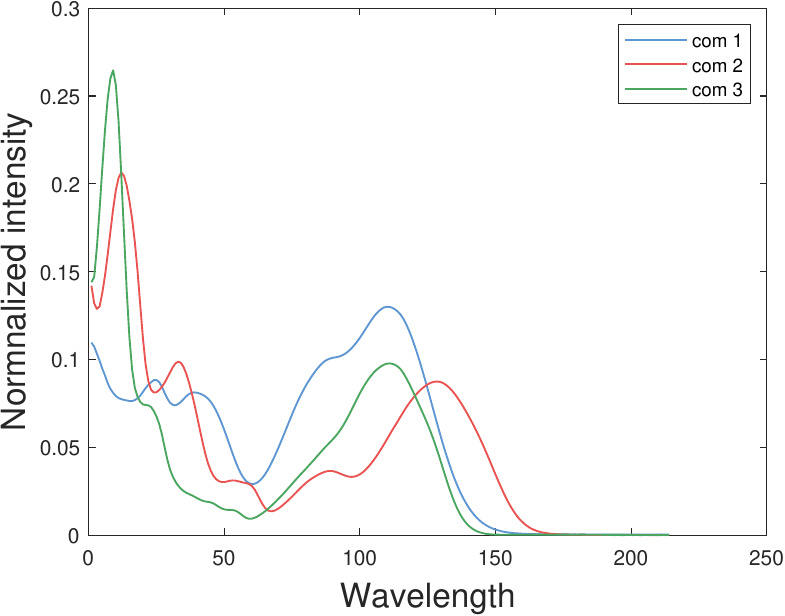}}
	\subfloat[INDAFAC]{\includegraphics[width=4.2cm,height=3.6cm]{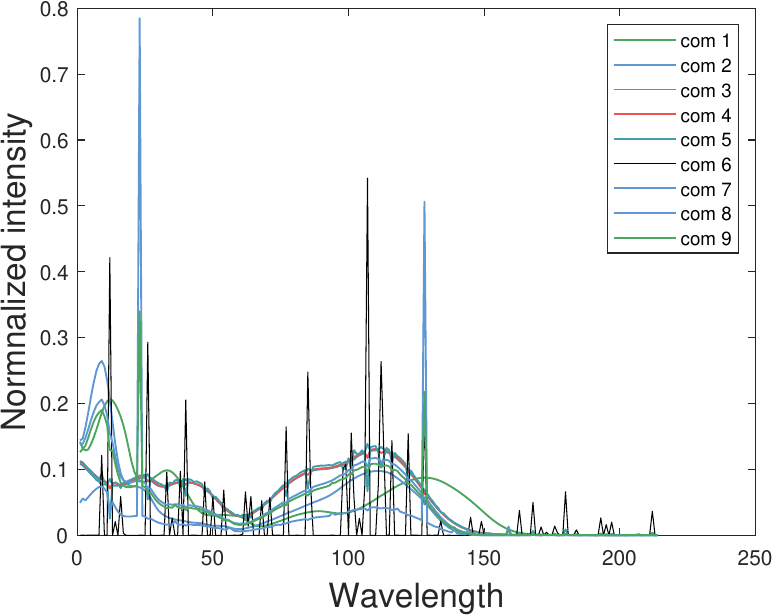}}
	\subfloat[PARAFACM]{\includegraphics[width=4.2cm,height=3.6cm]{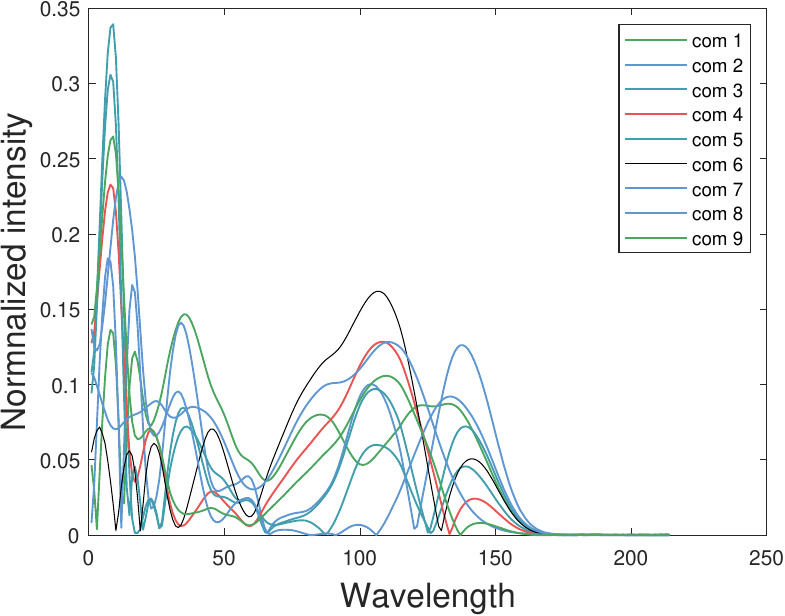}}
	
	\subfloat[LRTCPR]{ \includegraphics[width=4.2cm,height=3.6cm]{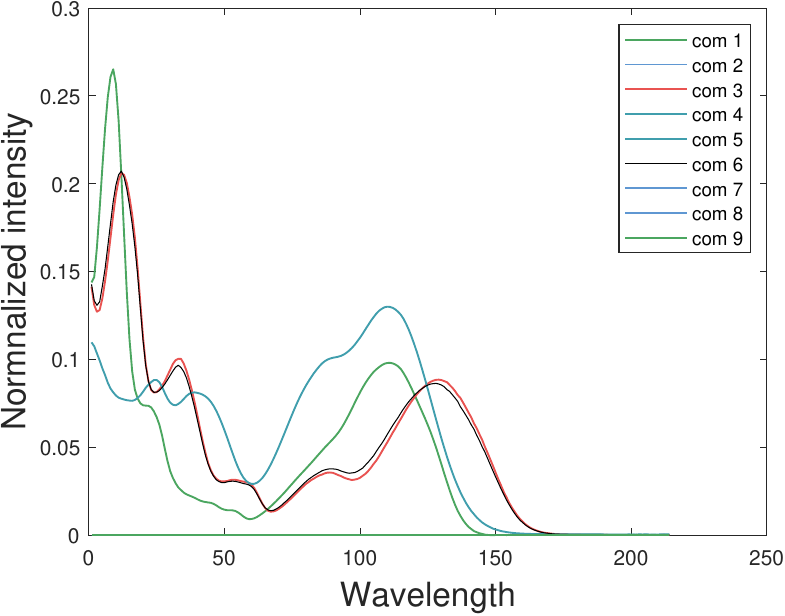}}
	\subfloat[BCPF]{\includegraphics[width=4.2cm,height=3.6cm]{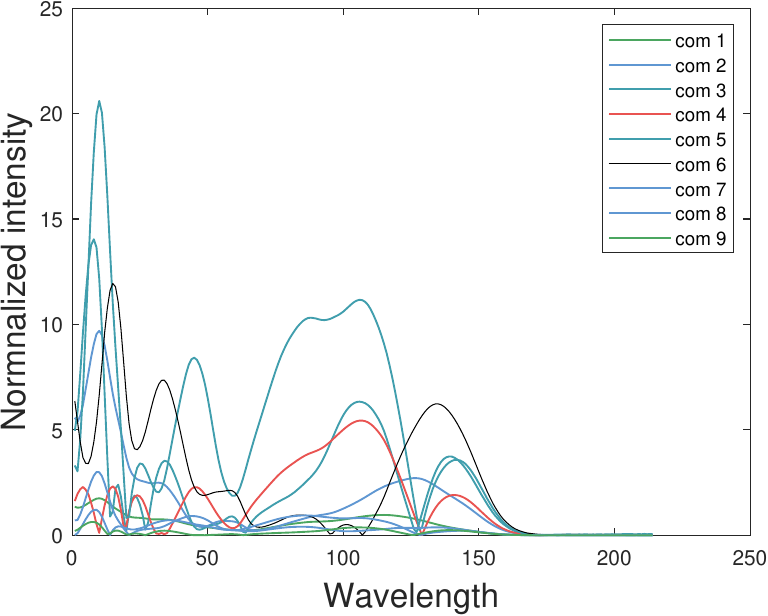}}
	\subfloat[BCPGD]{\includegraphics[width=4.2cm,height=3.6cm]{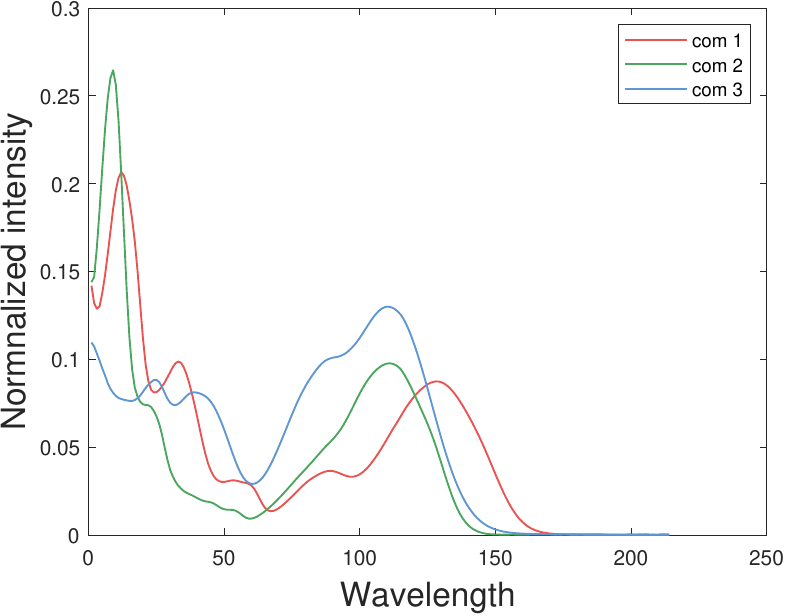}}
	\caption{Analytical results with $R=9$ and $SR=0.3 $ for the macrocephalae rhizoma data. The normalized spectra profiles resolved by (a) Real, (b) INDAFAC, (c) PARAFACM, (d) LRTCPR, (e) BCPF   and (f) BCPGD, respectively.} 
	\label{fig5}
\end{figure}

This data contains three components. We set the initial number of components to $ R = \{3,5,9,30\} $ under different sampling ratios $SR = \{0.5,0.3,0.1 \}$. The RMSEP values and computation times are reported in Table $\ref{table1}$.  Figure $ \ref{fig3} $ shows the real concentration profiles and the relative concentration profiles resolved by different methods with $R = 9$ and $SR = 0.3$, corresponding to C.
Figure \ref{fig4} presents the real chromatographic profiles and the normalized chromatographic profiles resolved by different methods with R = 9, corresponding to A. Figure \ref{fig5} shows the real spectra profiles and the normalized spectra profiles resolved by different methods with R = 9, corresponding to B.

As shown in Table \ref{table1}, only BCPGD can achieve satisfactory results when the initial number of components is overestimated across different sampling rates. Although LRTCPR can also achieve accurate component separation, it does so only when the assumed rank is slightly higher than the true rank. Furthermore, the performance of BCPGD is comparable to other methods when
the initial number of components is accurate.  These findings highlight the effectiveness and robustness of the proposed model and algorithm. As illustrated in Figure \ref{fig3}, \ref{fig4} and \ref{fig5}, the proposed method can recover accurate concentration profiles, chromatographic profiles, and spectra profiles with $R=9$ and  $SR=0.3$. 
\begin{table}[!htbp]
	\footnotesize
	\caption{ Comparison of root mean square error of prediction (RMSEP) of different algorithms and computation time (s) for the Macrocephalae rhizoma data in different sample radios $SR = \{0.5,0.3,0.1 \}$   with R = 3,5,9,30. 
		"-" shows that the components cannot be separated correctly. } 
	\label{table1}

	\scriptsize 
	\centering
	\begin{tabular}{|l|l|l|l|l|l|l|}
		\hline
		SR = 0.5 & Algorithm  & PARAFACM &  INDAFAC & LRTCPR &BCPF &   BCPGD\\ 
		\hline
		\multirow{1}*{R=3} 
		~ & RMSEP &  2.11e-07 &  \textbf{4.2848e-08} &1.82e-05 & 9.7675e-05  & 4.36e-07\\
		
		\hline
		\multirow{1}*{R=5}  
		~ & RMSEP & -    &-  &1.91e-05 & -  & \textbf{3.39e-07} \\
		\hline
		
		\multirow{1}*{R=9}  
		~ & RMSEP & -    &-  &- & -  & \textbf{8.03e-07} \\
		\hline

		\multirow{1}*{R=30}  
		~ & RMSEP & -   &-  &- & -& \textbf{3.41e-07}\\
		\hline
		
		\hline
		SR = 0.3  & Algorithm  & PARAFACM 
		& INDAFAC & LRTCPR &BCPF  & BCPGD\\
		\hline
		\multirow{1}*{R=3}  
		~ & RMSEP & 1.04e-06 & \textbf{3.47e-08} & 4.95e-05 & 1.4330e-04   & 5.99e-07 \\
		\hline
		\multirow{1}*{R=5} 
		~ & RMSEP & -    &-  &3.60e-05 & - & \textbf{7.10e-07}\\
		\hline
		
		\multirow{1}*{R=9}  
		~ & RMSEP & -    &-  &- & -  & \textbf{5.67e-07} \\
		\hline

		\multirow{1}*{R=30}  
		~ & RMSEP & -   &-  &- & - & \textbf{5.99e-07}\\
		\hline

		SR = 0.1 & Algorithm  & PARAFACM 
		& INDAFAC & LRTCPR & BCPF & BCPGD\\
		\hline
		\multirow{1}*{R=3}  
		~ & RMSEP &  2.61e-05 & \textbf{2.62e-08} & 1.09e-04 &0.0070  &1.95e-06\\
		\hline
		\multirow{1}*{R=5}  
		~ & RMSEP & -    &-  &1.22e-04&-  & \textbf{1.95e-06} \\
		\hline
		\multirow{1}*{R=9}  
		~ & RMSEP & -    &-  &- & -  & \textbf{3.00e-06} \\
		\hline
		\multirow{1}*{R=30}  
		~ & RMSEP & -   &-  &-& -&  \textbf{1.96e-06}\\
		\hline
	\end{tabular} 
\end{table}

\subsection{Experiments on randomly generated noisy low rank tensors}
We generate randomly a low rank tensor corrupted by Gaussian noise using the approach proposed, which has been widely adopted in chemometric studies. Let $\mathbf{A}, \mathbf{B}, \mathbf{C} $ be the factor matrices of $\mathcal{X} $.
\begin{equation}
	\mathbf{A} = \operatorname{rand}(I,R), \; \mathbf{B} = \operatorname{rand}(J,R), \;\mathbf{C} = \operatorname{rand}(K,R), \; \mathcal{\mathbf{X}}_{..k} =\mathbf{A} \text{diag}(\mathbf{c}_{k}^{T}) \mathbf{B}^T.
\end{equation}
Then we denote the error tensor and noise level by $E$ and $a $, respectively, and use $\hat{\mathcal{X}} $ to represent the resulting noisy data.
\begin{equation}
	\begin{array}{rl}
		& \mathbf{E}_{..k} = \max(\mathbf{x}_{..k}) \times randn(I,J) \times a,  \; k=1,\dots,K, \\
		& \hat{\mathcal{\mathbf{X}}}_{..k} = \mathcal{\mathbf{X}}_{..k} + \mathbf{E}_{..k},  \; k=1,\dots,K.
	\end{array}
	\nonumber
\end{equation}
where $rand $ generates a random matrix with entries uniformly distributed on \((0,1)\), whereas $randn $ generates a random matrix whose entries follow a standard normal distribution with zero mean and unit variance.

We generate a synthetic tensor \(\mathcal{X} \in \mathbb{R}^{300 \times 200 \times 20}\) with six underlying components. The noise level \(a\) is varied over \(\{0.01, 0.05, 0.1\}\), while the initial number of components $R $ is set to 30.

\begin{table}[!htbp]
	\footnotesize
	\caption{The performance of Algorithm \ref{algorithm1} on randomly generated low rank tensors corrupted by Gaussian noise.}
	\label{noise}
	\centering
	\begin{tabular}{|c|c|c|c|c|c|}
		\hline
		$a$ & SR & Initial-R & True-R & Compute-R & RMSEP \\
		\hline
		
		\multirow{4}{*}{$0.1$}
		& 1 & 30 & 6 & 6 & 0.0020 \\
		\hhline{|~||-|-|-|-|-|}
		& 0.8 & 30 & 6 & 6 & 0.0034 \\
		\hhline{|~||-|-|-|-|-|}
		& 0.5 & 30 & 6 & 6 & 0.0042 \\
		\hhline{|~||-|-|-|-|-|}
		& 0.3 & 30 & 6 & 6 & 0.0048 \\
		\hline
		
		\multirow{4}{*}{$0.05$}
		& 1 & 30 & 6 & 6 & 0.0013 \\
		\hhline{|~||-|-|-|-|-|}
		& 0.8 & 30 & 6 & 6 & 0.0017 \\
		\hhline{|~||-|-|-|-|-|}
		& 0.5 & 30 & 6 & 6 & 0.0019 \\
		\hhline{|~||-|-|-|-|-|}
		& 0.3 & 30 & 6 & 6 & 0.0021 \\
		\hline
		
		\multirow{4}{*}{$0.01$}
		& 1 & 30 & 6 & 6 & 2.4711e-04 \\
		\hhline{|~||-|-|-|-|-|}
		& 0.8 & 30 & 6 & 6 & 2.6417e-04 \\
		\hhline{|~||-|-|-|-|-|}
		& 0.5 & 30 & 6 & 6 & 3.5032e-04 \\
		\hhline{|~||-|-|-|-|-|}
		& 0.3 & 30 & 6 & 6 & 4.5692e-04 \\
		\hline
		
	\end{tabular}
\end{table}
Table $\ref{noise}$ shows that the proposed model and algorithm can successfully calculate the true rank of the randomly generated low rank tensors corrupted by Gaussian noise, which indicates their effectiveness and robustness.

\section{Conclusions}
In this paper, we developed a group sparse-based CP decomposition model for incomplete tensors, which achieves rank estimation, CP decomposition, and tensor completion simultaneously. We also established an error bound between the optimal solution  and the ground truth of the model under the true CP rank is unknown. Furthermore, we established a theory of adaptive parameter selection and obtained a closed-form expression of regularization parameter.
To solve this model, we employed a  block-coordinate proximal gradient descent algorithm (BCPGD) with extrapolation and provided the convergence analysis. Finally, we evaluated the performance of our proposed model and algorithm by using the randomly generated low rank tensors and real-word data in the field of chemometrics. Numerical experiments by showed that the Algorithm $\ref{algorithm1}$ can estimate the true rank and achieve CP decomposition for incomplete tensors. It is of interest whether the global recovery theory can be established for the proposed model. We leave this as a topic for future
research.

\begin{appendices}

\section{The proof of Lemma 6.3} 
\label{leamm_6_3}

\textit{Proof}. \textbf{(i)} At the $k$-th iteration of Algorithm $\ref{algorithm1} $, if the condition of step 10 is is not satisfied, then $(\ref{decrease_inequality})  $  holds automatically. Conversely, if the condition of Step 10 is satisfied, we can obtain the following inequalities from a basic descent property of proximal gradient steps 
\begin{equation}
	\begin{array}{rcl}
		
		&  F(\mathbf{X}^{k,u},\mathbf{Y}^{k},\mathbf{Z}^{k})   \geq  \frac{\epsilon}{2} \Vert \mathbf{X}^{k,u} - \mathbf{X}^{k,u+1}  \Vert^2 + F(\mathbf{X}^{k,u+1},\mathbf{Y}^{k},\mathbf{Z}^{k}),  \\
		
		& F(\mathbf{X}^{k+1},\mathbf{Y}^{k,u},\mathbf{Z}^{k})   \geq   \frac{\epsilon}{2} \Vert \mathbf{Y}^{k,u} - Y^{k,u+1}  \Vert^2 + F(\mathbf{X}^{k+1},\mathbf{Y}^{k,u+1},\mathbf{Z}^{k}), \\
		
		& F(\mathbf{X}^{k+1},\mathbf{Y}^{k+1},\mathbf{Z}^{k,u})   \geq  \frac{\epsilon}{2} \Vert \mathbf{Z}^{k,u} - \mathbf{Z}^{k,u+1}  \Vert^2 + F(\mathbf{X}^{k+1},\mathbf{Y}^{k+1},\mathbf{Z}^{k,u+1}). \\ 
	\end{array}
	\label{13}
\end{equation}
For each inequality in (\ref{13}), adding from $u=0$ to $u=m$ and utilizing $\sum_{i=1}^{m} (a_i)^2 \geq \frac{(\sum_{i=1}^{m} a_i)^2}{m}$, we have
\begin{equation}
	\begin{array}{rl}
		F(\mathbf{X}^{k},\mathbf{Y}^{k},\mathbf{Z}^{k})   &\geq  \frac{\epsilon}{2}( \sum_{u=0}^m  \Vert \mathbf{X}^{k,u} - \mathbf{X}^{k,u+1}  \Vert^2) + F(\mathbf{X}^{k+1},\mathbf{Y}^{k},\mathbf{Z}^{k})  \\
		&\geq  \frac{\epsilon}{2m}  \left( \sum_{u=0}^m \Vert \mathbf{X}^{k,u} - \mathbf{X}^{k,u+1}  \Vert \right)^2 + F(\mathbf{X}^{k+1},\mathbf{Y}^{k},\mathbf{Z}^{k})   \\
		&\geq   \frac{\epsilon}{2m} \Vert \sum_{u=0}^m (\mathbf{X}^{k,u} - \mathbf{X}^{k,u+1})  \Vert^2 + F(\mathbf{X}^{k+1},\mathbf{Y}^{k},\mathbf{X}^{k}) \\
		
	\end{array}
	\label{14}
\end{equation}
where $\mathbf{X}^{k,0} = \mathbf{X}^{k}   $ and $\mathbf{X}^{k+1} = \mathbf{X}^{k,m_{}+1}  $. Likewise, for the calculation of $\mathbf{Y}^k$ and $ \mathbf{Z}^k$, we obtain 
\begin{equation}
	\begin{array}{cc}
		& F(\mathbf{X}^{k+1},\mathbf{Y}^{k},\mathbf{Z}^{k}) \geq  \frac{\epsilon}{2m} \Vert \mathbf{Y}^{k} - \mathbf{Y}^{k+1}  \Vert^2 + F(\mathbf{X}^{k+1},\mathbf{Y}^{k+1},\mathbf{Z}^{k}),  \\
		& F(\mathbf{X}^{k+1},\mathbf{Y}^{k+1},\mathbf{Z}^{k}) \geq  \frac{\epsilon}{2m} \Vert \mathbf{Z}^{k} - \mathbf{Z}^{k+1}  \Vert^2 + F(\mathbf{X}^{k+1},\mathbf{Y}^{k+1},\mathbf{Z}^{k+1}).               \\
		
	\end{array}
	\label{15}
\end{equation}
Combining $(\ref{14})$ and $(\ref{15})$ and letting $\rho = \frac{\epsilon}{m} $, we have 
\begin{equation*}
	\begin{array}{lll}
		F( \! \mathbf{X}^{k} \!, \! \mathbf{Y}^{k} \!, \!\mathbf{Z}^{k} \! ) \! \geq \! \frac{\rho}{2} (\Vert  \mathbf{X}^{k} \! \! - \! \! \mathbf{X}^{k+1} \!  \Vert^2 \! + \! \! \Vert \mathbf{Y}^{k} \!\! -\! \! \mathbf{Y}^{k+1} \!  \Vert^2 \! +\! \!  \Vert \mathbf{Z}^{k} \! \! - \! \!  \mathbf{Z}^{k+1} \!  \Vert^2 )  
		\! +  \! \! F( \mathbf{X}^{k+1} \!, \!\mathbf{Y}^{k+1} \!, \! \mathbf{Z}^{k+1} \!),
	\end{array}
	\label{16}
\end{equation*} 
which implies  $(\ref{decrease_inequality})$ holds. 

\textbf{(ii)} Since $F(\mathbf{W}^k)$ is non-increasing and bounded from below, it follows from (\ref{decrease_inequality}) that $\sum_{k=1}^{+\infty} \Vert \mathbf{W}^k - \mathbf{W}^{k+1}  \Vert^2 < +\infty$, which implies $\lim\limits_{k \to +\infty} \Vert \mathbf{W}^k - \mathbf{W}^{k+1}  \Vert =0 $. From the relationship of Step 4, 6 and 8 of Algorithm $\ref{algorithm1}$, it can be deduced that
\begin{equation}
	\begin{array}{rcl}
		& \bar{\mathbf{X}}^{k+1} - {\mathbf{X}}^{k+1} = \alpha_1 (\bar{\mathbf{X}}^{k} - {\mathbf{X}}^{k}) + \alpha_1 (\mathbf{X}^{k+1} - \mathbf{X}^{k}), \\
		& \bar{\mathbf{Y}}^{k+1} - \mathbf{Y}^{k+1} = \alpha_2 (\bar{\mathbf{Y}}^{k} - \mathbf{Y}^{k}) + \alpha_2 (\mathbf{Y}^{k+1} - \mathbf{Y}^{k}), \\
		& \bar{\mathbf{Z}}^{k+1} - {\mathbf{Z}}^{k+1} = \alpha_3 (\bar{\mathbf{Z}}^{k} - {\mathbf{Z}}^{k}) + \alpha_3 (\mathbf{Z}^{k+1} - \mathbf{Z}^{k}). \\
		
	\end{array}
	\label{17}
\end{equation}
Therefore, one has
\begin{equation*}
	\begin{array}{rl}
		\Vert \mathbf{X}^{k+1} - \mathbf{X}^{k} \Vert^2 \! \! \! \! \! \!   & =  \Vert (\bar{\mathbf{X}}^{k} - \mathbf{X}^{k}) 
		- \frac{1}{\alpha_1}(\bar{\mathbf{X}}^{k+1} - \mathbf{X}^{k+1}) \Vert^2  \\
		& = \frac{1}{\alpha_1^2} \Vert \bar{\mathbf{X}}^{k+1} - \mathbf{X}^{k+1}   \Vert^2
		+ \Vert \bar{\mathbf{X}}^{k} -\mathbf{X}^{k}  \Vert^2 -\frac{2}{\alpha_1} \langle
		\bar{\mathbf{X}}^{k} \!-\! \mathbf{X}^{k}, \bar{\mathbf{X}}^{k+1}  \!-\!  \mathbf{X}^{k+1}  \rangle  \\
		& \geq  (1-\frac{1}{\alpha_1})\Vert \bar{\mathbf{X}}^{k} - \mathbf{X}^{k} \Vert^2 + (\frac{1}{\alpha_1^2} - \frac{1}{\alpha_1} )\Vert \bar{\mathbf{X}}^{k+1} - \mathbf{X}^{k+1} \Vert^2.  \\  
	\end{array}
	\label{18}
\end{equation*}
Summing from $k=0$ to $k=+\infty$,  it can be inferred that 
\vspace{-6pt}
\begin{equation}
	\begin{array}{rl}
		+\infty &> \sum_{k=0}^{+\infty} \Vert \mathbf{X}^{k+1} - \mathbf{X}^{k}  \Vert^2 \\
		& \geq  (1-\frac{2}{\alpha_1}+\frac{1}{\alpha_1^2}) \sum_{k=1}^{+\infty} \Vert \bar{\mathbf{X}}^{k} - \mathbf{X}^{k}  \Vert^2 + (1-\frac{1}{\alpha_1}) \Vert \bar{\mathbf{X}}^{0} - \mathbf{X}^{0}  \Vert^2 \\
		& =  (1-\frac{1}{\alpha_1})^2 \sum_{k=1}^{+\infty} \Vert \bar{\mathbf{X}}^{k} - \mathbf{X}^{k}  \Vert^2.   \\
	\end{array}
	\label{19}
\end{equation}
Parallelly, we have the following inequalities for $\mathbf{Y}$ and $\mathbf{Z}$
\begin{equation}
	\begin{array}{rl}
		+\infty &> (1-\frac{1}{\alpha_2})^2 \sum_{k=1}^{+\infty} \Vert \bar{\mathbf{Y}}^{k} - \mathbf{Y}^{k}  \Vert^2,  \\
		+\infty &> (1-\frac{1}{\alpha_3})^2 \sum_{k=1}^{+\infty} \Vert \bar{\mathbf{Z}}^{k} - \mathbf{Z}^{k}  \Vert^2.  \\ 
	\end{array}
	\label{20}
\end{equation}
(\ref{19}) and (\ref{20}) imply $\sum_{k=1}^{+\infty} \Vert \bar{\mathbf{W}}^{k} - \mathbf{W}^{k}  \Vert^2 <+ \infty$, so that $\lim\limits_{k \to +\infty} \Vert \bar{\mathbf{W}}^{k} - \mathbf{W}^{k}  \Vert=0.$ 
This completes this proof.    \hfill  $\square   $

\section{The proof of Theorem 6.4} 
\label{theorem6_4}

\textit{Proof}. Before that, we present a subgradient lower bound for the iteration gap. Suppose that Step 10 is skipped. Using the Fermat's rule to the optimization problem in Step 4, 6, 8 and applying [\cite{rockafellar1998variational},Proposition 10.5  10.10(subgradients of Lipschitzian sums)], for any $1\leq j \leq R$,  there exist three vectors $ \bm{\mu}_1^j,\bm{\mu}_2^j,\bm{\mu}_3^j$ belong to $ \partial r_1^j(\mathbf{X}^{k+1}), r_2^j(\mathbf{Y}^{k+1}) $ and $r_3^j(\mathbf{Z}^{k+1})$ respectively such that 
\begin{equation*}
	\begin{array}{cc}
		\nabla_{\mathbf{x}_j} f(\mathbf{x}_{<j}^{k,m_1+1},\mathbf{x}_{\geq j}^{k,m_1},\bar{\mathbf{Y}}^{k},\bar{\mathbf{Z}}^{k}) + (l_{x,j}^{k,m_1}+\epsilon)(\mathbf{x}_j^{k,m_1+1} -\mathbf{x}_j^{k,m_1} )+ \bm{\mu}_1^j=\mathbf{0}, \\
		
		\nabla_{\mathbf{y}_j} f(\bar{\mathbf{X}}^{k+1},\mathbf{y}_{<j}^{k,m_2+1},\mathbf{y}_{\geq j}^{k,m_2},\bar{\mathbf{Z}}^{k}) + (l_{y,j}^{k,m_2}+\epsilon)(\mathbf{y}_j^{k,m_2+1} -\mathbf{y}_j^{k,m_2} )+ \bm{\mu}_2^j=\mathbf{0},   \\
		\nabla_{\mathbf{z}_j} f(\bar{\mathbf{X}}^{k+1},\bar{\mathbf{Y}}^{k+1},\mathbf{z}_{<j}^{k,m_3+1},\mathbf{z}_{\geq j}^{k,m_3}) + (l_{z,j}^{k,m_3}+\epsilon)(\mathbf{z}_j^{k,m_3+1} -\mathbf{z}_j^{k,m_3} )+ \bm{\mu}_3^j=\mathbf{0}.  \\
	\end{array}
	\label{21}
\end{equation*}
Note that applying the Lipschitz continuity of the partial gradient, we have
\begin{equation*}
	\begin{array}{rl}
		\Vert \nabla_{\mathbf{x}_j} f(\mathbf{x}_{<j}^{k,m_1+1},\mathbf{x}_{\geq j}^{k,m_1},\bar{\mathbf{Y}}^k,\bar{\mathbf{Z}}^k) - \nabla_{\mathbf{x}_j} f(\mathbf{X}^{k+1},\bar{\mathbf{Y}}^k,\bar{\mathbf{Z}}^k)  \Vert   
		\! \!\!\!   &\leq l_{x}^+ \Vert (\mathbf{X}^{k+1}\! - \!\mathbf{X}^{k,m_1})_{ \geq j}  \Vert \\
		&\leq l_{x}^+ \Vert \mathbf{X}^{k+1}\! - \!\mathbf{X}^{k,m_1}   \Vert ,  \\
	\end{array}
	\label{22}
\end{equation*}
\begin{equation*}
	\begin{array}{ll}
		\Vert \nabla_{\mathbf{y}_j} f(\bar{\mathbf{X}}^{k+1},\mathbf{y}_{<j}^{k,m_2+1},\mathbf{y}_{\geq j}^{k,m_2},\bar{\mathbf{Z}}^k)  \! -  \! \nabla_{\mathbf{y}_j} f(\bar{\mathbf{X}}^{k+1},{\mathbf{Y}}^{k+1},\bar{\mathbf{Z}}^k)  \Vert \! \! \leq   l_y^+ \Vert \mathbf{Y}^{k+1} \! - \!\mathbf{Y}^{k,m_2}   \Vert ,  \\
		
		\Vert \nabla_{\mathbf{z}_j} f( \bar{\mathbf{X}}^{k+1}, \! \bar{\mathbf{Y}}^{k+1}, \!\mathbf{z}_{<j}^{k,m_3+1},\!\mathbf{z}_{\geq j}^{k,m_3})  \!-  \! \nabla_{\mathbf{z}_j} f(\bar{\mathbf{X}}^{k+1},\!\bar{\mathbf{Y}}^{k+1},\!\bar{\mathbf{Z}}^{k+1})  \Vert \! \!  \leq  l_z^+  \Vert \mathbf{Z}^{k+1}\! - \!\mathbf{Z}^{k,m_2}   \Vert  . \\
		
	\end{array}
	\label{23}
\end{equation*}
By combining the preceding three inequalities, we obtain
\begin{equation}
	\begin{array}{rl}
		& \Vert  \nabla_{\mathbf{x}_j} f(\mathbf{X}^{k+1},\bar{\mathbf{Y}}^{k},\bar{\mathbf{Z}}^{k}) \!+\! \mu_1^j \Vert  \\
		&= 
		\Vert  \nabla_{\mathbf{x}_{j}} f(\mathbf{X}^{k+1},\bar{\mathbf{Y}}^{k},\bar{\mathbf{Z}}^{k})\! - \! \nabla_{\mathbf{x}_j} f(\mathbf{x}_{<j}^{k,m_1+1},\mathbf{x}_{\geq j}^{k,m_1},\bar{\mathbf{Y}}^{k},\bar{\mathbf{Z}}^{k}) \\
		&- (l_{x,j}^{k,m_1}+\epsilon)(\mathbf{x}_j^{k,m_1+1} -\mathbf{x}_j^{k,m_1} )  \Vert \\
		&\leq \Vert  \nabla_{\mathbf{x}_j} f(\mathbf{X}^{k+1},\bar{\mathbf{Y}}^{k},\bar{\mathbf{Z}}^{k}) 
		-  \nabla_{\mathbf{x}_j} f(\mathbf{x}_{<j}^{k,m_1+1},\mathbf{x}_{\geq j}^{k,m_1},\bar{\mathbf{Y}}^{k},\bar{\mathbf{Z}}^{k})    \Vert   \\
		& + \Vert (l_{x,j}^{k,m_1}+\epsilon)(\mathbf{x}_j^{k,m_1+1} -\mathbf{x}_j^{k,m_1} )   \Vert   \\
		&\leq   l_x^+ \Vert \mathbf{X}^{k+1} - \mathbf{X}^{k,m_1} \Vert + (l_x^+ + \epsilon) \Vert \mathbf{X}^{k+1} - \mathbf{X}^{k,m_1} \Vert  \\
		&= (2 l_x^+ + \epsilon) \Vert \mathbf{X}^{k+1} - \mathbf{X}^{k,m_1}  \Vert.
	\end{array}
	\label{24}
\end{equation}
We have similar results for $\mathbf{Y}$ and $\mathbf{Z}$
\begin{equation}
	\begin{array}{cc}
		\Vert  \nabla_{\mathbf{y}_j} f(\bar{\mathbf{X}}^{k+1},{\mathbf{Y}}^{k+1},\bar{\mathbf{Z}}^{k}) + \bm{\mu}_2^j \Vert \leq (2 l_y^+ + \epsilon) \Vert \mathbf{Y}^{k+1} - \mathbf{Y}^{k,m_2} \Vert, \\
		
		\Vert  \nabla_{\mathbf{z}_j} f(\bar{\mathbf{X}}^{k+1},\bar{\mathbf{Y}}^{k+1},{\mathbf{Z}}^{k+1}) + \bm{\mu}_3^j \Vert \leq (2 l_z^+ + \epsilon) \Vert \mathbf{Z}^{k+1} - \mathbf{Z}^{k,m_3} \Vert. \\
		
	\end{array}
	\label{25}
\end{equation}
Therefore, from inequalities (\ref{24}) and (\ref{25}), there exist $\bm{\mu}_i = (\bm{\mu}_1^i,\dots, \bm{\mu}_i^R),i=1,2,3, \bm{\mu}_1 \in   
\partial r_1(\mathbf{X}^{k+1}), \bm{\mu}_2 \in \partial r_2(\mathbf{Y}^{k+1}) $, $\bm{\mu}_3\in \partial r_3(\mathbf{Z}^{k+1}) $ such that
\begin{equation}
	\begin{array}{rl}
		\Vert  \nabla_X f(\mathbf{X}^{k+1},\bar{\mathbf{Y}}^{k},\bar{\mathbf{Z}}^{k}) + \bm{\mu}_1 \Vert 
		&\leq \sum_{j=1}^R \Vert \nabla_{x_j}  
		f(\mathbf{X}^{k+1},\bar{\mathbf{Y}}^{k},\bar{\mathbf{Z}}^{k}) + \bm{\mu}_1^j  \Vert  \\
		&\leq  R (2 l_x^+ + \epsilon) \Vert \mathbf{X}^{k+1} - \mathbf{X}^{k,m_1} \Vert \\
		&= L_x \Vert \mathbf{X}^{k+1} - \mathbf{X}^{k,m_1} \Vert,   \\
		
		\Vert  \nabla_\mathbf{Y} f(\bar{\mathbf{X}}^{k+1},\mathbf{Y}^{k+1},\bar{\mathbf{X}}^{k}) + \bm{\mu}_2 \Vert &\leq L_y \Vert \mathbf{Y}^{k+1} - \mathbf{Y}^{k,m_2} \Vert , \\
		
		\Vert  \nabla_\mathbf{Z} f(\bar{\mathbf{X}}^{k+1},\bar{\mathbf{Y}}^{k+1},\mathbf{Z}^{k+1}) + \bm{\mu}_3 \Vert &\leq L_z \Vert \mathbf{Z}^{k+1} - \mathbf{Z}^{k,m_3} \Vert.  \\
		
	\end{array}
	\label{26}
\end{equation}
where $\mathbf{X}^{k,m_1} = (\mathbf{x}^{k,m_1}_1,\dots, \mathbf{x}^{k,m_1}_R), \mathbf{Y}^{k,m_2} = (\mathbf{y}^{k,m_2}_1,\dots, \mathbf{y}^{k,m_2}_R),\mathbf{Z}^{k,m_3} = (\mathbf{z}^{k,m_3}_1,\dots,
\mathbf{z}^{k,m_3}_R), 
L_i=R(2l_i^+ + \epsilon), i=x,y,z$.
Utilizing (\ref{26}) and assumption $\ref{a_4}$(iv),  let $ \alpha = \operatorname{min}\{ \alpha_1, \alpha_2, \alpha_3\} $, we get 
\begin{equation}
	\begin{array}{rl}
		&\Vert  \nabla_\mathbf{X} f(\mathbf{W}^{k+1}) + \bm{\mu}_1 \Vert   \\
		& = \Vert \nabla_\mathbf{X} f(\mathbf{W}^{k+1}) - \nabla_\mathbf{X} f(\mathbf{X}^{k+1},\bar{\mathbf{Y}}^{k},\bar{\mathbf{Z}}^{k}) +  \nabla_\mathbf{X} f(\mathbf{X}^{k+1},\bar{\mathbf{Y}}^{k},\bar{\mathbf{Z}}^{k}) +  \bm{\mu}_1  \Vert  \\
		&\leq \Vert \nabla_\mathbf{X} f(\mathbf{W}^{k+1}) - \nabla_\mathbf{X} f(\mathbf{X}^{k+1},\bar{\mathbf{Y}}^{k},\bar{\mathbf{Z}}^{k})   \Vert + \Vert \nabla_\mathbf{X} f(\mathbf{X}^{k+1},\bar{\mathbf{Y}}^{k},\bar{\mathbf{Z}}^{k}) +  \bm{\mu}_1   \Vert   \\
		&\leq L_x \Vert \mathbf{X}^{k+1} - \mathbf{X}^{k,m_1} \Vert + \frac{M}{\alpha} \Vert \bar{\mathbf{W}}^{k+1} - \mathbf{W}^{k+1} \Vert.
	\end{array}
	\label{27}
\end{equation}
In a similar way, we obtain
\begin{equation}
	\begin{array}{rl}
		\Vert  \nabla_Y f(\mathbf{W}^{k+1}) +  \bm{\mu}_2 \Vert &\leq  L_y \Vert \mathbf{Y}^{k+1} - \mathbf{Y}^{k,m_2} \Vert + \frac{M}{\alpha} \Vert \bar{\mathbf{W}}^{k+1} - \mathbf{W}^{k+1} \Vert, \\
		\Vert  \nabla_Z f(\mathbf{W}^{k+1}) +  \bm{\mu}_3 \Vert &\leq L_z \Vert \mathbf{Z}^{k+1} - \mathbf{Z}^{k,m_3} \Vert + \frac{M}{\alpha} \Vert \bar{\mathbf{W}}^{k+1} - \mathbf{W}^{k+1} \Vert. \\      
	\end{array}
	\label{28}
\end{equation}
In the reverse case, if Step 10 is executed, by repeating the above procedure to the optimization problem in Step 11, 13, 15, we have there exist $ \bm{\mu}_1 \in \partial r_1(\mathbf{X}^{k+1})$, $ \bm{\mu}_2 \in \partial r_2(\mathbf{Y}^{k+1})$ and $ \bm{\mu}_3 \in \partial r_3(\mathbf{Z}^{k+1}) $ such that
\begin{equation}
	\begin{array}{rcl}
		&\Vert  \nabla_\mathbf{X} f(\mathbf{X}^{k+1},\mathbf{Y}^{k},\mathbf{Z}^{k}) +  \bm{\mu}_1 \Vert
		\leq L_x \Vert \mathbf{X}^{k+1} - \mathbf{X}^{k,m_1} \Vert, \\
		
		&\Vert  \nabla_\mathbf{Y} f(\mathbf{X}^{k+1},\mathbf{Y}^{k+1},\mathbf{Z}^{k}) +  \bm{\mu}_2 \Vert \leq L_y \Vert \mathbf{Y}^{k+1} - \mathbf{Y}^{k,m_2} \Vert , \\
		
		&\Vert  \nabla_\mathbf{Z} f(\mathbf{X}^{k+1},\mathbf{Y}^{k+1},\mathbf{Z}^{k+1}) +  \bm{\mu}_3 \Vert \leq L_z \Vert \mathbf{Z}^{k+1} - \mathbf{Z}^{k,m_3} \Vert.   
	\end{array}
	\label{29}
\end{equation}
Therefore, 
\begin{equation}
	\begin{array}{rl}
		& \Vert \nabla_\mathbf{X} f(\mathbf{W}^{k+1}) +  \bm{\mu}_1  \Vert  \\
		& = \Vert \nabla_\mathbf{X} f(\mathbf{W}^{k+1}) - \nabla_\mathbf{X} f(\mathbf{X}^{k+1},\mathbf{Y}^{k},\mathbf{Z}^{k}) +  \nabla_\mathbf{X} f(\mathbf{X}^{k+1},\mathbf{Y}^{k},\mathbf{Z}^{k}) +   \bm{\mu}_1  \Vert  \\
		&\leq \Vert \nabla_X f(\mathbf{W}^{k+1}) - \nabla_\mathbf{X} f(\mathbf{X}^{k+1},\mathbf{Y}^{k},\mathbf{Z}^{k})   \Vert + \Vert \nabla_X f(\mathbf{X}^{k+1},\mathbf{Z}^{k},\mathbf{Z}^{k}) +   \bm{\mu}_1   \Vert   \\
		&\leq L_x \Vert \mathbf{X}^{k+1} - \mathbf{X}^{k,m_1} \Vert + M \Vert {\mathbf{W}}^{k+1} - {\mathbf{W}}^{k} \Vert   \\
		&\leq  L_x \Vert \mathbf{X}^{k+1} - \mathbf{X}^{k,m_1} \Vert + \frac{M}{\alpha} \Vert \bar{\mathbf{W}}^{k+1} - \mathbf{W}^{k+1} \Vert +
		M \Vert \bar{\mathbf{W}}^{k} - \mathbf{W}^{k}    \Vert.
	\end{array}
	\label{30}
\end{equation}
Repeating the above procedure for $\mathbf{Y}$ and $\mathbf{Z}$, we have
\begin{equation}
	\begin{array}{rl}
		\Vert  \nabla_\mathbf{Y} f(\mathbf{W}^{k+1}) \!+ \!  \bm{\mu}_2 \Vert \! \! \! \! \!   & \! \leq \!  L_y \Vert \mathbf{Y}^{k+1} - \mathbf{Y}^{k,m_2} \Vert \!  + \! \frac{M}{\alpha} \Vert \bar{\mathbf{W}}^{k+1} \!-\! \mathbf{W}^{k+1} \Vert    
		\!+\! M \Vert \bar{\mathbf{W}}^{k} \!- \!\mathbf{W}^{k}  \Vert,   \\
		\Vert  \nabla_\mathbf{Z} f(\mathbf{W}^{k+1})  \!+ \! \bm{\mu}_3 \Vert \! \! \! \! \! &\leq L_z \Vert \mathbf{Z}^{k+1} \!-\! \mathbf{Z}^{k,m_3} \Vert \!+\! \frac{M}{\alpha} \Vert \bar{\mathbf{W}}^{k+1} \!- \!\mathbf{W}^{k+1} \Vert  
		\! + \! M \Vert \bar{\mathbf{W}}^{k} \! - \!\mathbf{W}^{k}  \Vert.      \\  
	\end{array}
	\label{31}
\end{equation}
According to inequalities  $(\ref{27}) $, $(\ref{28}) $, $ (\ref{30})$,  and  $ (\ref{31}) $, we obtain
\begin{equation}
	\begin{array}{rl}
		\Vert  \nabla_\mathbf{X} f(\mathbf{W}^{k+1}) \! + \!  \bm{\mu}_1  \Vert & \! \! \! \! \leq  L_x \! \Vert \mathbf{X}^{k+1} \! - \! \mathbf{X}^{k,m_1} \Vert  \!+ \! \frac{M}{\alpha} \Vert \bar{\mathbf{W}}^{k+1}  \!- \!\mathbf{W}^{k+1} \Vert    
		\!+ \! M \Vert \bar{\mathbf{W}}^{k} \! - \! \mathbf{W}^{k}    \Vert,  \\
		
		\Vert  \nabla_\mathbf{Y} f(\mathbf{W}^{k+1})  \!+ \!  \bm{\mu}_2 \Vert & \! \! \! \! \leq  L_y \! \Vert \mathbf{Y}^{k+1} \!-\! \mathbf{Y}^{k,m_2} \Vert \! + \! \frac{M}{\alpha} \Vert \bar{\mathbf{W}}^{k+1} \!-  \!\mathbf{W}^{k+1} \Vert   
		\! + \! M \Vert \bar{\mathbf{W}}^{k} \!- \! \mathbf{W}^{k}  \Vert,  \\
		
		\Vert  \nabla_\mathbf{Z} f(\mathbf{W}^{k+1}) \! + \!  \bm{\mu}_3 \Vert & \! \! \! \! \!  \leq  L_z \! \Vert \mathbf{Z}^{k+1} - \mathbf{Z}^{k,m_3} \Vert \! + \! \frac{M}{\alpha} \Vert \bar{\mathbf{W}}^{k+1} \!- \!\mathbf{W}^{k+1} \Vert       
		\!+\! M \Vert \bar{\mathbf{W}}^{k}\! - \!\mathbf{W}^{k}  \Vert .  \\
	\end{array}
	\label{32}
\end{equation}
Next, we consider the upper bound of $\Vert \mathbf{X}^{k+1} - \mathbf{X}^{k,m_1} \Vert $, $\Vert \mathbf{Y}^{k+1} - \mathbf{Y}^{k,m_2} \Vert $, and $\Vert \mathbf{Z}^{k+1} - \mathbf{Z}^{k,m_3} \Vert$. Suppose that $\mathbf{W}^{k+1}$ is calculated by Step 4, 6 and 8, we obtain
\begin{equation}
	\begin{array}{rl}
		& \Vert \mathbf{X}^{k+1} \!  - \!  \mathbf{X}^{k,m_1} \Vert^2 \!\! \leq \frac{2}{\epsilon} (F(\mathbf{X}^{k,m_1},\bar{\mathbf{Y}}^{k},\bar{\mathbf{Z}}^{k})-F(\mathbf{X}^{k+1},\bar{\mathbf{Y}}^{k},\bar{\mathbf{Z}}^{k}))   \\
		
		& \leq  \frac{2}{\epsilon} (F(\mathbf{X}^{k},\bar{\mathbf{Y}}^{k},\bar{\mathbf{Z}}^{k})-F(\mathbf{X}^{k+1},\bar{\mathbf{Y}}^{k},\bar{\mathbf{Z}}^{k}))   \\
		& =  \frac{2}{\epsilon} (f(\mathbf{X}^{k},\bar{\mathbf{Y}}^{k},\bar{\mathbf{Z}}^{k})-f(\mathbf{X}^{k+1},\bar{\mathbf{Y}}^{k},\bar{\mathbf{Z}}^{k}))+\frac{2}{\epsilon}(r_1(\mathbf{X}^{k})-r_1(\mathbf{X}^{k+1}))  \\
		& =  \frac{2}{\epsilon}(F(\mathbf{X}^{k},{\mathbf{Y}}^{k},\mathbf{Z}^{k})-F(\mathbf{X}^{k+1},({\mathbf{Y}}^{k},{\mathbf{Z}}^{k})))+ \frac{2}{\epsilon} (f(\mathbf{X}^{k+1},({\mathbf{Y}}^{k},{\mathbf{Z}}^{k}))  \\
		&- \frac{2}{\epsilon}f(\mathbf{X}^{k},{\mathbf{Y}}^{k},{\mathbf{Z}}^{k})    
		+ \frac{2}{\epsilon} (f(\mathbf{X}^{k},\bar{\mathbf{Y}}^{k},\bar{\mathbf{Z}}^{k})-f(\mathbf{X}^{k+1},\bar{\mathbf{Y}}^{k},\bar{\mathbf{Z}}^{k}))     \\
		&\leq   \frac{2}{\epsilon}(F(\mathbf{X}^{k},({\mathbf{Y}}^{k},{\mathbf{Z}}^{k}))-F(\mathbf{X}^{k+1},{\mathbf{Y}}^{k},{\mathbf{Z}}^{k}) + l_x^+ \Vert \mathbf{X}^{k+1}-\mathbf{X}^{k}   \Vert^2   \\
		&+\langle \mathbf{X}^{k+1}-\mathbf{X}^{k}, \nabla_\mathbf{X} f(\mathbf{X}^{k},{\mathbf{Y}}^{k},{\mathbf{Z}}^{k}) -\nabla_\mathbf{X} f(\mathbf{X}^{k+1},\bar{\mathbf{Y}}^{k},\bar{\mathbf{Z}}^{k})   \rangle ) \\
		& \leq  \frac{2}{\epsilon}(F(\mathbf{X}^{k},{\mathbf{Y}}^{k},{\mathbf{Z}}^{k})-F(\mathbf{X}^{k+1},{\mathbf{Y}}^{k},{\mathbf{Z}}^{k}) + (\frac{1}{2}+l_x^+) \Vert \mathbf{X}^{k+1}-\mathbf{X}^{k}  \Vert^2  \\
		&+\frac{1}{2} \Vert \nabla_\mathbf{X} f(\mathbf{X}^{k},{\mathbf{Y}}^{k},{\mathbf{Z}}^{k}) -\nabla_\mathbf{X} f(\mathbf{X}^{k+1},\bar{\mathbf{Y}}^{k},\bar{\mathbf{Z}}^{k})    \Vert^2 )  \\
		&\leq \frac{2}{\epsilon}(F(\mathbf{X}^{k},{\mathbf{Y}}^{k},{\mathbf{Z}}^{k})-F(\mathbf{X}^{k+1},({\mathbf{Y}}^{k},{\mathbf{Z}}^{k})) + (\frac{1}{2}) \Vert \mathbf{X}^{k+1}-\mathbf{X}^{k}  \Vert^2    \\
		&+l_x^+\Vert \mathbf{X}^{k+1}-\mathbf{X}^{k}  \Vert^2  + \frac{M^2}{2} \Vert (\mathbf{X}^{k+1}-\mathbf{X}^k,\mathbf{Y}^k-\bar{\mathbf{Y}}^k,\mathbf{Z}^k-\bar{\mathbf{Z}}^k)   \Vert^2 )   \\
		&=   \frac{2}{\epsilon}(F(\mathbf{X}^{k},{\mathbf{Y}}^{k},{\mathbf{Z}}^{k})-F(\mathbf{X}^{k+1},{\mathbf{Y}}^{k},{\mathbf{Z}}^{k}) + (\frac{1+M^2}{2}) \Vert \mathbf{X}^{k+1}-\mathbf{X}^{k}  \Vert^2    \\
		&+l_x^+ \Vert \mathbf{X}^{k+1}-\mathbf{X}^{k}  \Vert^2 + \frac{M^2}{2} \Vert\mathbf{Y}^k-\bar{\mathbf{Y}}^k  \Vert^2 + \frac{M^2}{2}
		\Vert \mathbf{Z}^k-\bar{\mathbf{Z}}^k  \Vert^2 ).  \\ 
	\end{array}
	\label{33}
\end{equation}
The first two inequalities follow from a basic descent property of proximal gradient steps. The third inequality holds by Lemma $\ref{dec}  $. The fourth inequality utilizes $ \langle \mathbf{a},\mathbf{b} \rangle \leq \frac{1}{2} \Vert \mathbf{a}\Vert^2+\frac{1}{2} \Vert \mathbf{b}\Vert^2$. The fifth inequality is the result of assumption $\ref{a_4} $(iv).

Suppose that $\mathbf{W}^{k+1}$ is calculated by Step 11, 13, 15, we have
\begin{equation}
	\begin{array}{rl}
		\Vert \mathbf{X}^{k+1} - \mathbf{X}^{k,m_1} \Vert^2 
		&\leq \frac{2}{\epsilon} (F(\mathbf{X}^{k,m_1},{\mathbf{Y}}^{k},{\mathbf{Z}}^{k})-F(\mathbf{X}^{k+1},{\mathbf{Y}}^{k},{\mathbf{Z}}^{k}))    \\
		& \leq  \frac{2}{\epsilon} (F(\mathbf{X}^{k},{\mathbf{Y}}^{k},{\mathbf{Z}}^{k})-F(\mathbf{X}^{k+1},{\mathbf{Y}}^{k},{\mathbf{Z}}^{k})) .\\
	\end{array}
	\label{34}
\end{equation}
Combining (\ref{33}) and (\ref{34}) yields
\begin{equation}
	\begin{array}{rl}
		\Vert \mathbf{X}^{k+1} - \mathbf{X}^{k,m_1} \Vert^2 
		&\leq \frac{2}{\epsilon}(F(\mathbf{X}^{k},{\mathbf{Y}}^{k},{\mathbf{Z}}^{k})-F(\mathbf{X}^{k+1},{\mathbf{Y}}^{k},{\mathbf{Z}}^{k})    \\
		& + \frac{M^2}{2} \Vert \mathbf{W}^k-\bar{\mathbf{W}}^k  \Vert^2 + (\frac{1+M^2}{2}+l_x^+) \Vert \mathbf{W}^{k+1}-\mathbf{W}^{k}  \Vert^2) .    \\
	\end{array}
	\label{35}
\end{equation}
Then the following inequalities for $\mathbf{Y}$ and $\mathbf{Z}$ hold 
\begin{equation}
	\begin{array}{rl}
		\Vert \mathbf{Y}^{k+1} - \mathbf{Y}^{k,m_2} \Vert^2 
		&\leq   \frac{2}{\epsilon} ( F(\mathbf{X}^{k+1},{\mathbf{Y}}^{k},{\mathbf{Z}}^{k})-F(\mathbf{X}^{k+1},{\mathbf{Y}}^{k+1},{\mathbf{Z}}^{k})  \\
		& +  (\frac{1+M^2}{2}+l_y^+) \Vert \mathbf{W}^{k+1}-\mathbf{W}^{k}  \Vert^2 +   \frac{M^2}{2} \Vert \mathbf{W}^{k+1}-\bar{\mathbf{W}}^{k+1}  \Vert^2   \\
		&+    \frac{M^2}{2} \Vert \mathbf{W}^k-\bar{\mathbf{W}}^k \Vert^2 ) , \\
		
		\Vert \mathbf{Z}^{k+1} - \mathbf{Z}^{k,m_3} \Vert^2 
		&\leq \frac{2}{\epsilon} (F(\mathbf{X}^{k+1},\mathbf{Y}^{k+1},\mathbf{Z}^k)-F(\mathbf{X}^{k+1},\mathbf{Y}^{k+1},\mathbf{Z}^{k+1})    \\
		&+ (\frac{1+M^2}{2}+l_z^+) \Vert \mathbf{W}^{k+1}-\mathbf{W}^{k}  \Vert^2 + \frac{M^2}{2} \Vert \mathbf{W}^{k+1}-\bar{\mathbf{W}}^{k+1}  \Vert^2 ).  \\
	\end{array}
	\label{36}
\end{equation}
Combining (\ref{35}) and (\ref{36}), we can deduce that
\vspace{-3mm}
\begin{equation*}
	\begin{array}{rl}
		\Vert \mathbf{W}^{k+1} - \mathbf{W}^{k,{m}} \Vert^2 &= \Vert \mathbf{X}^{k+1} - \mathbf{X}^{k,m_1} \Vert^2 +   \Vert \mathbf{Y}^{k+1} - \mathbf{Y}^{k,m_2} \Vert^2+\Vert \mathbf{Z}^{k+1} - \mathbf{Z}^{k,m_3} \Vert^2  \\
		&\leq  \frac{2}{\epsilon}(F(\mathbf{X}^{k},{\mathbf{Y}}^{k},{\mathbf{Z}}^{k})-F(\mathbf{X}^{k+1},\mathbf{Y}^{k+1},\mathbf{Z}^{k+1})  \\
		&+ (\frac{3(1+M^2)}{2}+ l_x^+  + l_y^+ +l_z^+) \Vert \mathbf{W}^{k+1}-\mathbf{W}^{k}  \Vert^2   \\
		&+ M^2( \Vert \mathbf{W}^k-\bar{\mathbf{W}}^k  \Vert^2+ \Vert \mathbf{W}^{k+1}-\bar{\mathbf{W}}^{k+1}  \Vert^2)).
	\end{array}
\end{equation*}
Applying Lemma $\ref{l_2}$  $\textbf{(ii)}$, we have $ \lim\limits_{k \to +\infty} \Vert \mathbf{W}^{k+1} -\mathbf{W}^{k}  \Vert = 0 $ and $ \lim\limits_{k \to +\infty} \Vert  \bar{\mathbf{W}}^{k} -\mathbf{W}^{k}  \Vert = 0 $. Therefore, it can be inferred that $ \lim\limits_{k \to +\infty} \Vert \mathbf{W}^{k+1} - \mathbf{W}^{k,{m}} \Vert=0$, which indicates that the limit of $ \Vert \mathbf{X}^{k+1} - \mathbf{X}^{k,m_1} \Vert,\Vert \mathbf{Y}^{k+1} - \mathbf{Y}^{k,m_2} \Vert $ and $\Vert \mathbf{Z}^{k+1} - \mathbf{Z}^{k,m_3} \Vert $ all converge to 0. Finally, we have $ \lim\limits_{k \to +\infty} \Vert  \nabla_X f(\mathbf{W}^{k+1}) +  \bm{\mu}_1  \Vert,\lim\limits_{k \to +\infty} \Vert  \nabla_Y f(\mathbf{W}^{k+1}) +  \bm{\mu}_2 \Vert$ and $\lim\limits_{k \to +\infty} \Vert  \nabla_Z f(\mathbf{W}^{k+1}) +  \bm{\mu}_3 \Vert $  all equal to 0 by (\ref{32}).

From [\cite{rockafellar1998variational} , Proposition 10.5  10.10(subgradients of Lipschitzian sums)], we obtain the following equation
\begin{equation*}
	\begin{array}{lll}
		\nabla_{\mathbf{W}} f(\mathbf{W}^{k+1}) \! \!  + \! \! \bm{\mu} \!= \! (\nabla_{ \mathbf{X}} f(\mathbf{W}^{k+1}), \! \nabla_{ \mathbf{Y}} f(\mathbf{W}^{k+1}), \! \nabla_{ \mathbf{Z}} f(\mathbf{W}^{k+1})) \! \! + \!  \!
		( \bm{\mu}_1 , \bm{\mu}_2 , \bm{\mu}_3) \! \in  \! \partial_{\mathbf{W}} F(\mathbf{W}^{k+1}) .
	\end{array}
	\label{38}
\end{equation*}
Hence, it can be concluded that  $\lim\limits_{k \to +\infty} \Vert \nabla_{\mathbf{W}} f(\mathbf{W}^{k+1}) +  \bm{\mu} \Vert =0.$

Based on the assumption and constraints, we have the sequence $ \{\mathbf{W}^k\}_{k \in N}  $ is bounded. Let $\mathbf{W}^*$  be  an accumulation point of $ \{\mathbf{W}^k\}_{k \in N}  $, there exists a subsequence $ \{ \mathbf{W}^{k_p } \} $ of $ \{\mathbf{W}^k\}_{k \in N}  $ which converges to $\mathbf{W}^*$. It follows that $\Vert \nabla_{\mathbf{W}} f(\mathbf{W}^{*}) +  \bm{\mu}^* \Vert =  \lim\limits_{k \to +\infty} \Vert \nabla_{\mathbf{W}} f(\mathbf{W}^{k+1}) +  \bm{\mu} \Vert =0$, which implies $\mathbf{0} \in \partial F(\mathbf{W}^*). $ Therefore, one can deduce that $\mathbf{W}^*$ is a stationary point of the optimization problem (\ref{gene_Cp_MODEL}). This concludes the proof.    \hfill  $\square   $ \\

\noindent
\textbf{Funding} This work was supported by the National Natural Science Foundation of China (Grant No.12671362) \\

\noindent
\textbf{Data availability} The data that support the findings of this study are openly available at the following URL/DOI:\url{https://github.com/V-Geler/Conv2dPA/tree/main/Simulator}.

\section*{Declarations}

\noindent
\textbf{Conflict of interest} The authors declare that they have no conflict of interest.

\end{appendices}

\bibliography{reference}

\end{document}